\documentclass[10pt]{amsart}
\usepackage{amsmath,amssymb}
\usepackage[top=1in,bottom=1in,left=1in,right=1in]{geometry}
\usepackage{amsthm}
\usepackage{newlfont}
\usepackage{booktabs}
\usepackage{caption}
 \usepackage{subcaption}
\usepackage{hyperref}
\usepackage{tikz,float}
\usepackage{graphicx}
\usepackage{array}
\usepackage{tikz}
\usetikzlibrary{arrows,shapes,positioning}
\usetikzlibrary{backgrounds,patterns,calc}
\usetikzlibrary{intersections}
\usepackage{pgfplots}
\usepackage{thmtools}

\usepackage{pstricks}
\usepackage[T1]{fontenc}
\usepackage{relsize}

\usepackage{amsfonts}
\usepackage{mathrsfs}
\usepackage{amsmath}
\usepackage{amscd}
\usepackage{amsthm}

\usepackage{dsfont}
\usepackage[all,cmtip]{xy}
\usepackage{pgflibraryarrows}
\usepackage{pgflibrarysnakes}

\usepackage[utf8]{inputenc}

\usepackage{mathtools}
\usepackage{accents}
\usepackage{hyperref}

\newcommand\ub[1]{%
  \underaccent{\bar}{#1}}
\newcommand\ul[1]{%
  \underline{#1}}
\newenvironment{tikz1}{\vspace*{0pt}\begin{center}\begin{tikzpicture}}
{\
\end{tikzpicture}\end{center}\vspace*{-5pt}

}
\numberwithin{equation}{section}
\newtheorem{thm}{Theorem}
\newtheorem{definition}{Definition}[section]
\newtheorem{question}{Question}

\newtheorem*{question*}{Question}

\newtheorem{lem}[definition]{Lemma}
\newtheorem{pro}[definition]{Proposition}
\newtheorem{cor}[definition]{Corollary}
\newtheorem{claim}[definition]{Claim}

\newtheorem{rem}[definition]{Remark}

\newtheorem*{conj*}{Conjecture}

\newtheorem*{pres*}{Presumption}

\newcommand{\bl}{\boldsymbol}

\newcommand{\rin}{\mathrm{rin}}

\newcommand{\pic}{\mathrm{Pic}}
\newcommand{\amp}{\mathrm{Amp}}

\newcommand{\dnef}{\mathfrak{D} }
\newcommand{\el}{\mathfrak{L} }

\newcommand{\nef}{\mathrm{Nef}}

\newcommand{\coker}{\mathrm{coker}}
\newcommand{\wid}{\mathrm{wid}}

\begin{document}

\title{Some Finiteness Results on Oda's Question for Pairs of Line Bundles}
\author{He Xin}
\thanks{xinhe@qymail.bhu.edu.cn}
\date{}
\maketitle

\section*{Abstract}

For a projective toric variety $X$ and nef line bundles $L_1, L_2$ on it we prove some finiteness results related to Oda's question, which concerns the surjectivity of the multiplication map $\Phi_{L_1, L_2}: H^0(X, L_1)\otimes H^0(X, L_2) \rightarrow H^0(X, L_1\otimes L_2)$. When  $X$ is general in the sense that there are no paralleling edges on the polytopes associated to ample line bundles on it we show $\dim\coker\:\Phi_{L_1, L_2}$ can be uniformly bounded. When $L_1$ is an ample line bundle lying outside a certain finite subset of $\amp(X)$ for a smooth toric threefold $X$, we show $\Phi_{L_1, L_2}$ is surjective for any nef line bundle $L_2$. Explicit bounds for pairs that might not fulfill the surjectivity property are also obtained in various settings.
In the proofs we make use of coverings and quasi-coverings of polytopes together with a filtered structure on $\nef(X)$  gained from the restrictions of nef line bundles to the invariant curves.

\section{Introduction}

Let $X$ be a smooth projective toric variety, $L_1$ be an ample and $L_2$  a nef line bundle on it. In \cite
{odaproblems} Oda asks whether the following multiplication map is surjective 
\begin{equation}\Phi_{L_1, L_2}: H^0(X, L_1)\otimes H^0(X, L_2)\rightarrow H^0(X, L_1\otimes L_2).\label{sp}\end{equation}

Recall that one can associate a polytope $P_L$ to a nef line bundle $L$ on  $X$ and the lattice points in $P_L$ just form a basis of $H^0(X, L)$. In terms of polytopes, the surjectivity of (\ref{sp}) is equivalent to the surjectivity of the following map defined by Minkowski sum 
\begin{equation}(P_{L_1}\cap M) \times (P_{L_2}\cap M)\rightarrow P_{L_1\otimes L_2}\cap M\label{ssp},\end{equation}
where $M$ is the group of characters of the torus acting on $X$.

 Obviously, when $L_1, L_2$ take values in $\{L^k\}_{k\geq1}$ for some given ample line bundle $L$, the surjectivity of (\ref{sp}) just implies $L$ is projectively normal.  Conversely, if one could prove $\Phi_{L, L}$ is surjective for \emph{any} ample line bundle on \emph{any} smooth toric variety $X$ \cite[Remark 4.3]{hmp2010}, then the surjectivity of (\ref{sp}) for pairs of ample line bundles follows readily. 
 
Concerning the projective normality of very ample line bundles on projective varieties, the first positive result can be traced back to Castelnuovo
and Mumford \cite{mumford1969varieties}, where it is proved a line bundle $L$ on a curve of genus $g$  is projectively normal when $\deg L\geq 2g+1$. For projective normality on  varieties of higher dimensions, a complete survey can be found in \cite[1.8D]{lazarsfeld2004p1}. When it comes to toric varieties, projective normality and higher syzygies are investigated by Ewald et al. \cite{ewald1991ampleness} and Hering et al. \cite{hering2006syzygies} respectively. Similar as other varieties, the line bundles that are verified to be projectively normal are required to be powers of some other ample line bundles.  In addition, the dimension of the base variety also appears in these powers for bounding the positivity of the ample line bundle.

 The map (\ref{sp}) is obviously surjective for projective spaces, which are among the simplest toric varieties. As a generalization of projective spaces in another vein flag varieties are also proved to be projectively normal \cite{ramanan1985projective}, in which the technique of diagonal splitting of Frobenius plays a key role. Unfortunately, this method cannot be applied to toric varieties as they are not diagonally split in general \cite{payne2009frobenius}. Probably due to the fact that toric varieties are much more varied, a proof of the surjectivity of (\ref{sp}) in its final form is more difficult. As far, Oda's question
is still widely open, see \cite{gubeladze2023normal} for a recent survey.   

  On the positive side, the equivalent reformulation (\ref{ssp}) makes Oda's question more approachable if treated from an intuitive point of view. The surjectivity of (\ref{sp}) is proved to hold for smooth toric surfaces in \cite{ramanan1985projective} and general toric surfaces not necessarily smooth in \cite{haase2008lattice}. The methods used in these works relies heavily on the dimension of the base, which makes it difficult to generalize. A more recent major breakthrough is achieved in \cite{gubeladze2012convex}, where an ample line bundle on a toric variety of dimension $d$ is proved to be projectively normal provided the lattice length of the associated polytope is no less than $4d(d+1)$.  This bound is improved to $2d(d+1)$ in \cite{haase2017convex}.

In spite of the importance of the accomplishment made in \cite{gubeladze2012convex}, Oda's question is still unanswered fully. Firstly, in  \cite{gubeladze2012convex} only the special case when $L_1, L_2$ are powers of a third line bundle is considered. Towards this direction Haase and Hofmann make progress in \cite{haase2017convex} by requiring that one of the polytopes has each of its edge at least $d$ times larger than the corresponding edge of the other. Secondly, for a given toric variety there might be infinitely many ample line bundles out of the scope of these works loc. cit. , as the lattice lengths of their associated polytopes might not reach the lower bounds aforementioned. With these questions in mind, we prove the following results.

\begin{thm}\label{t1}  If $X$ is general in the sense that no two edges on the polytope associated to an ample line bundle on it are parallel, then $\dim\coker\:\Phi_{L_1, L_2}$ can be uniformly bounded for any pair of nef line bundles $L_1, L_2$.
\end{thm}

\begin{thm}\label{t2} For a smooth projective toric threefold $X$ there exists a finite subset of $\amp(X)$ such that the morphism (\ref{sp}) is  surjective for any ample line bundle $L_1$ not in that subset  and any nef line bundle $L_2$. In addition, $\dim\coker\:\Phi_{L_1, L_2}$ can be uniformly bounded. 

\end{thm}

\subsection*{Strategy of the Proof and Outline of the Paper}\hfill

First of all, instead of (\ref{ssp}) in this work we will mainly investigate the following map
\begin{equation}\Psi_{L_1, L_2}: P_{L_1}\times (P_{L_2}\cap M)\rightarrow P_{L_1\otimes L_2}\label{os},\end{equation}
where $L_1$ is ample and $L_2$ is nef as before. It's evident the surjectivity of (\ref{os}) implies that of (\ref{ssp}). Our starting point is the following easy observation. For each vertex  of $P_{L_1\otimes L_2}$ we have a lattice translation of $P_{L_1}$ lying inside $P_{L_1\otimes L_2}$ such that  this vertex and its corresponding vertex of $P_{L_1}$ coincide. If the difference of the lattice lengths of the edges of $P_{L_1}$ and $P_{L_1\otimes L_2}$ are small compared with their sizes, the points of $P_{L_1\otimes L_2}$ not in a translation of $P_{L_1}$ as above could only be found near the boundary of $P_{L_1\otimes L_2}$. By applying Carath\'eodory's theorem we can prove (see Lemma \ref{lem1}) in this setting all these translations of $P_{L_1}$ obtained via fitting the corresponding vertices will suffice to cover $P_{L_1\otimes L_2}$. Following this line we can show (Proposition \ref{p1}) the map (\ref{sp}) is surjective as long as $L_1$ is sufficiently ample but not necessarily a power of another line bundle.

Another essential ingredient in our proofs lies in the fact that the nef cone of a projective toric variety is a rational polyhedral cone and it is endowed with a natural partially ordered filtration consisting of \emph{linear subsets} (see Definition \ref{lss}). In the proof of Theorem \ref{t1} the polytopes associated to those sufficiently ample line bundles can be used as fillers to `cover' polytopes associated to other nef line bundles. Albeit such coverings are not perfect the space left out can be uniformly bounded, and we will call them \emph{quasi-coverings} (see Claim \ref{qcl3}). The first main result just follow from our proof that the polytopes associated to all but finitely many nef line bundles on a general toric variety admit such quasi-coverings.

In the case of smooth projective toric threefolds we first show for a given ample line bundle the dimension of the cokernel of the map (\ref{sp}) is bounded.
The proof of our second main result is built on Proposition \ref{fe}, which says that for a given nef line bundle
there always exists a finite subset of $\amp(X)$ such that the map (\ref{sp}) is surjective for any ample line bundle not in this subset and the nef line bundle. In the proof of the latter we make use of  Theorem \ref{sfhn}, which is nothing but the surjectivity of (\ref{os}) in dimension 2.

After proving our main results we elaborate the dimension lowering process, which already appears in the proof of Lemma \ref{pfr}. With this decisive input, we can make an explicit estimation (see Proposition \ref{loqr} and Proposition \ref{lopr}) as complement to the finiteness results.

Theorem \ref{sfhn} is proved by induction on the Picard number of  the base toric surface in the final section. The most delicate part of the proof occurs while choosing a suitable subset of $P_{L_1^{-1}\otimes L_2}\cap M$ (we follow the notations in the proof) such that  the union of the translations of $P_{L_1}$
by the vectors in that subset can cover the region missed in the inductive step. In two situations, i.e. (I1) and (I2), where such choices should be made, the numbers of such translation vectors turn out to be one and four. In the first case, although one translation of $P_{L_1}$  
would sustain our induction, we need to select a set of translation vectors first then prove one of them will suffice. A case-by-case argument for determining admissible types of adjacent translation vectors in this situation seems unavoidable and we make frequent use of the condition that the underlying toric surface is smooth.




\subsection*{Conventions and Notations}\hfill

$X=X(\Sigma)$: a toric variety $X$ defined by $\Sigma$;

$\Sigma(k)$: the set of cones of $\Sigma$ with dimension $k$;

$|\Sigma|$: the support of $\Sigma$;

$\Sigma_L$: the normal fan of the polytope $P_L$ associated to a nef line bundle $L$;

$\rho\in\Sigma_X(1)$: a one-dimensional cone of $\Sigma_X$ and the primitive vector paralleling with this cone
 
$D_{\rho}$: the invariant divisor corresponding to $\rho\in\Sigma_X(1)$;

$F_{\alpha}(P_L)$: the facet of $P_L$ corresponding to a cone $\alpha$ of $\Sigma_L$;

$H_{\alpha}(P_L)$: the minimal affine subspace which contains $F_{\alpha}(P_L)$;


$v_{\sigma}(P_L)$: the vertex of $P_L$ corresponding to a full-dimension cone $\sigma$ of $\Sigma_L$;

$e_{\tau}(P_L)$: the edge of $P_L$ corresponding to a codimension one cone $\tau$ of $\Sigma_L$;

$u_{\tau}$: a primitive vector in $M$ paralleling with the edge $e_{\tau}(P_L)$;

$u_{\tau}(\sigma)$: the primitive vector in $M$ paralleling with $e_{\tau}(P_L)$ such that $v_{\sigma}(P_L)+\epsilon u_{\tau}(\sigma)\in e_{\tau}(P_L)$ for all $\epsilon>0$ sufficiently small;


$AB$: the line passing through $A$ and $B$;

$\overline{AB}$: the segment with endpoints $A$ and $B$;

$l_{\bl{v}}(p)$: the line passing through $p$ and paralleling with $\bl{v}$;

$r_{\bl{v}}(p)$: the ray with initial at $p$ and and paralleling with  $\bl{v}$;

$\rin(S)$: the relative interior of a point set $S$.

For convenience we will occasionally endow $M_{\mathds{R}}$ with an Euclidean norm $\|\cdot\|$. Following \cite{gubeladze2012convex}, the lattice length $\mathbb{E}(L)$ of a nef line bundle $L$ on $X$ is defined to be the minimum of the lattice lengths of the edges of $P_L$. 
 For two polytopes $P_1$ and $P_2$ with the same normal fan the number 
$\mathbb{E}(\tfrac{P_1}{P_2})$ is defined to be the minimum of the ratios of the (lattice) lengths of the corresponding edges of $P_1$ and $P_2$. When $P_1=P_{L_1}, P_2=P_{L_2}$ for two line bundles $L_1, L_2$ we will denote  $\mathbb{E}(\tfrac{P_{L_1}}{P_{L_2}})$ by $\mathbb{E}(\tfrac{L_1}{L_2})$.

\section{Preliminary Materials}

In this section we collect some standard facts as well as some less standard definitions and results that will be used in our proofs.
\subsection{Partial Orders on the Nef Cone}\hfill

For a projective toric variety $X$, recall that its nef cone is a rational strongly convex polyhedral cone in $\pic(X)_{\mathds{R}}=N^1(X)$ and each lattice point in the nef cone corresponds to a nef line bundle on $X$. In this work we will denote by $\nef(X)$ the semigroup of nef line bundles, $\amp(X)\subset\nef(X)$ the set of ample line bundles and $\nef(X)_{\mathds{Q}}$ the nef cone.  

The semigroup $\nef(X)$ is endowed a natural partial order `$\prec$'  given by
\[L\prec \tilde{L}\Leftrightarrow  \tilde{L}\in L+\nef(X).\]
More generally, for a convex subcone $\mathfrak{C}\subseteq\nef(X)_{\mathds{Q}}$, we can define a partial order $\prec_{\mathfrak{C}}$ on $\nef(X)_{\mathds{Q}}$ as follows
\begin{equation}\label{zzpx}L\prec_{\mathfrak{C}} \tilde{L}\Leftrightarrow  \tilde{L}\in L+\mathfrak{C}.\end{equation}
To define another partial order `$\prec_c$' on $\nef(X)$, recall that to any nef line bundle $L$
we can associate to it a convex polytope $P_L$, which is the convex hull of the following set of characters 
\[\{\chi\in M\:|\:H^0(X,L)_{\chi}\neq0\}.
\]
\noindent Information of the line bundle can be read from that of the 
associated polytope and vice versa. For instance, when $L$ is big $P_L$ has the highest possible dimension. When we translate $P_L$ in $M$ by a character $\chi$, the corresponding line bundle would be $L(-D_{\chi})$ with $D_{\chi}=\langle\chi,\rho\rangle D_{\rho}$, which is linearly equivalent to $L$.

Now we define for any $L, \tilde{L}\in\nef(X)$
\begin{center}
$L\prec_c \tilde{L}\Leftrightarrow$ each point of $P_{\tilde{L}}$ lies in some lattice translation of $P_L$ which is contained in $P_{\tilde{L}}$. 
\end{center}
On the set $\amp(X)$ we can define another partial order $\prec_o$ by $L\prec_o \tilde{L}$ iff $L\prec \tilde{L}$ and $\coker\:\Phi_{L, L^{-1}\otimes \tilde{L}}=0$. Then it is obvious that we have the following implication
\[L\prec_c \tilde{L}\Rightarrow L\prec_o \tilde{L}\Rightarrow L\prec \tilde{L}.
\]
and Oda's question can be reformulated as $L\prec \tilde{L}\Rightarrow L\prec_o \tilde{L}$.
\begin{rem} 

Unlike the partial order $\prec$, we do not whether the following implication is true \begin{equation}\label{tpcp}L\prec_c\tilde{L}\Rightarrow L\otimes L'\prec_c \tilde{L}\otimes L'.\end{equation} It is not know either whether (\ref{tpcp}) can obtained with one of or both of the relations $L\prec_c L\otimes L'$, $L\prec_c \tilde{L}\otimes L'$.
More generally, there are also implications without tensors unclear to us such as 
\begin{equation}\notag L_1\prec_c L_2, L_1\prec_c L_3, L_2\prec L_3\Rightarrow L_2\prec_c L_3\end{equation}
for $L_1, L_2, L_3\in\amp(X)$. What is even more obscure is when we consider implications as above with some of `$\prec$' on the left or/and `$\prec_c$' on the right replaced by $\prec_o$.

\end{rem}


The following lemma is the starting point of all our results in this work, the proof of it is apparent and we omit it.
\begin{lem}
Let $L, \tilde{L}$ be two ample line bundles such that $L\prec_c \tilde{L}$, then the map (\ref{sp}) is surjective for the pair $(L, L^{-1}\otimes \tilde{L})$.
\end{lem}
Now it is clear to see Oda's question has an affirmative answer if for any pair of ample line bundles $L, \tilde{L},$ we can prove 
\begin{equation}\label{o1}L\prec \tilde{L}\Rightarrow L\prec_c \tilde{L}.
\end{equation}
However, this implication is not always true as one notices by taking $L=\mathcal{O}_{\mathbb{P}^2}(1)$.

The following observations for a pair of nef line bundles $\ul{L}\prec L$ will be used later. 

Let $\Sigma_{\ul{L}}$ (resp. $\Sigma_L$) be the normal fan of $P_{\ul{L}}$ (resp. $P_L$) and $n_{\ul{L}}=\dim\Sigma_{\ul{L}}$, $n_{L}=\dim\Sigma_{L}$. Note that for any $t>0$, $\Sigma_L$ is also the normal fan of $P_{\ul{L}}+tP_L$. 
Then for any $\ub{\sigma}\in\Sigma_{\ul{L}}(n_{\ul{L}})$ we can define the following subset $\Sigma_L(n_L)_{\ub{\sigma}}$ of $\Sigma_L(n_L)$
\begin{equation}\label{pv}\sigma\in\Sigma_L(n_L)_{\ub{\sigma}}\Longleftrightarrow \lim_{t\rightarrow 0^+}v_{\sigma}(P_{\ul{L}}+tP_L)=v_{\ub{\sigma}}(P_{\ul{L}}),
\end{equation}
from which one sees easily \[\Sigma_L(n_L)=\bigsqcup_{\sigma\in\Sigma_{\ul{L}}(n_{\ul{L}})}\Sigma_L(n_L)_{\ub{\sigma}}.\]



\subsection{Linear Subsets of the Nef Cone}\hfill

\begin{definition}\label{lss}Let $\{C_{\tau}\}_{\tau\in I}$ be the set of invariant curves on $X$, then a  linear subset of $\nef(X)_{\mathds{Q}}$ is a subset of the following form
\[\el^J_{\bl{b}} =\{L\in\nef(X)_{\mathds{Q}}\:|\:L.C_{\tau}=b_{\tau}\;\text{if}\; \tau\in J \;\text{and}\; L.C_{\tau}\geq b_{\tau} \;\text{if}\; \tau\notin J\},
\]
where $\bl{b}=(b_{\tau})_{\tau\in I}\in\mathds{Z}_{\geq0}^I$ and $J\subseteq I$. 
A nonempty linear subset $\mathfrak{L}^J_{\bl{b}}$ or its index $(J, \bl{b})$ is said to be reduced if  

i) there exists $L_0\in\nef(X)$ such that $b_{\tau}=L_0.C_{\tau}$ for all $\tau\in I$;

ii) for each $\tau\in I\backslash J$ the image of the following map is infinite \begin{equation}
\label{ii}\mathfrak{L}^J_{\bl{b}}\cap\nef (X)\rightarrow \mathds{Z}, \;\;L\rightarrow L.C_{\tau}.
\end{equation}
\end{definition}


\begin{lem}\label{dd}Let $\mathfrak{L}^J_{\bl{b}}$ be a nonempty linear subset, then

1) $\mathfrak{L}^J_{\bl{b}}$ is reduced iff it can be written in the form $L+\mathfrak{L}^J_{\bl{0}}$, where $L$ is a nef line bundle and 
$\mathfrak{L}^J_{\bl{0}}$ is a reduced subcone of $\nef(X)$;

2) there exists finitely many reduced linear subsets $\mathfrak{L}^{J_k}_{\bl{b}_k}$ such that
\[\mathfrak{L}^J_{\bl{b}}\cap\nef(X)=\bigcup_k\mathfrak{L}^{J_k}_{\bl{b}_k}\cap\nef(X).
\]


\end{lem}
\begin{proof} The first part follows from definition and the second part from \cite[Theorem 2.12]{bg2009}.
\end{proof}

The following lemma will be used in the proofs of  our main results. Note that it is still valid if we replace $\nef(X)$ by $\amp(X)$.

\begin{lem}\label{pfr} 
Let $\mathscr{P}$ be a property defined on $\nef(X)$. Suppose for each reduced linear subset $\mathfrak{L}^J_{\bl{b}}\subseteq \nef(X)_{\mathds{Q}}$  with $\bl{b}=(\bl{b}_J, \bl{b}_{I\backslash J})$  we can find  $\bl{b}'$ with $\bl{b}'=(\bl{b}_J, \bl{b}'_{I\backslash J})$, $\bl{b}'_{I\backslash J}>\bl{b}_{I\backslash J}$ such that   the property $\mathscr{P}$ is true for all line bundles in $\mathfrak{L}^J_{\bl{b}'}$, then $\mathscr{P}$ is true for all but finitely many nef line bundles.
\end{lem}

\begin{proof} To prove the lemma, we note that
\begin{equation}\label{ddp}\mathfrak{L}^J_{\bl{b}}\backslash \mathfrak{L}^J_{\bl{b}'}\subseteq\bigcup_{\tau\in I\backslash J}\bigcup_{b_{\tau}\leq b''_{\tau}\leq b'_{\tau}}\mathfrak{L}^{J\cup\{\tau\}}_{(\bl{b}_J, b''_{\tau}, \bl{b}_{I\backslash (J\cup\{\tau\})})}.
\end{equation}
By Lemma \ref{dd} each linear subset on the right side can be written as a union of finitely many reduced ones. On the other hand, again by Lemma \ref{dd} any reduced linear subset is a translation of a reduced subcone of $\nef(X)$, hence the dimension of each of these reduced linear subsets mentioned above is strictly smaller than that of $\mathfrak{L}^J_{\bl{b}}$.  We start with $J=\emptyset$ and $\bl{b}=\bl{0}$ and apply the condition in the lemma and (\ref{ddp}) repeatedly. Since the dimension of $\mathfrak{L}^{\emptyset}_{\bl{0}}=\nef(X)$  is finite, after finitely many steps we will arrive at  finitely many reduced linear subsets each of which has zero dimension, i.e. contains only one nef line bundle. It is easy to see any line bundle belonging to one of the reduced linear subsets with positive dimensions appearing in the process satisfies $\mathscr{P}$.

\end{proof}


\subsection{A Class of General Toric Varieties}\hfill

In this work we will consider the following class of general toric varieties. Note that our definition differs from \cite[Definition 2.1]{fs2009}, where the condition is exerted  on  codimension one facets.

\begin{definition}\label{gtv} A projective toric variety $X$ is said to be general if the polytope associated to an ample line bundle on it has no paralleling edges.
\end{definition}

To prove the proposition above, we need the following
\begin{lem}\label{sdrc}
Let $P\subset M_{\mathds{R}}$ be a strongly convex lattice polyhedron with $\dim P=\dim M_{\mathds{R}}=n$ such that no two edges of $P$ are parallel, then

1) there exists a one-to-one correspondence between the infinie edges (resp. codimension one facets) of $P$ and its recession cone $C$;

2) $\dim C=n$.

\end{lem}

\begin{proof} We first prove 1) for the case when $\dim C=n$. Let $\Sigma_P, \Sigma_C\subseteq N_{\mathds{R}}$ be the normal fans of $P$ and $C$ respectively, then by \cite[Theorem 7.1.6]{cox2011toric} $|\Sigma_C|=|\Sigma_P|=C^{\vee}$ and the cones of $\Sigma_C$ with lower dimensions are contained in $\partial C^{\vee}$. In particular, we have $\bigcup_{\tau\in\Sigma_C(n-1)}\tau=\partial C^{\vee}$. On the other hand, since $\Sigma_P$ is a refinement of $\Sigma_C$ and there are no paralleling edges on $P$, $\Sigma_C(n-1)$ can be regarded as a subset of $\Sigma_P(n-1)$. Then each infinite edge of $C$ corresponds uniquely to an infinite edge of $P$. Next we prove $\Sigma_C(n-1)=\Sigma_P(n-1)$. Let $\tau\in\Sigma_P(n-1)$ such that it is not contained in $\partial C^{\vee}$, then we can find $\rho\in\Sigma_P(1)$ in the interior of $C^{\vee}$ such that $\rho\preceq\tau$. However, for such $\rho$ it is easy to see the corresponding facet is bounded, hence the edge corresponding to $\tau$ contained in it is also finite, a contradiction. 

It remains to show there also exists a correspondence between the codimension one facets. We only need to consider the one-dimensional cones of $\Sigma_P$ lying on the boundary as those in the interior of $|\Sigma_P|$ correspond to bounded facets. By abuse of notation, let $\rho\in\Sigma_P(1)\cap\partial|\Sigma_P|$, then since $\partial C^{\vee}=\bigcup_{\tau\in\Sigma_C(n-1)}\tau$, one  can find $\tau\in\Sigma_C(n-1)\subseteq\Sigma_P(n-1)$ such that $\rho$ is contained in $\tau$. As $\Sigma_P$ is a refinement of $\Sigma_C$ and there are no paralleling edges on $C$ we must have $\rho\preceq\tau$ and $\rho\in\Sigma_C(1)$.

Next we prove $C$ must be full-dimensional. If $\dim C<\dim C^{\vee}$, then $|\Sigma_C|=\mathfrak{U}+\mathfrak{N}$, where $\mathfrak{U}$ is a strongly convex cone and $\mathfrak{N}\subset N_{\mathds{R}}$ is the nonzero subspace  consisting of vectors vanishing on $C$.  On the other hand, by the proof loc. cit. we have 
\[C^{\vee}=|\Sigma_P|=\bigcup_{v\in V}\sigma_v,
\] 
where $V$ is the set of vertices of $P$ and $\sigma_v$ is dual of the cone generated by $P-v$. Then we get 
\[\mathfrak{N}=\bigcup_{v\in V}(\sigma_v\cap\mathfrak{N}).
\]
We claim
\begin{equation}\label{mfkn}\mathfrak{N}=\bigcup_{\substack{\sigma\in\Sigma_P\\\sigma\subset\mathfrak{N}}}\sigma.
\end{equation}
Indeed, let $\text{rin}(\sigma)$ be the relative interior of $\sigma$, then it suffices to show $\rin(\sigma)\cap\mathfrak{N}\neq\emptyset$ implies $\sigma\subset\mathfrak{N}$. Let $\sigma(1)=\{\rho_1, \rho_2,\cdots, \rho_k\}$ and $\rho\in\rin(\sigma)\cap\mathfrak{N}$, then $\rho=\sum_{1\leq i\leq k}a_i\rho_i$ with $a_i>0, 1\leq i\leq k$. Suppose further $\sigma\preceq\sigma_v$, then we have $\langle x,\rho_i\rangle\geq0$ for all $x\in C$ as $\sigma_v\subseteq C^{\vee}$. If $\sigma\not\subset\mathfrak{N}$ we have $\langle x,\rho_i\rangle>0$ for some $1\leq i\leq k$ and $x\in C$. However, this would imply $\langle x, \rho\rangle>0$, which contradicts with our assumption that $\rho\in\mathfrak{N}$, i.e. $\langle x,\rho \rangle=0$ for all $x\in C$.

From (\ref{mfkn}) one deduces $\mathfrak{N}$ is spanned by \[\{\rho\:|\:\rho\in\Sigma_P(1)\cap\mathfrak{N}\}.\]

In particular, there are at least two cones $\sigma_1, \sigma_2 \in\Sigma_P$ contained in $\mathfrak{N}$ such that $\dim\sigma_1=\dim\sigma_2=\dim \mathfrak{N}$. Then $\dim F_{\sigma_1}(P)=\dim F_{\sigma_2}(P)=\dim C$ and the minimal affine subspaces containing $F_{\sigma_1}(P)$ and $F_{\sigma_2}(P)$ respectively are translations of each other. On the other hand, one sees easily $C$ is the recession cone of $F_{\sigma_i}(P), i=1, 2$, then 
by what we have proved above for 1), both of the infinite edges of $F_{\sigma_1}(P)$ and $F_{\sigma_1}(P)$ have a one-to-one correspondence with those of $C$. By condition, however, no two edges of $P$ are parallel, hence $F_{\sigma_1}(P)$ and $F_{\sigma_2}(P)$ must coincide and $\dim C=n$. As a result, part 1) of the lemma is also proved. 

\end{proof}

\begin{pro}\label{sf}Let $X$ be a general toric variety of dimension $n$, then 

1) any nontrivial nef line bundle on $X$ is big;

2) for nontrivial nef line bundles $\ul{L}, L$ such that $\ul{L}\prec L$,  $\Sigma_{\ul{L}}(n-1)$ can be regarded as a subset of  $\Sigma_{\ul{L}}(n-1)$;

3) for $\tau\in\Sigma_{\ul{L}}(n-1)\subseteq\Sigma_{L}(n-1)$, there is a one-to-one correspondence between the following sets
\[\{\ub{\rho}\in\Sigma_{\ul{L}}(1)\;|\;\ub{\rho}\prec \tau\},\;\;\;\{\rho\in\Sigma_{L}(1)\;|\;\rho\prec \tau\};
\]

\end{pro}

\begin{proof} Let $L_1$ be an ample and $L_2$ a nontrivial nef line bundle on $X$. Take a full-dimensional cone $\sigma$ of $\Sigma_{L_2}$ and let $P=P_{L_1}+\sigma^{\vee}$, then by Lemma \ref{sdrc} 2) we get $\dim P=\dim\sigma^{\vee}=\dim\Sigma_{L_2}=n$.

For the second part, let $\ub{\sigma}$, $\tau\in\Sigma_{\ul{L}}(n-1)$ such that $\tau\preceq\ub{\sigma}$. Recall that (\ref{pv}) $\Sigma_{\ul{L}}(n)$ defines a partition of $\Sigma_L(n)$, then by Lemma \ref{sdrc} 1) there exists a unique $\sigma\in\Sigma_L(n)_{\ub{\sigma}}$ and an edge of $P_{L, \ub{\sigma}}=P_L+\ub{\sigma}^{\vee}$ through $v_{\sigma}(P_{L, \ub{\sigma}})$ that parallels with $e_{\tau}(\ub{\sigma}^{\vee})$. Thus $\tau$ corresponds uniquely to an cone of $\Sigma_L(n-1)$. As there are no paralleling edges on $P_L$ or $P_{\ul{L}}$, one sees easily this map from $\Sigma_{\ul{L}}(n-1)$ to  $\Sigma_L(n-1)$ is well-defined and injective.

To prove the third part,  it suffices to observe that by Lemma \ref{sdrc} 1) there exists a one-to-one correspondence between the codimension one facets of $P_{L,\ub{\sigma}}=P_L+\ub{\sigma}^{\vee}$ containing the edge $e_{\tau}(P_{L,\ub{\sigma}})$ and those of $\ub{\sigma}^{\vee}$ that contain $e_{\tau}(\ub{\sigma}^{\vee})$.

\end{proof}

\section{Coverings of Polytopes and Proof of Theorem \ref{t1}}
In this section we will prove Theorem \ref{t1} first. Afterwards, we will set forth a local variant of Oda's question, which is evoked by the proof of the theorem.  We will see there some essential difficulties will naturally arise for non-general toric varieties even when we consider this weaker question. 

\subsection{Coverings of Polytopes via Vertex Fitting}\hfill

Next we will show the implication (\ref{o1}) is true as long as $\mathbb{E}(L)$ is large enough. That is, we have the following 
\begin{pro}\label{p1} For a given sufficiently ample line bundle $L_1$, the multiplication map $\Phi_{L_1, L_2}$ is surjective for any nef line bundle $L_2$. 
\end{pro}
In order to prove the proposition above, we need the following 

\begin{lem}\label{lem1} Let $P$ be a convex polytope with dimension $d$, then for any $\tfrac{d}{d+1}<c<1$, we have  
\begin{equation}\label{c}P=\bigcup_{v\in V}cP-cv+v,\end{equation}
where $V$ is the set of vertices of $P$. 
\end{lem}

\begin{proof} It is easy to observe that the lemma is true when $P$ is a simplex. More general case can be deduced from this observation by using Carath\'eodory's theorem \cite[Theorem 1.53]{bg2009}. 
\end{proof}

\begin{proof}[Proof of Proposition \ref{p1}]  Since $\nef(X)_{\mathbb{Q}}$ is a strongly convex rational polyhedral cone we can find a Hilbert basis
$\{B_j\}_{1\leq j\leq m}$ for it such that any nef line bundle can be written as a non-negative integral linear combination of $B_j$. Now for each invariant curve $C$ of $X$ we take
\begin{equation}n_C=\max_jB_j.C
\end{equation}
and let 
\begin{equation}\label{dx}\mathfrak{D}(X)=\{L\in\nef(X)\;|\;\tfrac{L.C}{L.C+n_C}>\tfrac{d}{d+1}\;\text{for any invariant curve}\; C\}.
\end{equation}
Then for any $L\in\mathfrak{D}(X)$ and $1\leq j\leq m$, we have  $\mathbb{E}(\tfrac{L}{L\otimes B_j})>\tfrac{d}{d+1}$ hence $P_L\prec_c P_{L\otimes B_j}$ by Lemma \ref{lem1}. Note that $L\otimes B_j\in\mathfrak{D}(X)$, then one can deduce $L\prec_c L\otimes L'$ for any $L'\in\nef(X)$, which just implies the map $\Phi_{L, L'}$ is surjective.
\end{proof}

\begin{pro}\label{p3}There are finitely many pairs of nef line bundles on $X$ such that  the surjectivity of (\ref{sp}) for each of these pairs will yield an affirmative answer to Oda's question.
\end{pro}

\begin{proof} By a similar argument as in the proofs of Proposition \ref{p1}, one can show for any linear subset $\mathfrak{L}^J_{\bl{b}}$ there exists $\bl{b}'>\bl{b}$ such that 
\begin{equation}\label{sd5}L_1\prec L_2\Rightarrow L_1\prec_c L_2
\end{equation} 
for any $L_1, L_2\in\mathfrak{L}^J_{\bl{b}'}$.
Then by Lemma \ref{pfr}  one can find finitely many reduced linear subsets $\{\mathfrak{L}^{J_k}_{\bl{b}_k}\}_k$ of $\nef(X)$ such that (\ref{sd5}) is valid  for any two line bundles $L_1, L_2$ in one such subset. In particular, let 
\[\mathfrak{L}^{J_k}_{\bl{b}_k}=L_k+\mathfrak{L}^{J_k}_{\bl{0}},\] 
then $L\in \mathfrak{L}^{J_k}_{\bl{b}_k}$ implies $L_k\prec_c L$, hence for any nef line bundle $L$, we have $L\succ_c L_k$
for some $L_k$ from this finite set $\{L_k\}_k$.

Next suppose the map (\ref{sp}) is surjective for any pair consisting of a nef and an ample line bundle from $\{L_k\}_k$, we will prove the same conclusion is also true for other such pairs of line bundles from $\nef(X)$ and $\amp(X)$ respectively. Suppose $L_1$ is a ample and $L_2$ is an nef line bundle, then we can find $L_{k_1}, L_{k_2}\in\{L_k\}_k$ such that $L_{k_1}\prec_c L_1, L_{k_2}\prec_c L_2$. Then
\begin{equation}\label{12pq}P_{L_1}=\bigcup_{m_1\in M_1}m_1+P_{L_{k_1}}, \;\;\;P_{L_2}=\bigcup_{m_2\in  M_2}m_2+P_{L_{k_2}},
\end{equation}
where $M_i=P_{L_i\otimes L_{k_i}^{-1}}\cap M, i=1,2 $. By \cite[Lemma 3.1]{haase2008lattice}, to show  (\ref{sp}) is surjective for $L_1, L_2$ we need only to show for any $m\in M$, $P_{L_1}\cap (m+P^-_{L_2})\neq\emptyset$ implies $P_{L_1}\cap (m+P^-_{L_2})\cap M\neq\emptyset$. By (\ref{12pq}) if $P_{L_1}\cap (m+P^-_{L_2})\neq\emptyset$, then 
\[(m_1+P_{L_{k_1}})\cap (m+m_2+P^-_{L_{k_2}})\neq\emptyset
\]
for some $m_1\in M_1$ and $m_2\in M_2$. Since the map (\ref{sp}) is surjective for $L_{k_1}$ and $L_{k_2}$ by assumption, by using \cite[Lemma 3.1]{haase2008lattice} again one deduces
\[(m_1+P_{L_{k_1}})\cap (m+m_2+P^-_{L_{k_2}})\cap M\neq\emptyset,
\]
which in turn implies  $P_{L_1}\cap (m+P^-_{L_2})\cap M\neq\emptyset$. 



\end{proof}

\subsection{Quasi-coverings of Polytopes via Edge Sliding (along the Long Edges)}\hfill

Given a polytope $P$, by a \emph{quasi-covering} of $P$ we mean a finite set of polytopes $P_i$ contained in $P$ with the same dimension such that $P \backslash\bigcup_i P_i$  is contained in a union of small neighborhoods (compared with the size of $P$) of the vertices of $P$.
The word `small' would be clear from the statement of Claim \ref{qcl3}, in which setting it just means bounded for a family of infinitely many polytopes.

Next we will show the following 

\begin{pro}\label{qcd} With at most finitely many exceptions, the associated polytope of any nef line bundle on a general toric variety $X$ admits a quasi-covering by polytopes corresponding to line bundles from $\dnef(X)$ (see (\ref{dx}) for the definition). 
\end{pro}


By Lemma \ref{pfr}, we only need to show for any reduced linear subset $\mathfrak{L}^J_{\bl{b}}$ there exists $\bl{b}'>\bl{b}$ such that the conclusion of the proposition above is valid for all line bundles in $\mathfrak{L}^J_{\bl{b}'}$. 
Before proving this assertion, we introduce some notations.

 Let $\ul{L}$ (resp. $L_0$) be a line bundle 
in the relative interior of $\mathfrak{L}^J_{\bl{0}}$ (resp. $\mathfrak{L}^J_{\bl{b}}$), $\Sigma_J$ be the normal fan of $P_{\ul{L}}$, which is evidently independent of the choice of $\ul{L}$. 
Let $\hat{L},\tilde{L}$ be two nef line bundles such that $L_0\prec \hat{L}$, $L_0\prec\tilde{L}$, and $ \Sigma_{\hat{L}},\Sigma_{\tilde{L}}$ be the normal fans of $P_{\hat{L}}$ and $P_{\tilde{L}}$ respectively, then by Proposition \ref{sf} 1) we have $\dim\Sigma_J=\dim\Sigma_{\hat{L}}=\dim\Sigma_{\tilde{L}}=d=\dim X$.
Let $\ub{\sigma}\in\Sigma_J(d)$, $\tau\in\Sigma_J(d-1)$ such that $\tau\prec\ub{\sigma}$, then by Proposition \ref{sf} 2) and (\ref{pv}), there exists a unique $\hat{\sigma}\in\Sigma_{\hat{L}}(d)_{\ub{\sigma}}$ (resp. $\tilde{\sigma}\in \Sigma_{\tilde{L}}(d)_{\ub{\sigma}}$) such that $\tau\preceq\hat{\sigma}$ (resp. $\tau\preceq\tilde{\sigma}$).  For saving notations, in the following we will take 
\begin{equation}\label{tp}P_{\tilde{L}}(v_{\ub{\sigma}}(P_{\hat{L}}))_{\tau}=v_{\hat{\sigma}}(P_{\hat{L}})-v_{\tilde{\sigma}}(P_{\tilde{L}})+P_{\tilde{L}}.\end{equation}
When $\hat{L}, \tilde{L}$ are contained in $\mathfrak{L}^J_{\bl{b}}$, (\ref{tp}) is independent of $\tau$, we will just write $P_{\tilde{L}}(v_{\ub{\sigma}}(P_{\hat{L}}))$ for $P_{\tilde{L}}(v_{\ub{\sigma}}(P_{\hat{L}}))_{\tau}$.


\begin{proof}[Proof of Proposition \ref{qcd}] Let $\tilde{L}_0\in\mathfrak{D}(X)$ such that $\tilde{L}_0\succ L_0$,  we will choose $L_2\succ L_0$ such that $P_{L_2}$ admits a quasi-covering by translations of $P_{\tilde{L}_0}$ and prove for any line bundle $L$ in $\mathfrak{L}^J_{\bl{b}}$ such that $L\succ L_2$ the polytope $P_L$  admits a quasi-covering by translations of $P_{\tilde{L}_0\otimes L\otimes L_2^{-1}}$. The proof will be given in three steps. 

Step 1. We will show the following
\begin{claim}\label{qcl1}There exists $L_1\in\mathfrak{L}^J_{\bl{b}}$, $L_1\succ L_0$ such that for each $\tau\in\Sigma_J(d-1), \ub{\sigma}\in\Sigma_J(d)$ with $\tau\preceq\ub{\sigma}$, one can find a non-negative integer $N(\tau)$ satisfying
 \begin{equation}\label{jm}P_{\tilde{L}_0}(v_{\ub{\sigma}}(P_{L_1}))_{\tau}+N(\tau)u_{\tau}(\ub{\sigma})\subseteq P_{L_1}.\end{equation}
\end{claim}
It suffices to prove the existence of $N(\tau)$ for given $\tau$ and $\ub{\sigma}$ when $\mathbb{E}(P_{L_1\otimes L_0^{-1}})\gg0$.  Given $\alpha\in\Sigma_{L_1}(1)=\Sigma_{L_0}(1)$ such that $\alpha\preceq\sigma$ for some $\sigma\in\Sigma_{L_1}(d)_{\ub{\sigma}}=\Sigma_{L_0}(d)_{\ub{\sigma}}$, let $H_{\alpha}(P_{L_1})$ be the hyperplane containing the facet  $F_{\alpha}(P_{L_1})$. If $\alpha\preceq \tau$, then by Proposition \ref{sf} 3) $P_{\tilde{L}_0}(v_{\ub{\sigma}}(P_{L_1}))_{\tau}+N(\tau)u_{\tau}(\ub{\sigma})$ is located on the same side of $H_{\alpha}(P_{L_1})$ with $P_{L_1}$ for any $N(\tau)$. If $\alpha\npreceq \tau$, by Proposition \ref{sf} 3) we have
\[d(H_{\alpha}(P_L), P_{\tilde{L}_0}(v_{\ub{\sigma}}(P_{L_1}))_{\tau}+N(\tau)u_{\tau}(\ub{\sigma}))\gg0
\]
when $N(\tau)$ is sufficiently large, hence in this case the translation $P_{\tilde{L}_0}(v_{\ub{\sigma}}(P_{L_1}))_{\tau}+N(\tau)u_{\tau}(\ub{\sigma})$ also lies on the same side of $H_{\alpha}(P_{L_1})$ with $P_{L_1}$.  
 
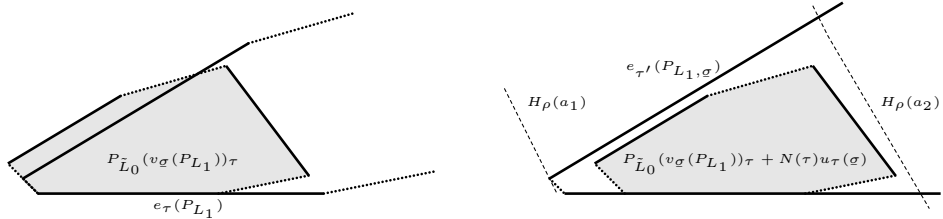
\begin{figure}[H]
\begin{tikzpicture}

\draw [line width=1pt] (-8,0)--(-4.2,0);

\draw [line width=1pt] (-8.2,0.2)--(-5.2,2);

\draw[line width=1pt, line cap=round, dash pattern=on 0pt off 1.5\pgflinewidth](-8.2,0.2)--(-8,0);

\draw[line width=1pt, line cap=round, dash pattern=on 0pt off 1.5\pgflinewidth](-5.2,2)--(-3.8,2.4);

\draw[line width=1pt, line cap=round, dash pattern=on 0pt off 1.5\pgflinewidth](-4.2,0)--(-2.7,0.3);


\draw[line width=1pt, line cap=round, dash pattern=on 0pt off 1.5\pgflinewidth](-8,0)--(-8.4,0.4);

\draw [line width=1pt] (-8.4,0.4)--(-6.9,1.3);

\draw [line width=1pt] (-5.5,1.7)--(-4.4, 0.24);

\draw[line width=1pt, line cap=round, dash pattern=on 0pt off 1.5\pgflinewidth](-6.9,1.3)--(-5.5,1.7);

\draw[line width=1pt, line cap=round, dash pattern=on 0pt off 1.5\pgflinewidth](-5.6,0)--(-4.4, 0.24);


\draw [line width=1pt] (-1,0)--(4,0);

\draw [line width=1pt] (-1.2,0.2)--(2.7,2.54);

\draw[line width=1pt, line cap=round, dash pattern=on 0pt off 1.5\pgflinewidth](-1.2,0.2)--(-1,0);



\draw[dash pattern=on 1.5pt off 1pt](2.3, 2.5)--(2.5,2.1)--(3.85,-0.2);

\draw[dash pattern=on 1.5pt off 1pt](-1.8, 1.4)--(-1.2,0.2)--(-1.1,0);
\draw[line width=1pt, line cap=round, dash pattern=on 0pt off 1.5\pgflinewidth](-0.2,0)--(-0.6,0.4);

\draw [line width=1pt] (-0.6,0.4)--(0.9,1.3);

\draw [line width=1pt] (2.3,1.7)--(3.4, 0.24);

\draw[line width=1pt, line cap=round, dash pattern=on 0pt off 1.5\pgflinewidth](0.9,1.3)--(2.3,1.7);

\draw[line width=1pt, line cap=round, dash pattern=on 0pt off 1.5\pgflinewidth](2.2,0)--(3.4, 0.24);

\node at(1.4,0.4) (0) {\relsize{-100}$P_{\tilde{L}_0}(v_{\ub{\sigma}}(P_{L_1}))_{\tau}+N(\tau)u_{\tau}(\ub{\sigma})$}; 

\node at(-6.2,0.4) (0) {\relsize{-100}$P_{\tilde{L}_0}(v_{\ub{\sigma}}(P_{L_1}))_{\tau}$}; 

\node at(0.45,1.65) (0) {\relsize{-10}$e_{\tau'}(P_{L_1,\ub{\sigma}})$};

\node at(-6,-0.2) (0) {\relsize{-10}$e_{\tau}(P_{L_1})$};

\node at(3.6,1.2) (0) {\relsize{-10}$H_{\rho}(a_2)$};

\node at(-1.1,1.2) (0) {\relsize{-10}$H_{\rho}(a_1)$}; 

\fill[fill opacity=0.1] (-0.2, 0)--(-0.6, 0.4)--(0.9,1.3)--(2.3,1.7)--(3.4,0.24)--(2.2,0);

\fill[fill opacity=0.1] (-8, 0)--(-8.4, 0.4)--(-6.9,1.3)--(-5.5,1.7)--(-4.4,0.24)--(-5.6,0);
\end{tikzpicture}
\end{figure}

Let $P_{L_1,\ub{\sigma}}=P_{L_1}+\ub{\sigma}^{\vee}$, then by our argument above we have
 \[P_{\tilde{L}_0}(v_{\ub{\sigma}}(P_{L_1}))_{\tau}+N(\tau)u_{\tau}(\ub{\sigma})\subset P_{L_1,\ub{\sigma}}.\]
Take a nonzero vector $\rho$ in the interior of $\ub{\sigma}\subset N_{\mathds{R}}$ and let $H_{\rho}(a)=\{x\in M_{\mathds{R}}\:|\:\langle x,\rho \rangle=a\}$ for each $a\in\mathds{R}$. Then since the recession cone of $P_{L_1,\ub{\sigma}}$ is $\ub{\sigma}^{\vee}$, we can find a minimal number $a_1\in\mathds{R}$ such that $H_{\rho}(a_1)\cap P_{L_1,\ub{\sigma}}\neq\emptyset$ and $a_2>a_1$ such that  $P_{\tilde{L}_0}(v_{\ub{\sigma}}(P_{L_1}))_{\tau}+N(\tau)u_{\tau}(\ub{\sigma})$ is lying between the $H_{\rho}(a_1), H_{\rho}(a_2)$, as is shown in the figure on the right above.

For each $\tau'\in\Sigma_J(d-1), \tau'\preceq\ub{\sigma}$ by our choice of $\rho$ the infinite edge $e_{\tau'}(P_{L_1,\ub{\sigma}})$ which contains $e_{\tau'}(P_{L_1})$ will intersect with $H_{\rho}(a_2)$ when

\[\tfrac{\|e_{\tau'}(P_{L_1})\|}{\|u_{\tau'}\|}|\langle u_{\tau'}, \rho\rangle|\geq a_2-a_1\]
or equivalently
\begin{equation}\label{jil}L_1.C_{\tau'}\geq\tfrac{a_2-a_1}{|\langle u_{\tau'},\rho \rangle|}.\end{equation}
Then $P_{\tilde{L}_0}(v_{\ub{\sigma}}(P_{L_1}))_{\tau}+N(\tau)u_{\tau}$ is the contained in the polytope bounded by $H_{\rho}(a_2)$ and the facets $F_{\alpha}(P_{L_1})$, where $\alpha\preceq\sigma$ for some $\sigma\in\Sigma_{L_1}(d)_{\ub{\sigma}}$. This polytope is obviously contained in $P_{L_1}$  when $\mathbb{E}(L_1\otimes L_0^{-1})\gg0$.
As a result, (\ref{jm}) will follow when the inequality (\ref{jil}) is satisfied for  all $\tau'\in\Sigma_J(d-1), \tau'\preceq\ub{\sigma}$. 

\medskip

Take a line bundle $L\in\mathfrak{L}^J_{\bl{b}}, L\succ L_1$ and add $P_{L\otimes L^{-1}_1}$ to both sides of (\ref{jm}), we get 
\[P_{\tilde{L}_0\otimes L\otimes L^{-1}_1}(v_{\ub{\sigma}}(P_L))_{\tau}+N(\tau)u_{\tau}(\ub{\sigma})\subset P_L.
\]
Note that the line bundle $\tilde{L}_0\otimes L\otimes L^{-1}_1$ falls in $\mathfrak{D}(X)$. If we can find a suitable line bundle $L_2\succ L_1$ such that for all $L\succ L_2$
\begin{equation}\label{ul}\bigcup_{\ub{\sigma}\in\Sigma_J(d)}\bigcup_{\substack{\tau\preceq\ub{\sigma}\\\tau\in \Sigma_J(d-1)}}P_{L_0\otimes L\otimes L^{-1}_2}(v_{\ub{\sigma}}(P_L))+N(\tau)u_{\tau}(\ub{\sigma})\end{equation}
constitutes a quasi-covering of $P_L$, then so does
\begin{equation}\label{u}\bigcup_{\ub{\sigma}\in\Sigma_J(d)}\bigcup_{\substack{\tau\preceq\ub{\sigma}\\\tau\in \Sigma_J(d-1)}}P_{\tilde{L}_0\otimes L\otimes L^{-1}_2}(v_{\ub{\sigma}}(P_L))_{\tau}+N(\tau)u_{\tau}(\ub{\sigma})\end{equation}
since $\tilde{L}_0\otimes L\otimes L^{-1}_2\succ L_0\otimes L\otimes L^{-1}_2$.

Step 2.  Before finding a line bundle $L_2$ as described above, we will prove the distances between the points of $P_L$ in the complement of (\ref{ul}) and the border of $P_L$ can be uniformly bounded.  

\begin{claim}\label{qcl2} 
There exists a constant $c>0$ such that for all $L\in\mathfrak{L}^J_{\bl{b}}$,    $L\succ L_1$, $\tau\in\Sigma_J(d-1), \ub{\sigma}\in\Sigma_J(d), \tau\preceq\ub{\sigma}$ and $x\in P_L\backslash (P_{L_0\otimes L\otimes L^{-1}_1}(v_{\ub{\sigma}}(P_L))+N(\tau)u_{\tau}(\ub{\sigma}))$ we have
\[d(x, \partial P_L)<c.\]
\end{claim}

Since $P_L=P_{L_0\otimes L\otimes L^{-1}_1}+P_{L^{-1}_0\otimes L_1}$ and by our construction $P_{L_0\otimes L\otimes L^{-1}_1}(v_{\ub{\sigma}}(P_L))+N(\tau)u_{\tau}(\ub{\sigma})$ is contained in $P_L$, this translation of $P_{L_0\otimes L\otimes L^{-1}_1}$ corresponds to a lattice point $p$ in 
$P_{L^{-1}_0\otimes L_1}$. Then we have 
\begin{equation}\label{fpq}P_L=\bigcup_{q\in P_{L_0^{-1}\otimes L_1}}\overrightarrow{pq}+P_{L_0\otimes L\otimes L^{-1}_1}(v_{\ub{\sigma}}(P_L))+N(\tau)u_{\tau}(\ub{\sigma}).
\end{equation}
Therefore, for $x\in P_L\backslash(P_{L_0\otimes L\otimes L^{-1}_1}(v_{\ub{\sigma}}(P_L))+N(\tau)e_{\tau}(\ub{\sigma}))$ we can find $z\in P_{L_0\otimes L\otimes L^{-1}_1}(v_{\ub{\sigma}}(P_L))+N(\tau)e_{\tau}(\ub{\sigma})$ such that both $x$ and $z$ lie in some (not necessarily a lattice) translation of $P_{L_0^{-1}\otimes L_1}$ which is contained in $P_L$. Obviously one can assume further $z\in\partial (P_{L_0\otimes L\otimes L^{-1}_1}(v_{\ub{\sigma}}(P_L))+N(\tau)e_{\tau}(\ub{\sigma}))$. As \[d(x, \partial P_L)\leq d(z, \partial P_L)+\max_{x, y\in P_{L_0^{-1}\otimes L_1}}\|\overline{xy}\|,\]
 it suffices to bound $d(z, \partial P_L)$ for any $z\in \partial(P_{L_0\otimes L\otimes L^{-1}_1}(v_{\ub{\sigma}}(P_L))+N(\tau)e_{\tau}(\ub{\sigma}))$. 

Note that  by (\ref{fpq}) for each $\ub{\rho}\in\Sigma_{L_0}(1)$ we can find $q\in P_{L^{-1}_0\otimes L_1}$ such that the affine hyperplanes containing $F_{\ub{\rho}}(P_L)$ and  $\overrightarrow{pq}+F_{\ub{\rho}}(P_{L_0\otimes L\otimes L^{-1}_1}(v_{\ub{\sigma}}(P_L))+N(\tau)e_{\tau}(\ub{\sigma}))$ coincide. Therefore, one gets 
\[d(z, \partial P_L)\leq\max_{q\in P_{L^{-1}_0\otimes L_1}}\|\langle\overrightarrow{pq},\rho\rangle\|,
\]
 for any $z\in F_{\ub{\rho}}(P_{L_0\otimes L\otimes L^{-1}_1}(v_{\ub{\sigma}}(P_L))+N(\tau)e_{\tau}(\ub{\sigma}))$, hence Claim \ref{qcl2} is proved.


Step 3. To complete the proof of the proposition, we will show the following
\begin{claim}\label{qcl3} 
There exists $L_2\in\mathfrak{L}^J_{\bl{b}}, L_2\succ L_1$  such that for any $L\in\mathfrak{L}^J_{\bl{b}}, L\succ L_2$ the complement of 
\begin{equation}\label{u}\bigcup_{\ub{\sigma}\in\Sigma_J(d)}\bigcup_{\substack{\tau\preceq\ub{\sigma}\\\tau\in \Sigma_J(d-1)}}P_{L_0\otimes L\otimes L^{-1}_2}(v_{\ub{\sigma}}(P_L))+N(\tau)u_{\tau}(\ub{\sigma})\end{equation}
in $P_L$ is contained in a union of neighorhoods  of the vertices of $P_L$ whose radii can be uniformly bounded. 

\end{claim}

Let $\delta=\tfrac{d}{d+1}$ and $(L\otimes L_0^{-1})^{\delta}\in\nef(X)_{\mathds{Q}}$ be the line bundle corresponding to the polytope $\delta P_{L\otimes L_0^{-1}}$, one checks easily
\[\mathbb{E}(\tfrac{L_0\otimes (L\otimes L_0^{-1})^{\delta}}{L})>\tfrac{d}{d+1}.
\]
Then by the proof of Proposition \ref{p1} we have
\begin{equation}\label{u2}P_L=\bigcup_{\ub{\sigma}\in\Sigma_J(d)}P_{L_0\otimes(L\otimes L_0^{-1})^{\delta}}(v_{\ub{\sigma}}(P_L)).
\end{equation}
By Claim \ref{qcl2} for any $x$ in the complement of the union (\ref{u}), we have $d(x, \partial P_L)<c$. Let 
\[B(L, \ub{\sigma}, c)=\{x\in P_{L_0\otimes(L\otimes L_0^{-1})^{\delta}}(v_{\ub{\sigma}}(P_L))\;|\;d(x, \partial P_L)<c\},\]
then to prove Claim \ref{qcl3}, we need only to show for each $\ub{\sigma}\in\Sigma_J(d)$ the intersection of  $B(L, \ub{\sigma}, c)$ with the complement of (\ref{u})
 is contained in a bounded neighborhood of $v_{\sigma}(P_L)$ for some hence all $\sigma\in\Sigma_{L_0}(d)_{\ub{\sigma}}=\Sigma_L(d)_{\ub{\sigma}}$. 
 
Now we define
\[\Sigma_{L_0}(1)_{\ub{\sigma}}=\{\rho\in\Sigma_{L_0}(1)\;|\;\rho\prec\sigma\;\text{for some}\; \sigma\in\Sigma_{L_0}(d)_{\ub{\sigma}}\}\]
and
\[\Sigma_{L_0}(1)^{\ub{\sigma}}=\Sigma_{L_0}(1)\backslash\Sigma_{L_0}(1)_{\ub{\sigma}}.\]
It is easy to see for any $x\in P_{L_0\otimes (L\otimes L^{-1}_0)^{\delta}}(v_{\ub{\sigma}}(P_L))$ and all $\alpha\in\Sigma_{L_0}(1)^{\ub{\sigma}}$, 
\[d(x, H_{\alpha}(P_L))
>d(H_{\alpha}(P_{L_0\otimes(L\otimes L_0^{-1})^{\delta}}(v_{\ub{\sigma}}(P_L))), H_{\alpha}(P_L))>c
\]
when $\mathbb{E}(L\otimes L_1^{-1})\gg0$. Therefore we can choose $L_2$ such that for all $L\in\mathfrak{L}^{J}_{\bl{b}}, L\succ L_2$,
\begin{equation}\label{lr}B(L,\ub{\sigma}, c)=\{x\in P_{L_0\otimes (L\otimes L^{-1}_0)^{\delta}}(v_{\ub{\sigma}}(P_L))\;|\;d(x,  F_{\rho}(P_L))<c \text{ for some}\;\rho\in\Sigma_{L_0}(1)_{\ub{\sigma}}\}.\end{equation}

Let $P_{L_0\otimes L\otimes L^{-1}_2, \ub{\sigma}}=P_{L_0\otimes L\otimes L^{-1}_2}(v_{\ub{\sigma}}(P_L))+\ub{\sigma}^{\vee}$, then one sees easily
\[B(L,\ub{\sigma}, c)\backslash (P_{L_0\otimes L\otimes L^{-1}_2}(v_{\ub{\sigma}}(P_L))+N(\tau)u_{\tau}(\ub{\sigma}))=B(L,\ub{\sigma}, c)\backslash (P_{L_0\otimes L\otimes L^{-1}_2, \ub{\sigma}}+N(\tau)u_{\tau}(\ub{\sigma})).\]
Therefore we need only to prove the complement 
\[B(L,\ub{\sigma}, c)\backslash \bigcup_{\substack{\tau\preceq\ub{\sigma},\\\tau\in \Sigma_J(d-1)}}(P_{L_0\otimes L\otimes L^{-1}_2, \ub{\sigma}}+N(\tau)u_{\tau}(\ub{\sigma}))\]
is contained in a neighborhood of $v_{\sigma}(P_L) $ with $\sigma\in\Sigma_{L_0}(d)_{\ub{\sigma}}$ of bounded size. As we have \[B(L, \ub{\sigma}, c)\subseteq P_{L_0\otimes(L\otimes L_0^{-1})^{\delta}}(v_{\ub{\sigma}}(P_L))\subseteq P_{L_0\otimes L\otimes L^{-1}_2, \ub{\sigma}},\]
it suffices to prove the assertion above for the following complement
\[P_{L_0\otimes L\otimes L^{-1}_2, \ub{\sigma}}\backslash \bigcup_{\substack{\tau\prec\ub{\sigma},\\\tau\in \Sigma_J(d-1)}}(P_{L_0\otimes L\otimes L^{-1}_2, \ub{\sigma}}+N(\tau)u_{\tau}(\ub{\sigma})).
\]

Note that $P_{L_0\otimes L\otimes L^{-1}_2, \ub{\sigma}}=P_{L_0}+{\ub{\sigma}}^{\vee}$, then this can be further reduced to prove the following complement is contained in  a neighborhood of $y$ of bounded size for all $y\in P_{L_0}$
\begin{equation}\label{compy}y+\ub{\sigma}^{\vee}\backslash \bigcup_{\substack{\tau\preceq\ub{\sigma},\\\tau\in \Sigma_J(d-1)}}(y+\ub{\sigma}^{\vee}+N(\tau)u_{\tau}(\ub{\sigma})).
\end{equation}
Moreover, since $y+\ub{\sigma}^{\vee}\backslash \bigcup_{\tau}(y+\ub{\sigma}^{\vee}+N(\tau)u_{\tau}(\ub{\sigma}))$ is a translation of the following complement, we only need to prove the same conclusion for it
\begin{equation}\label{comp}\ub{\sigma}^{\vee}\backslash \bigcup_{\substack{\tau\preceq\ub{\sigma},\\\tau\in \Sigma_J(d-1)}}(\ub{\sigma}^{\vee}+N(\tau)u_{\tau}(\ub{\sigma})).
\end{equation}

Now as in the proof of Claim \ref{qcl1}, let $\rho$ be a nonzero vector in the interior of $\ub{\sigma}\subset N_{\mathds{R}}$ and for each $a\in\mathds{R}$ let $H_{\rho}(a)=\{x\in M_{\mathds{R}}\:|\:\langle x,\rho \rangle=a\}$. Then by our choice of $\rho$, the translation $H_{\rho}(0)+\lambda u_{\tau}(\ub{\sigma})$ has nonempty intersection with all infinite edges of $\ub{\sigma}^{\vee}+N(\tau)u_{\tau}(\ub{\sigma})$ for each $\tau$ when $\lambda$ is sufficiently large. Moreover, for all large $\lambda_1, \lambda_2$ the intersection $(H_{\rho}(0)+\lambda_1 u_{\tau}(\ub{\sigma}))\cap \ub{\sigma}^{\vee}$ is a homothety of $(H_{\rho}(0)+\lambda_2 u_{\tau}(\ub{\sigma}))\cap \ub{\sigma}^{\vee}$ by a factor of $\tfrac{\lambda_1}{\lambda_2}$.  In particular, we have 
\[\mathbb{E}\Big(\tfrac{(H_{\rho}(0)+\lambda u_{\tau}(\ub{\sigma}))\cap (\ub{\sigma}^{\vee}+N(\tau)u_{\tau}(\ub{\sigma}))}{(H_{\rho}(0)+\lambda u_{\tau}(\ub{\sigma}))\cap \ub{\sigma}^{\vee}}\Big)=\mathbb{E}\Big(\tfrac{(H_{\rho}(0)+(\lambda-N(\tau))u_{\tau}(\ub{\sigma}))\cap \ub{\sigma}^{\vee}}{(H_{\rho}(0)+\lambda u_{\tau}(\ub{\sigma}))\cap \ub{\sigma}^{\vee}}\Big)>\tfrac{d-1}{d},\]
for large $\lambda$. Therefore, by the proof of Proposition \ref{p1} we have 
\[(H_{\rho}(0)+\lambda u_{\tau}(\ub{\sigma}))\cap \ub{\sigma}^{\vee}\subseteq\bigcup_{\tau}(H_{\rho}(0)+\lambda u_{\tau}(\ub{\sigma}))\cap (\ub{\sigma}^{\vee}+N(\tau)u_{\tau}(\ub{\sigma}))\]
for all sufficiently large $\lambda$, hence the complement (\ref{comp}) is contained in a bounded neighborhood of the vertex of $\ub{\sigma}^{\vee}$. Therefore, the union of (\ref{compy}) for all $y\in P_{L_0}$ is also contained in the bounded neighborhood of $v_{\sigma}(P_L)$, hence Claim \ref{qcl3} and the proposition are proved.

\end{proof}

\begin{proof}[Proof of Theorem \ref{t1}]
By Proposition \ref{qcd} and Lemma \ref{pfr}, there exist finitely many linear subsets $\{\mathfrak{L}^{J_k}_{\bl{b}_k}\}$ such that for any line bundle from one of them its associated polytope admits a quasi-covering by polytopes associated to line bundles from $\mathfrak{D}(X)$. Note that (see the proof of Lemma \ref{pfr}) still some hence finitely many of these linear subsets might contain only one line bundle, in this case no polytopes will be jammed into the associated polytope to form a quasi-covering. 

Given $L_i\in\mathfrak{L}^{J_{k_i}}_{\bl{b}_{k_i}}, i=1,2$, we can find a quasi-covering of $P_{L_i}$, say, by $\{Q_i+x_i\}_{x_i\in M_i}$, where $Q_i$ corresponds to a line bundle in $\mathfrak{D}(X)$ and $M_i\subset M$ is a  finite subset. By Proposition \ref{p1} we have 
\[(Q_1\cap M)+(P_{L_2}\cap M)=(Q_1+P_{L_2})\cap M,\;\;\;(Q_2\cap M)+(P_{L_1}\cap M)=(Q_2+P_{L_1})\cap M.\]
Now given $u_1\in P_{L_1}, u_2\in P_{L_2}$ such that $u_1+u_2\in P_{L_1\otimes L_2}\cap M$, suppose there exist no $v_1\in P_{L_1}\cap M$ and $v_2\in P_{L_2}\cap M$ such that $v_1+v_2=u_1+u_2$, then $u_1$ and $u_2$ must be contained in $P_{L_1}\backslash\bigcup_{x_1} (Q_1+x_1)$ and $P_{L_2}\backslash\bigcup_{x_2} (Q_2+x_2)$ respectively. By Proposition \ref{qcd}, $P_{L_1}\backslash\bigcup_{x_1} (Q_1+x_1)$ and $P_{L_2}\backslash\bigcup_{x_2} (Q_2+x_2)$ are contained in a union of neighborhoods of the vertices bounded by some fixed number. It is easy to see one can deduce the same conclusion when one of $L_1, L_2$ is the only line bundle in the linear subset that contains it. Therefore, we can get finally a finite number bounding $\dim\coker\:\Phi_{L_1, L_2}$ for all $L_1, L_2$. 

\end{proof}

\subsection{Local Variants of Oda's Question}\hfill

It should be perceptible that the following question can be viewed as a local variant of Oda's original one. 

\begin{question}\label{joq}Let $L_1$ be an ample line bundle on a smooth projective toric variety $X$,  $L_2$ and $\ul{L}$ be nef line bundles and $\ub{\sigma}$ be a full dimensional cone of $\Sigma_{\ul{L}}$, is the following map surjective
\begin{equation}\label{floq}(P_{L_1}\cap M)\times ((P_{L_2}+\ub{\sigma}^{\vee})\cap M)\rightarrow (P_{L_1\otimes L_2}+\ub{\sigma}^{\vee})\cap M?
\end{equation}
\end{question}

The following strengthened version of Question \ref{joq} can be regardeded as a counterpart of (\ref{os}) in the local setting. 

\begin{question}\label{qjoq}
Under the same conditions as Question \ref{joq}, is the following map surjective
\begin{equation}\label{sloq}P_{L_1}\times ((P_{L_2}+\ub{\sigma}^{\vee})\cap M)\rightarrow P_{L_1\otimes L_2}+\ub{\sigma}^{\vee}?
\end{equation}
\end{question}

Note that $\coker\:\Phi_{L_1, L_2\otimes \ul{L}^k}=0$ for all $k$ sufficiently large would imply the surjectivity of (\ref{floq}). On the other hand, by the proof of Theorem \ref{t1}, the lattice points of $P_{L_1\otimes L_2}+\ub{\sigma}^{\vee}$ staying away from the image (\ref{floq}) is finite if $X$ is  general in the sense of Definition \ref{gtv}. However, a smooth projective variety is general only if each of the two-dimensional facets of the polytope associated to an ample line bundle is a triangle.    

Next we consider to what extent  the map (\ref{floq}) might fail to be surjective in a setting less restrictive than being general. Firstly, we assume $\dim\Sigma_{\ul{L}}=\dim\Sigma_X$. Moreover, we require that each two-dimensional facet of $\ub{\sigma}^{\vee}$ corresponds uniquely to a two-dimensional facet of $P_{L_1}+\ub{\sigma}^{\vee}$, where the former is the recession cone of the latter.
To make the second condition more explicit, let $\Sigma_{X,\ub{\sigma}}$ be the normal fan of $P_{L_1}+\ub{\sigma}^{\vee}$ and for each $\ub{\rho}\in\ub{\sigma}(n-2)$ where $n=\dim\Sigma_X$, we define 
\begin{equation}\label{ac2f}\Sigma_{X,\ub{\sigma}}(n-2)_{\ub{\rho}}=\{\rho\in\Sigma_{X,\ub{\sigma}}(n-2)\;|\:\lim_{t\rightarrow 0^+}F_{\rho}(\ub{\sigma}^{\vee}+tP_{L_1})=F_{\ub{\rho}}(\ub{\sigma}^{\vee})\},
\end{equation}
then this condition just means $\#\Sigma_{X,\ub{\sigma}}(n-2)_{\ub{\rho}}=1$ for any $\ub{\rho}\in\ub{\sigma}(n-2)$.

\begin{pro} \label{cas3}Suppose $\dim\Sigma_X=\dim\Sigma_{\ul{L}}=n$ and for each $\ub{\rho}\in\ub{\sigma}(n-2)$ we have $\#\Sigma_{X,\ub{\sigma}}(n-2)_{\ub{\rho}}=1$,
\noindent then the lattice points of $P_{L_1\otimes L_2}+\ub{\sigma}^{\vee}$ lying outside the image of  (\ref{floq}) are contained in certain tubular neighborhoods of its infinite edges. 

\end{pro}
\begin{proof}
Let  $\iota_i, 1\leq i\leq t$ be the $(n-2)$-dimensional cones of $\Sigma_{X,\ub{\sigma}}$ corresponding to 2-dimensional infinite facets of $P_{L_1}+\ub{\sigma}^{\vee}$ which contain two infinite edges with different directions. Since $X$ is smooth and $L_1$ is ample, each $\iota_i$ is the convex hull of $n-2$ one dimensional cones of $\Sigma_X$. In particular, it can be regarded as an element of $\bigwedge^{n-2}N_{\mathds{R}}$ and by abuse of notation we still denote this element by $\iota_i$. By condition the elements $\{\iota_i\}_{1\leq i\leq t}$ of $\bigwedge^{n-2}N_{\mathds{R}}$ are mutually different. 

Take $x$ in the interior of $\ub{\sigma}$ and let $v_x\in N_{\mathds{R}}$ be corresponding vector. Then it is easy to see by choosing $x$ suitably we may assume any two of the vectors $v_x\wedge\iota_i, 1\leq i\leq t$  are non-colinear. For simplicity, we write $v$ for $v_x$ from now on  and for each number $a$ let \[H_v(a)=\{z\in M_{\mathds{R}}\:|\:\langle z, v\rangle=a\}.\]
Then it is easy to see for all $a$ larger than some number $a_0$, the hyperplane $H_v(a)$ intersects with each infinite edge of $P_{L_1\otimes L_2}+\ub{\sigma}^{\vee}$. Moreover, for any $a_0$ the complement of $\bigcup_{a\geq a_0}H_v(a)\cap(P_{L_1\otimes L_2}+\ub{\sigma}^{\vee})$ in $P_{L_1\otimes L_2}+\ub{\sigma}^{\vee}$ is contained in some bounded neighorhood of the vertices of $P_{L_1\otimes L_2}+\ub{\sigma}^{\vee}$. 

Note that each edge of 
$H_v(a)\cap(P_{L_1\otimes L_2}+\ub{\sigma}^{\vee})$ is an intersection of $H_v(a)$ with a two-dimensional infinite facet of $P_{L_1\otimes L_2}+\ub{\sigma}^{\vee}$. Moreover, we have
\[\lim_{a\rightarrow\infty}\|H_v(a)\cap F_{\iota_i}(P_{L_1\otimes L_2}+\ub{\sigma}^{\vee})\|=\infty, \;\;1\leq i\leq t
\]
and the edges of 
$H_v(a)\cap(P_{L_1\otimes L_2}+\ub{\sigma}^{\vee})$ whose lengths keep invariant with $a$ are just those lying on some two-dimensional facet of $P_{L_1\otimes L_2}+\ub{\sigma}^{\vee}$ with a one-dimensional recession cone. 

By our assumption, except for those with fixed lengths all other edges of 
$H_v(a)\cap(P_{L_1\otimes L_2}+\ub{\sigma}^{\vee})$ have different directions. By a similar argument as in the proof of Theorem \ref{t1} and particularly Claim \ref{qcl3}, one can show that for all $a$ sufficiently large the following union forms a quasi-covering of $H_v(a)\cap(P_{L_1\otimes L_2}+\ub{\sigma}^{\vee})$ 
\[\bigcup_{\substack{m_1\in P_{L_1}\cap M\\m_2\in P_{L_2}\cap M}}H_v(a)\cap (m_1+m_2+\ub{\sigma}^{\vee}).
\]
Then the conclusion of the proposition readily follows from this result.  

\end{proof}

\begin{rem}\label{psi}The case when $\dim\Sigma_X>\dim\Sigma_{\ul{L}}$ is even more intricate. Recall that the equivariant morphism $f: X\rightarrow Y=Y(\Sigma_{\ul{L}})$ can be induced by a map of lattices \[\psi_{\ul{L}}: N\rightarrow N'\] 
such that the image of each cone of $\Sigma_X$ under $\psi_{\ul{L}}$ is contained in some cone of $\Sigma_{\ul{L}}\subset N'_{\mathds{R}}$. Let 
$\Sigma_Z=\ker\psi_{\ul{L}}\cap\Sigma_X$, then $\Sigma_Z$ is a subfan of $\Sigma_X$, hence by \cite[Theorem 3.3.4]{cox2011toric}
the inclusion $\Sigma_Z\rightarrow \Sigma_X$ 
defines a morphism
\[g: Z=Z(\Sigma_Z)\rightarrow X.
\]
Since $\psi_{\ul{L}}$ is surjective, by \cite[(3.3.5)]{cox2011toric} $Z$ is isomorphic to the generic fiber of $f$, then the restriction of $L_1$ to $Z$ is also ample.  From the inclusion $\Sigma_Z(1)\subset \Sigma_X(1)$ we have the following surjective map
\[\tilde{g}: P_{L_1} \rightarrow P_{L_1|_Z}.\]
Suppose the following map is surjective as Question \ref{joq} indicates
\[P_{L_1}\cap M+\sigma^{\vee}\cap M\rightarrow (P_{L_1}+\sigma^{\vee})\cap M,
\]
then for each $z\in P_{L_1|_Z}$ the following map is necessarily surjective
\begin{equation}\label{ofr}\tilde{g}^{-1}(z)\cap M+\sigma^{\vee}\cap M\rightarrow (\tilde{g}^{-1}(z)+\sigma^{\vee})\cap M.
\end{equation}
Note that the fiber $\tilde{g}^{-1}(z)$ is in general not a lattice polytope. On the other hand, one can show $\wid_{\rho}(\tilde{g}^{-1}(z))$ is a piecewise linear convex function of $z$ for any $\rho\in N'$.  If $\mathbb{E}(\tilde{g}^{-1}(z))$ is sufficiently large for any vertex $z$ of $P_{L_1|_Z}$, then by \cite[Corollary 2.11]{kannan1988covering} and the convexity of width function one can deduce $\tilde{g}^{-1}(z)$ is not lattice-free for any $z\in P_{L_1|_Z}$. However, since even the normal fan of the convex hull of the lattice points of $\tilde{g}^{-1}(z)$ is unclear to us, the surjecitivty of \ref{ofr} in this general setting cannot be touched by our methods.

If the map $\psi_{\ul{L}}$ splits $\Sigma_X$ in the sense of  \cite[Definition 3.3.18]{cox2011toric}, then $f$ is a locally trivial fibration hence $Z$ is also smooth. Let $m=\dim\Sigma_Z$, then for each $\sigma\in\Sigma_Z(m)$, there is a section $Y_{\sigma}$ of $f$ passing through the point of the generic fiber of $f$ corresponding to $\sigma$. Moreover, by \cite[Proposition 3.8]{dirocco2014linear} the polytope $P_{L_1}$ is a twisted Cayley
sum of the polytopes $P_{L_1|_{Y_{\sigma}}}$, $\sigma\in\Sigma_Z(m)$. In this case  for a lattice point $z\in P_{L_1|_Z}$ one can show  the fiber $\tilde{g}^{-1}(z)$ is a lattice polytope with the same normal fan as $Y$. Then the surjectivity of (\ref{ofr}) in this setting can be reduced to answering the same question for lower dimensions.  




\end{rem}

\section{Proof of Theorem \ref{t2}}\hfill

To prove (\ref{o1}) in the setting of smooth toric threefolds, the method used in the previous section is inadequate \textemdash
here we need to move $P_L$ around while keeping one of its two-dimensional facets touching the corresponding facet of $P_{\tilde{L}}$. In other words, the translations of $P_L$ corresponding to lattice points on two-dimensional facets of $P_{L^{-1}\otimes \tilde{L}}$ should also be taken into account (see the proof of Proposition \ref{fe}).

\subsection{The Boundedness of $\dim\coker\:\Phi_{L_1, L_2}$ for Given $L_1$}\hfill

\begin{pro}\label{b2}
For a given ample line bundle $L_1$ and any nef line bundle $L_2$ on a smooth projective toric threefold $X$ the dimension of  $\coker\:\Phi_{L_1, L_2}$ can be bounded by some number independent of $L_2$.
\end{pro}

To prove this proposition we will need the following lemmas, the proof of the first of which is trivial and we omit it. 

\begin{lem}\label{upp}
Two neighboring edges of a smooth lattice polygon which contains at least four lattices are also neighboring edges of a unimodular parallelogram that is contained in this polygon.
\end{lem}

\begin{lem} \label{smw}Let $P$ be a pointed infinite polyhedron, $\partial_fP$ be the union of its finite facets, then we have
 $P=\partial_fP+C$. 
\end{lem}
\begin{proof} The lemma will be proved by induction on the dimension of $P$. The case when $\dim P=1$ is trivial. Let $Q$ be the convex hull of the vertices of $P$, then by Minkowski-Weyl decomposition we have $P=Q+C$. If $Q\subset\partial P$, then it is easy to see $Q\subseteq\partial_fP$ hence the conclusion follows from $P=Q+C$. In the case when $Q\cap P^{\circ}\neq\emptyset$ we will find for any $p\in Q\cap P^{\circ}$ a point $q\in\partial_fP$ such that $p\in q+C$.

Let $l$ be a line through $p$ and parallel to a ray $r^+$ in $C$, then there exists a unique $q\in\partial P$ such that $\{q\}=\partial P\cap l$ and $p\in q+r^+$. If $q\in\partial_fP$, then we are done. Otherwise, $q$ is contained in an infinite facet $F$ of $P$. By induction, we have $x\in\partial_fF$ and $y\in C(F)$ such that $q=x+y$, where $C(F)$ is the recession cone of $F$. Since each finite facet of $F$ is also a finite facet of $P$, we have $x\in\partial_fP$. Moreover, $y\in C$ as $C(F)\subseteq C$ hence the lemma is proved.

\end{proof}

\begin{lem}\label{ccas3}
Let $\Sigma_1$ be a two dimensional complete fan and $\Sigma_2$ be a noncomplete subfan of $\Sigma_1$ with a convex support. Let $P_1$ (resp. $P_2$) be a lattice polygon (resp. infinite lattice polygonal region) whose normal fan is $\Sigma_1$ (resp. $\Sigma_2$). Suppose the length of each finite edge of $P_2$ is no smaller than that of the corresponding edge of $P_1$, then each lattice point in $P_2$ can be found in some lattice translation of $P_1$ which is itself contained in $P_2$. 
\end{lem}
\begin{proof} Let $C$ be the recession cone of $P_2$, then the normal fan of $P_1+C$ is also $\Sigma_2$. Moreover, one sees easily there are only finitely many translations of $P_1+C$ contained in $P_2$ such that one of the edges with finite lengths of such a translation is contained in the corresponding edge of $P_2$. If we denote these translations by $P_1+C+v_i, 1\leq i\leq k$, where $v_i$ is an integral vector, then by Lemma \ref{smw} we get
\begin{equation}\label{cew}\bigcup_i(P_1+C+v_i)=P_2.\end{equation}
On the other hand, we can find a nef and big line budle $L_3$ on the toric surface defined by $\Sigma_1$ such that $C=\ub{\sigma}^{\vee}$ for some $\ub{\sigma}\in\Sigma_{L_3}(2)$. Then $P_1+C=\bigcup_{k\geq1}(P_1+kP_3)$, by the main results of \cite{fakhruddin2002multiplication} and \cite{haase2008lattice} one deduces easily any lattice point of $P_1+C+v_i$ is contained in certain integral translation of $P_1$ which itself is contained in $P_1+C+v_i$. Then the claim in the lemma follows from this result and the equality (\ref{cew}).

\end{proof}

\begin{proof}[Proof of Proposition \ref{b2}] By Lemma \ref{pfr} we need only to show for any reduced linear subset of the form $L_0+\mathfrak{L}^J_{\bl{0}}$ there exists a line bundle $\ul{L}\in \mathfrak{L}^J_{\bl{0}}$ such that $\dim\coker\:\Phi_{L_1, L_2}$ can be uniformly bounded for any $L_2\in\ul{L}\otimes L_0+\mathfrak{L}^J_{\bl{b}}$. 
Moreover, by a similar argument as the proof of Theorem \ref{t1} especially Claim \ref{qcl3}, this can be further reduced to a local question (see Question \ref{joq}). To be more precise, let $\Sigma_J$ be the normal fan of the polytope associated to an line bundle in the interior of $\mathfrak{L}^J_{\bl{0}}$, $n_J=\dim\Sigma_J$ and $\ub{\sigma}\in\Sigma_J(n_J)$. Then we only need to show the lattice point on the right of the map below not contained in the image can be uniformly bounded.
\begin{equation}\label{ssgm}(P_{L_1}\cap M)\times ((P_{L_0}+\ub{\sigma}^{\vee})\cap M)\rightarrow (P_{L_1\otimes L_0}+\ub{\sigma}^{\vee})\cap M.
\end{equation}

Next we will prove our claim for $n_J=1, 2$ and $3$ respectively. 


Case 1: $n_J=1$. Recall that (see Remark \ref{psi}) a nef line bundle $\ul{L}$ in the interior of $\mathfrak{L}^J_{\bl{0}}$ defines a morphism $X\rightarrow Y=Y(\Sigma_{\ul{L}})$, which is induced by a surjective map from $N$ to $N'$. We will denote this map by $\psi_J$ and the subfan $\ker\psi_J\cap\Sigma_X\subset\Sigma_X$ by $\Sigma^J_X$. Then for any $\tau\in\Sigma_X(2)\cap \Sigma^J_X(2)$ we have 
 \[\|e_{\tau}(P_{L_1})\|\geq\|u_{\tau}\|,\]
where $u_{\tau}$ is a primitive vector paralleling with $e_{\tau}(P_{L_1})$. Since $u_{\tau}$ parallels with $u_{\tau'}$ for any $\tau'\in\Sigma_X(2)\cap\Sigma_X^J$, then $\|u_{\tau}\|=\|u_{\tau'}\|$.  As a consequence, by convexity for any $x\in P_{L_1}$ we have
\begin{equation}\label{flp1}\|l_{\tau}(x)\cap P_{L_1}\|\geq \|u_{\tau}\|,
\end{equation}
where $l_{\tau}(x)$ is the line passing throught $x$ and paralleling with $e_{\tau}(P_{L_1})$.

To proceed further, let $\Pi\subset M$ be a two-dimensional sublattice such that  $\mathds{Z}u_{\tau}\oplus\Pi=M$. Let $\pi: M_{\mathds{R}}\rightarrow \Pi_{\mathds{R}}$ be the projection with kernel $\mathds{R}u_{\tau}$. As $\Sigma_X$ is smooth, so is its subfan $\Sigma^J_X$. Then $\pi(P_{L_1})$ and $\pi(P_{L_1\otimes L_0}+\sigma^{\vee})$ are both lattice polygons in $\Pi_{\mathds{R}}$.

If $\pi(P_{L_1})$ contains only three lattice points, then $X$ must be isomorphic to $\mathbb{P}^1\times\mathbb{P}^2$ and the claim can be proved directly. Otherwise by the main results of \cite{fakhruddin2002multiplication} and \cite{haase2008lattice}, for any lattice point $y\in\pi(P_{L_1\otimes L_0}+\ub{\sigma}^{\vee})$, one can find an integral vector $v\in\Pi$ such that \[y\in\pi(P_{L_1})+v\subseteq \pi(P_{L_1\otimes L_0}+\ub{\sigma}^{\vee}).\] 
Note that one can always find an integral translation of $P_{L_1}$ contained in the following semi-infinite tubular region \[\pi^{-1}(\pi(P_{L_1})+v)\cap (P_{L_1\otimes L_0}+\ub{\sigma}^{\vee}).\] By abuse of notations, if we still denote such translation by $P_{L_1}$, then one has \[\pi^{-1}(y)\cap P_{L_1}\subset \pi^{-1}(y)\cap(P_{L_1\otimes L_0}+\ub{\sigma}^{\vee}).\]
By (\ref{flp1}) the segment $\pi^{-1}(y)\cap P_{L_1}$ contains a lattice point, thus the lattice points of  $\pi^{-1}(y)\cap(P_{L_1\otimes L_0}+\ub{\sigma}^{\vee})$ that are not contained in lattice translations of $P_{L_1}$ are located in certain bounded neighborhoods of the vertices of $P_{L_1\otimes L_0}+\ub{\sigma}^{\vee}$.


\medskip

Case 2: $n_J=2$. In this case we have $\dim\ker\psi_J=1$, let $\alpha\in\Sigma_X(1)$ such that $\Sigma_X(1)\cap\ker\psi_J=\{\alpha, -\alpha\}$. 

\bigskip
Let $H$ be a hyperplane lying between $H_{\alpha}(P_{L_1\otimes L_0})$ and $H_{-\alpha}(P_{L_1\otimes L_0})$ such that $H\cap M\neq\emptyset$. It is clear that there are only finitely many such hyperplanes. Next we will show for any such hyperplane $H$ the lattice points of $(P_{L_1\otimes L_0}+\ub{\sigma}^{\vee})\cap H$ missed out by integral translations of $P_{L_1}$ are contained in certain bounded neighorhoods of the vertices of $P_{L_1\otimes L_0}$. 

Following the notations above, let $\ub{\sigma}(1)=\{\beta, \beta'\}$ with $\beta, \beta'\in\Sigma_J(1)$ and
\[\psi_J^{-1}(\beta)\cap\Sigma_X(1)=\{\alpha, -\alpha, \alpha_1,\cdots,\alpha_n\}.\]
Let $\alpha_0=\alpha$ and $\alpha_{n+1}=-\alpha$, then for each $0\leq i\leq n$ the two-dimensional convex cone $\tau_i$ spanned by $\alpha_i, \alpha_{i+1}$ 
is an element of $\Sigma_X(2)$. In particular, $\psi_J(\tau_i)=\beta$. Then for the hyperplane $H$ aforementioned, one can find some $1\leq i\leq n$ such that $F_{\alpha_i}(P_{L_1\otimes L_0})\cap H$ is an edge of the  polygon $P_{L_1\otimes L_0}\cap H$. Denoting this edge by $e_{\beta}(P_{L_1\otimes L_0}\cap H)$ and let $u_{\beta}$ be a primitive vector paralleling with the $e_{\beta}(P_{L_1\otimes L_0}\cap H)$, then one sees easily $u_{\beta}$ also parallels with $e_{\tau_i}(P_{L_1\otimes L_0}), 0\leq i\leq n$. Note that $\|e_{\tau_{i-1}}(P_{L_1\otimes L_0})\|\geq\|u_{\beta}\|$, $\|e_{\tau_i}(P_{L_1\otimes L_0})\|\geq\|u_{\beta}\|$, then by convexity
\[\|e_{\beta}(P_{L_1\otimes L_0}\cap H)\|\geq\|u_{\beta}\|.
\]

Therefore, $e_{\beta}(P_{L_1\otimes L_0}\cap H)$ either contains at least one or none lattice point. It is easy to see in either case we can find an integral translation of $P_{L_1}$ with $F_{\alpha_i}(P_{L_1})$ contained in $F_{\alpha_i}(P_{L_1\otimes L_0})$ and moreover it has nonempty intersection with $H$. By abuse of notation, we just denote this translation by $P_{L_1}$, then $e_{\beta}(P_{L_1}\cap H)\subseteq e_{\beta}(P_{L_1\otimes L_0}\cap H)$, where $e_{\beta}(P_{L_1}\cap H)$ is an edge of $P_{L_1}\cap H$ defined similarly as $e_{\beta}(P_{L_1\otimes L_0}\cap H)$. 

If $e_{\beta}(P_{L_1\otimes L_0}\cap H)$ contains a lattice point then so does $e_{\beta}(P_{L_1}\cap H)$ as $\|e_{\beta}(P_{L_1}\cap H)\|\geq\|u_{\beta}\|$.

If $e_{\beta}(P_{L_1\otimes L_0}\cap H)$ contains no lattice points, we will show some lattice point can be found on some chord of $P_{L_1}\cap H$ paralleling with  $e_{\beta}(P_{L_1}\cap H)$. 

First we give our proof under the condition that both $F_{\alpha}(P_{L_1})$ and $F_{-\alpha}(P_{L_1})$ contain more than three lattice points.  Note that the edge $e_{\tau_0}(P_{L_1})$ (resp. $e_{\tau_n}(P_{L_1})$) is also an edge of $F_{\alpha}(P_{L_1})$ (resp. $F_{-\alpha}(P_{L_1})$) and it parallels with $u_{\beta}$. Let $\tilde{e}_{\beta}(F_{\alpha}(P_{L_1}))$ (resp. $\tilde{e}_{\beta}(F_{-\alpha}(P_{L_1}))$) be the chord of $F_{\alpha}(P_{L_1})$ (resp. $F_{-\alpha}(P_{L_1})$) containing the opposite side of $e_{\tau_0}(P_{L_1})$ (resp. $e_{\tau_n}(P_{L_1})$) in the unimodular parallelogram as described in Lemma \ref{upp}. 

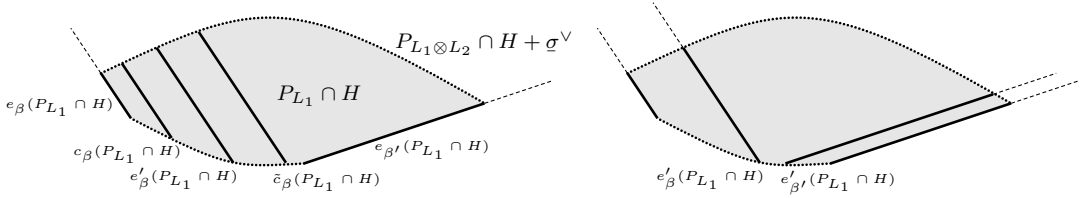
 \begin{figure}[H]
\begin{tikzpicture}
  \newcommand{\pathA}{(1.8,-0.4) .. controls (0.9,-0.45) and (0.8,-0.48) .. (-0.5,0.2)}
 \newcommand{\pathB}{(-0.9,0.8) .. controls (1.2,1.8) and (1.6,1.9) .. (4.2,0.4)}

\begin{scope}
 \fill[fill opacity=0.1] \pathB -- (1.8,-0.4) -- \pathA  -- (-0.9,0.8)--cycle;
\end{scope}

  \draw[line width=1pt, line cap=round, dash pattern=on 0pt off 1.5\pgflinewidth] \pathA;
  \draw[line width=1pt, line cap=round, dash pattern=on 0pt off 1.5\pgflinewidth] \pathB;

\draw[dash pattern=on 1.5pt off 1pt](-0.9,0.8)--(-1.3,1.4);

\draw[dash pattern=on 1.5pt off 1pt](4.2,0.4)--(5.1,0.7);
  
\draw[line width=1pt] (-0.5,0.2)--(-0.9,0.8);

\draw[line width=1pt] (1.56,-0.39)--(0.4,1.35);

\draw[line width=1pt] (0.86,-0.39)--(-0.16,1.14);

\draw [line width=1pt](0.04,-0.06)--(-0.62,0.93);

\draw [line width=1pt](1.8,-0.4)--(4.2,0.4);

  \newcommand{\pathC}{(8.8,-0.4) .. controls (7.9,-0.45) and (7.8,-0.48) .. (6.5,0.2)}
 \newcommand{\pathD}{(6.1,0.8) .. controls (8.2,1.8) and (8.6,1.9) .. (11.2,0.4)}

\begin{scope}
 \fill[fill opacity=0.1] \pathD -- (8.8,-0.4) -- \pathC  -- (6.1,0.8)--cycle;
\end{scope}

  \draw[line width=1pt, line cap=round, dash pattern=on 0pt off 1.5\pgflinewidth] \pathC;
  \draw[line width=1pt, line cap=round, dash pattern=on 0pt off 1.5\pgflinewidth] \pathD;

\draw[dash pattern=on 1.5pt off 1pt](6.1,0.8)--(5.7,1.4);

\draw[dash pattern=on 1.5pt off 1pt](6.84,1.14)--(6.44,1.74);

\draw[dash pattern=on 1.5pt off 1pt](11.2,0.4)--(12.1,0.7);

\draw[dash pattern=on 1.5pt off 1pt](10.96,0.52)--(11.8,0.8);

\draw[line width=1pt] (6.5,0.2)--(6.1,0.8);

\draw[line width=1pt] (7.86,-0.39)--(6.84,1.14);

\draw[line width=1pt] (8.2,-0.4)--(10.96,0.52);

\draw [line width=1pt](8.8,-0.4)--(11.2,0.4);

\node at(2,0.5) (0) {\relsize{-1}$P_{L_1}\cap H$}; 

\node at(4.2,1.2) (0) {\relsize{-1}$P_{L_1\otimes L_2}\cap H+\ub{\sigma}^{\vee}$}; 
\node at(2.1,-0.65) (0) {\relsize{-10}$\tilde{c}_{\beta}(P_{L_1}\cap H)$}; 

\node at(-0.55,-0.28) (0) {\relsize{-10}$c_{\beta}(P_{L_1}\cap H)$}; 

\node at(0.2,-0.6) (0) {\relsize{-10}$e'_{\beta}(P_{L_1}\cap H)$};

\node at(7.2,-0.6) (0) {\relsize{-10}$e'_{\beta}(P_{L_1}\cap H)$};

\node at(8.9,-0.65) (0) {\relsize{-10}$e'_{\beta'}(P_{L_1}\cap H)$};

\node at(-1.45,0.35) (0) {\relsize{-10}$e_{\beta}(P_{L_1}\cap H)$}; 

\node at(3.5,-0.22) (0) {\relsize{-10}$e_{\beta'}(P_{L_1}\cap H)$};

\end{tikzpicture}
\end{figure}

Let $c_{\beta}(P_{L_1}\cap H)$ (resp. $\tilde{c}_{\beta}(P_{L_1}\cap H)$) be the chord of $P_{L_1}\cap H$ that is contained in the intersection of $H$ with the hyperplane passing through $e_{\tau_0}(P_{L_1})$ and $e_{\tau_n}(P_{L_1})$ (resp. through $\tilde{e}_{\beta}(F_{\alpha}(P_{L_1}))$ and $\tilde{e}_{\beta}(F_{-\alpha}(P_{L_1}))$). Then one sees easily $\|c_{\beta}(P_{L_1}\cap H)\|\geq \|u_{\beta}\|$, $\|\tilde{c}_{\beta}(P_{L_1}\cap H)\|\geq \|u_{\beta}\|$ and the region bounded by these chords together with the border of $P_{L_1}\cap H$ contains a lattice point. 

Let $e'_{\beta}(P_{L_1}\cap H)$ be the chord of $P_{L_1}\cap H$ nearest to $e_{\beta}(P_{L_1}\cap H)$ such that it parallels with $u_{\beta}$ and contains a lattice point. Then this chord lies between $e_{\beta}(P_{L_1}\cap H)$ and $e'_{\beta}(P_{L_1}\cap H)$. Since the Euclidean lengths of all chords of $P_{L_1}\cap H$ with direction $u_{\beta}$ lying between $e_{\beta}(P_{L_1}\cap H)$ and $\tilde{c}_{\beta}(P_{L_1}\cap H)$ are no less than $\|u_{\beta}\|$, none of them or the line containing one of them has a lattice point by our definition of $e'_{\beta}(P_{L_1}\cap H)$.

By the same token, when  $e_{\beta'}(P_{L_1}\cap H)$ is lattice-free there exists a chord $e'_{\beta'}(P_{L_1}\cap H)$ containing  a lattice point and the open region bounded by the lines containing this chord and $e_{\beta'}(P_{L_1}\cap H)$ contains no lattice points. It's clear that our claim can be deduced from these results.

 It still remains to prove the case when $F_{\alpha}(P_{L_1})$ or $F_{-\alpha}(P_{L_1})$ contains exactly three lattice points. If $X$ fails to be be isomorphic to $\mathbb{P}^1\times \mathbb{P}^2$ and  $F_{\alpha}(P_{L_1})$ contains exactly three lattice points, then one sees easily the  lattice triangle $F'_{\alpha}(P_{L_1})\subset P_{L_1}$ paralleling with and nearest to $F_{\alpha}(P_{L_1})$ contains more lattice points. If we replace $F_{\alpha}(P_{L_1})$ by $F'_{\alpha}(P_{L_1})$ and $F_{-\alpha}(P_{L_1})$ by $F'_{-\alpha}(P_{L_1})$ if necessary, one sees easily the proof above still works.  Finally, the surjectivity of (\ref{ssgm}) can be proved directly when $X$ is isomorphic to $\mathbb{P}^1\times \mathbb{P}^2$.

Case 3: $n_J=3$. In this case by Proposition \ref{cas3}, we only need to show all but finitely many lattice points in a tubular neighborhood of the infinite edges of $P_{L_1\otimes L_2}+\ub{\sigma}^{\vee}$ with a fixed direction are contained in certain integral translations of $P_{L_1}$.

Let $U_{\tau}$ be a tubular neighborhood as described in Proposition \ref{cas3} such that the distance between points in it and the infinite edge $e_{\tau}(P_{L_1\otimes L_0}+\ub{\sigma}^{\vee})$ is bounded for some $\tau\in\ub{\sigma}(2)$.  Let $H_{\tau}=\{\xi\in N_{\mathds{R}}\:|\:\langle  u_{\tau}, \xi\rangle=0\}$, and $\alpha\in H_{\tau}\cap\Sigma_X(1)$ iff $F_{\alpha}(P_{L_1})$ contains an edge paralleling with $u_{\tau}$.


Since $\Sigma_X$ is smooth, we can take a $\rho\in N$ lying in the interior of $\ub{\sigma}$  such that \[\mathds{Z}\rho\oplus H_{\tau}\cap N=N.\] Now for each $a\geq 0$ let $H_{\rho}(a)$ be the hyperplane 
\[\{x\in M_{\mathds{R}}\:|\:\langle x, \rho\rangle=a\},
\] 
then it is easy to see for all sufficiently large integer $a_0$, the lattice points of $U_{\tau}$ that are not contained in the following set is finite
\[\bigcup_{a\geq a_0, a\in\mathds{Z}}H_{\rho}(a)\cap U_{\tau}.
\]
Besides, it is easy to see the intersections $H_{\rho}(a)\cap U_{\tau}$ for various $a\in\mathds{Z}$ are translations of each other by integral multiples of $u_{\tau}$, hence it suffices to prove for $a$ large enough any lattice point of $H_{\rho}(a)\cap U_{\tau}$ is contained in some translation of $P_{L_1}$. 
Note that $H_{\rho}(0)\cap M$ is a complement to $\mathds{Z}u_{\tau}$ in $M$. By abuse of notations, let $\pi: M_{\mathds{R}}\rightarrow H_{\rho}(0)$ be the projection with kernel $\mathds{R}u_{\tau}$, then $\pi(P_{L_1})$ is a lattice polygon on $H_{\rho}(0)$ which is contained in the lattice polygonal region $\pi(P_{L_1\otimes L_0}+\sigma^{\vee})$. Then by applying Lemma \ref{ccas3} any lattice point of $\pi(U_{\tau})\subset \pi(P_{L_1\otimes L_0}+\sigma^{\vee})$ is contained in some integral translation of $\pi(P_{L_1})$, from which our claim follows.

\end{proof}

\subsection{The Surjectivity of $\Phi_{L_1, L_2}$}\hfill



\begin{pro}\label{fe} For a given nef line bundle $L_2$, the map (\ref{sp}) is surjective for all but finitely many line bundles $L_1\in\amp(X)$.
\end{pro}

\begin{proof}Let $L_0+\mathfrak{L}^J_{\bl{b}}\subseteq\amp(X)$ be a reduced linear subset, by Lemma \ref{pfr} it suffices  to show the map $\Phi_{L_0\otimes\ul{L}, L_2}$ is surjective for all $\ul{L}\in \mathfrak{L}^J_{0}$ such that $\mathbb{E}(\ul{L})\gg0$. We first prove this claim for the case when $\dim\Sigma_{L_2}=3$ by showing the map below is surjective for all $\ul{L}$ lying sufficiently deep in $\mathfrak{L}^J_{\bl{0}}$
\begin{equation}\label{sm0}
(P_{L_0\otimes\ul{L}}\cap M)\times (\partial P_{L_2}\cap M)\rightarrow P_{L_0\otimes L_2\otimes\ul{L}}\cap M.
\end{equation}


Let  
$\Sigma_J$ be the normal fan of a polytope associated to a line bundle in the relative interior of $\mathfrak{L}^J_{\bl{0}}$ and $n_J=\dim\Sigma_J$. Take $\delta\in(\tfrac{3}{4},1)$ and denote by $\ul{L}^{\delta}\in\nef(X)_{\mathds{Q}}$ the line bundle whose associated polytope is $\delta P_{\ul{L}}$. 
Similar as (\ref{u2}) we have
\begin{equation}\label{sm1}P_{L_0\otimes L_2\otimes \ul{L}}=\bigcup_{\ub{\sigma}\in\Sigma_J(n_J)}P_{L_0\otimes L_2\otimes \ul{L}^{\delta}}(v_{\ub{\sigma}}(P_{L_0\otimes L_2\otimes\ul{L}})).
\end{equation}
Recall that (see (\ref{pv})) $\Sigma_J(n_J)$ defines a partition of $\Sigma_X(3)$, then since $\Sigma_{L_2}(1)\subseteq\Sigma_X(1)$ for any $\ub{\sigma}\in\Sigma_J(n_J)$ we can take 
\[S_{\ub{\sigma}}=\{\rho\in \Sigma_{L_2}(1)\;|\;\rho\preceq\sigma\;\text{for some}\;\sigma\in\Sigma_X(3)_{\ub{\sigma}}
\}.
\]Let $x_0\in\partial P_{L_2}\cap M$ and $P_{L_0\otimes \ul{L}}(x_0)$ be the integral translation of $P_{L_0\otimes \ul{L}}$ contained in $P_{L_0\otimes L_2\ul{L}}$ corrsponding to $x_0$, then by combining (\ref{sm0}) and (\ref{sm1}) it suffices to show for any $\ub{\sigma}\in\Sigma_J(n_J)$
\begin{equation}\label{ttl331}P_{L_0\otimes L_2\otimes \ul{L}^{\delta}}(v_{\sigma}(P_{L_0\otimes L_2\otimes \ul{L}}))\cap M\subseteq\bigcup_{\substack{\rho\in S_{\ub{\sigma}}}
} \bigcup_{x\in F_{\rho}(P_{L_2})\cap M}(\overrightarrow{x_0x}+P_{L_0\otimes \ul{L}}(x_0))\cap M.
\end{equation}
The proof of  (\ref{ttl331}) given below follows a similar line as our argument in Claim \ref{qcl3}. Let 
\[P_{L_0\otimes \ul{L}, \ub{\sigma}}(x_0)=P_{L_0\otimes \ul{L}}(x_0)+\ub{\sigma}^{\vee},\;\;P_{L_0\otimes L_2\otimes \ul{L}, \ub{\sigma}}=P_{L_0\otimes L_2\otimes \ul{L}}+\ub{\sigma}^{\vee},\]
and 
\[S'_{\ub{\sigma}}=\{\rho\in S_{\ub{\sigma}}\:|\: \text{the facet of}\; P_{L_0\otimes \ul{L}, \ub{\sigma}}\;\text{corresponding to}\; \rho \;\text{is bounded}\}.\] Then it is easy to see (\ref{ttl331}) follows from
\begin{equation}\label{ttl332}P_{L_0\otimes L_2\otimes \ul{L},\ub{\sigma}} \cap M\subseteq \bigcup_{\rho\in S'_{\ub{\sigma}}} \bigcup_{x\in F_{\rho}(P_{L_2})\cap M}(\overrightarrow{x_0x}+P_{L_0\otimes \ul{L}, \ub{\sigma}}(x_0))\cap M
\end{equation}
and for any $x\in F_{\rho}(P_{L_2})\cap M, \rho\in S'_{\ub{\sigma}}$,
\begin{equation}\label{ttl333}P_{L_0\otimes L_2\otimes \ul{L}^{\delta}}(v_{\ub{\sigma}}(P_{L_0\otimes L_2\otimes\ul{L}}))\cap\big((\overrightarrow{x_0x}+P_{L_0\otimes \ul{L}, \ub{\sigma}}(x_0))\backslash (\overrightarrow{x_0x}+P_{L_0\otimes \ul{L}}(x_0)\big)=\emptyset.\end{equation}

\bigskip
We first prove (\ref{ttl332}). By Lemma \ref{smw}, we only need to prove  for all $\rho\in  S'_{\ub{\sigma}}$ the following equality
\begin{equation}\label{ttl334}(F_{\rho}(P_{L_0\otimes L_2\otimes \ul{L}})+\ub{\sigma}^{\vee})\cap M=\bigcup_{x\in F_{\rho}(P_{L_2})\cap M}(\overrightarrow{x_0x}+F_{\rho}(P_{L_0\otimes \ul{L}}(x_0))+\ub{\sigma}^{\vee})\cap M
\end{equation}
If  $D_{\rho}$ is not isomorphic to $\mathbb{P}^2$ or if $D_{\rho}$ is isomorphic to $\mathbb{P}^2$ and $F_{\rho}(P_{L_0\otimes \ul{L}}(x_0))$ contains more than three lattice points the equality above can be deduced from Theorem \ref{sfhn}. Now we consider the case $D_{\rho}$ is isomorphic to $\mathbb{P}^2$ and $F_{\rho}(P_{L_0\otimes \ul{L}}(x_0))$ contains only three lattice points.
If $\#(F_{\rho}(P_{L_0\otimes \ul{L}\otimes\mathcal{O}_X(-D_{\rho})})\cap M)>3$, we can still apply Theorem \ref{sfhn} and get
\begin{equation}\label{ttl335}(F_{\rho}(P_{L_0\otimes L_2\otimes \ul{L}\otimes\mathcal{O}_X(-D_{\rho})})+\ub{\sigma}^{\vee})\cap M=\bigcup_{x\in F_{\rho}(P_{L_2})\cap M}(\overrightarrow{x_0x}+F_{\rho}(P_{L_0\otimes \ul{L}\otimes\mathcal{O}_X(-D_{\rho})})+\ub{\sigma}^{\vee})\cap M.
\end{equation}
As a consequence (\ref{ttl334}) can be obtained. If $\#(F_{\rho}(P_{L_0\otimes \ul{L}\otimes\mathcal{O}_X(-D_{\rho})})\cap M)=3$, we will have $n_J=\dim\ub{\sigma}=1$, from which we can deduce $X\cong \mathbb{P}^2\times\mathbb{P}^1$ and the conclusion can be proved directly. 

To prove (\ref{ttl333}), we first observe the normal fans of $P_{L_0\otimes \ul{L},\ub{\sigma}}(x_0)$ and
$P_{L_2\otimes L_0\otimes \ul{L},\ub{\sigma}}$ are the same and let it be $\Sigma_{\ub{\sigma}}$, then $\Sigma_{\ub{\sigma}}(1)\subset \Sigma_X(1)$. For each $\rho\in\Sigma_X(1)\backslash\Sigma_{\ub{\sigma}}(1)$, let $H^-_{\rho}(P_{L_0\otimes \ul{L}}(x_0))$ be the open half space such that it and $P_{L_0\otimes \ul{L}}(x_0)$ lie on different sides of $F_{\rho}(P_{L_0\otimes \ul{L}}(v_{\sigma}))$, then $(\overrightarrow{x_0x}+P_{L_0\otimes \ul{L}, \ub{\sigma}}(x_0))\backslash (\overrightarrow{x_0x}+P_{L_0\otimes \ul{L}}(x_0))$ is contained in 
\[\bigcup_{\rho\in\Sigma_X( 1)\backslash\Sigma_{\ub{\sigma}}(1)}(\overrightarrow{x_0x}+H^-_{\rho}(P_{L_0\otimes \ul{L}}(x_0)))\cap(\overrightarrow{x_0x}+P_{L_0\otimes \ul{L}, \ub{\sigma}}(x_0)).
\]
It suffices to see  for any $x\in F_{\rho}(P_{L_2})\cap M, \rho\in \Sigma_X(1)\backslash\Sigma_{\ub{\sigma}}(1)$ we have 
\begin{equation}\label{shna}P_{L_0\otimes L_2\otimes \ul{L}^{\delta}}(v_{\ub{\sigma}}(P_{L_0\otimes L_2}))\cap (\overrightarrow{x_0x}+H^-_{\rho}(P_{L_0\otimes \ul{L}}(x_0)))=\emptyset\end{equation}
when $\mathbb{E}(\ul{L})\gg0$, hence (\ref{sm0}) is valid when $\dim\Sigma_{L_2}=3$. The proofs for the cases when $\dim\Sigma_{L_2}=1$ or $2$ are similar as above except that we will replace $\partial P_{L_2}$ by $P_{L_2}$ in (\ref{sm0}), then the proof for the proposition is complete.

\end{proof}
\begin{proof}[Proof of Theorem \ref{t2}]First of all, by the proof of Proposition \ref{p3} we can write $\nef(X)$ as a union of finitely many reduced linear subsets $\{\mathfrak{L}^{J_k}_{\bl{b}_k}\}_{1\leq k\leq s}$ such that 
\begin{equation}\label{ijb}L\prec L'\Rightarrow L\prec_c L'
\end{equation}
 for any $L, L'\in\mathfrak{L}^{J_k}_{\bl{b}_k}$ (\ref{sd5}). 
Let $\mathfrak{L}^{J_k}_{\bl{b}_k}=L_k+\mathfrak{L}^{J_k}_{\bl{0}}$, $\{B_j\}_{1\leq j\leq m}$ be a Hilbert basis of $\mathfrak{L}^{J_k}_{\bl{0}}$ and 
take
\[\mathfrak{B}^{J_k}_{\bl{b}_k}=\{\bigotimes_j B_j^{i_j}\;|\; 0\leq i_j\leq 3\}.
\]
Next we will prove the following 
\begin{claim}\label{lab} For given $\tilde{L}\in\amp(X)$,  if the map $\Phi_{\tilde{L}, L_k\otimes \ul{L}}$ is surjective for all $\ul{L}\in \mathfrak{B}^{J_k}_{\bl{b}_k}$, then $\Phi_{\tilde{L}, \bar{L}}$ is surjective for all $\bar{L}\in L_k+\mathfrak{L}^{J_k}_{\bl{0}}$.

\end{claim} Indeed, $\bar{L}$ can be written  as 
\[\bar{L}=\bigotimes_jB_j^{k_j}\otimes \ul{L}\otimes L_k,
\]
for some $\ul{L}\in\mathfrak{B}^J_{\bl{b}}$ such that the power $k_j\neq 0$ iff the power of $B_j$ in $\ul{L}$ is equal to 3. The claim will be proved by induction on the powers $k_j, 1\leq j\leq m$. Firstly, if $k_j=0$ for all $j$, then we are done by assumption. Now assume the map $\Phi_{\tilde{L}, \bar{L}}$ is surjective, we will prove so is true for $\Phi_{\tilde{L}, \bar{L}\otimes B_j}$, where $j$ satisfies the condition that the power of $B_j$ in $\ul{L}$ is 3.
Given an invariant curve $C_{\tau}$, if $B_j.C_{\tau}\neq0$ then we have
\[\frac{(\tilde{L}\otimes \bar{L}). C_{\tau}}{(\tilde{L}\otimes \bar{L}\otimes B_j). C_{\tau}}>\frac{B^{k_j+3}_j.C_{\tau}}{B^{k_j+4}_j.C_{\tau}}\geq\tfrac{3}{4}.\]
If $B_j.C_{\tau}=0$, then since $\tilde{L}$ is ample we have
\[\frac{(\tilde{L}\otimes \bar{L}). C_{\tau}}{(\tilde{L}\otimes \bar{L}\otimes B_j). C_{\tau}}=1.\]
By the proof of Proposition \ref{p1}, the following map is surjective 
\[(P_{\tilde{L}\otimes \bar{L}}\cap M)\times (P_{B_j}\cap M)\rightarrow P_{\tilde{L}\otimes \bar{L}\otimes B_j}\cap M.
\]
On the other hand, by induction the following map is also surjective
\[(P_{\tilde{L}}\cap M)\times (P_{\bar{L}}\cap M)\rightarrow P_{\tilde{L}\otimes \bar{L}}\cap M,
\]
hence so is the following map
\[(P_{\tilde{L}}\cap M)\times (P_{\bar{L}}\cap M)\times (P_{B_j}\cap M)\rightarrow P_{\tilde{L}\otimes \bar{L}\otimes B_j}\cap M.
\]
Moreover, since $\bar{L}\prec \bar{L}\otimes B_j$ and both of them are contained in $\mathfrak{L}^{J_k}_{\bl{b}_k}$, by (\ref{ijb}) the map $\Phi_{\bar{L}, B_j}$ is surjective. As the following diagram is commutative, 
\[\xymatrix{
(P_{\tilde{L}}\cap M)\times (P_{\bar{L}}\cap M)\times (P_{B_j}\cap M) \ar[d] \ar[r] &P_{\tilde{L}\otimes \bar{L}\otimes B_j}\cap M\\
(P_{\tilde{L}}\cap M)\times (P_{\bar{L}\otimes B_j}\cap M) \ar[ur] &}
\]
the map $\Phi_{\tilde{L}, \bar{L}\otimes B_j}$ is surjective. Then by induction Claim \ref{lab} is proved.

Now by Proposition \ref{fe}, to each line bundle $\hat{L}$ from the following finite set
\[\mathfrak{S}_1=\{L_k\otimes \ul{L}\:|\:  \ul{L}\in\mathfrak{B}^{J_k}_{\bl{b}_k}, 1\leq k\leq s\},\]
we can assign a finite subset $\mathfrak{S}_{\hat{L}}$ of $\amp(X)$ such that the map $\Phi_{\tilde{L}, \hat{L}}$ is surjective as long as $\tilde{L}$ is contained in $\amp(X)\backslash \mathfrak{S}_{\hat{L}}$. Let $\mathfrak{S}_2=\bigcup_{\hat{L}\in\mathfrak{S}_1}\mathfrak{S}_{\hat{L}}$, then $\Phi_{\tilde{L}, \hat{L}}$ is surjective for any $\tilde{L}\in\amp(X)\backslash\mathfrak{S}_2$ and $\hat{L}\in\mathfrak{S}_1$. By Claim \ref{lab}, the map $\Phi_{\tilde{L}, L}$ is surjective for $\tilde{L}\in\amp(X)\backslash\mathfrak{S}_2$ and any nef line bundle $L$.

\end{proof}

\begin{rem} If each codimension one facet of an ample line bundle on $X$ has at least four edges, then by a modification of the proof above, one can show $L_1\prec_c L_1\otimes L_2$ for any ample line bunlde $L_1$ not in a certain finite subset of $\amp(X)$ and any nef line bunlde $L_2$. 
\end{rem}

\begin{cor} \label{palb}With at most finitely many exceptions,  any ample line bundle on a given smooth toric threefold is  projectively normal. 

\end{cor}

\section{Explicit Estimations related to Pairs Potentially Violating the Surjectivity of (\ref{sp})}\label{q2v}\hfill

In the first part of this section we reexamine the proof of Lemma \ref{pfr} in more detail
and formalize the dimension lowering process, which will act as the cornerstone of our estimation. In the second part, we provide some explicit bounds for the pairs of line bundles that fall outside some of the positive results obtained in previous sections.

\subsection{The Dimension Lowering Process for Linear Subsets}\hfill

The proof of all our finiteness results obtained as far are based on Lemma \ref{pfr}, where the number of the line bundles that might violate the property $\mathscr{P}$ is implicit. Recall that to apply Lemma \ref{pfr}, we need to find for a reduced linear subset $\mathfrak{L}^J_{\bl{b}}$ a new integer vector $\bl{b}'$ with $\bl{b}'>\bl{b}$ such that all line bundles in $\mathfrak{L}^J_{\bl{b}'}$ satisfy the property $\mathscr{P}$. Then we come to
\begin{equation}\label{rn}\mathfrak{L}^J_{\bl{b}}\backslash \mathfrak{L}^J_{\bl{b}'}=\bigcup_{\bar{J}, \bar{\bl{b}}}\mathfrak{L}^{\bar{J}}_{\bar{\bl{b}}},
\end{equation}
where each of the finitely many linear subsets $\mathfrak{L}^{\bar{J}}_{\bar{\bl{b}}}$ on the right is reduced and moreover $\dim\mathfrak{L}^{\bar{J}}_{\bl{0}}<\dim\mathfrak{L}^J_{\bl{0}}$. This dimension lowering process for the linear subsets to acquire the property $\mathscr{P}$ will continue until $\dim\mathfrak{L}^{\bar{J}}_{\bl{0}}=0$, exactly when we arrive at the line bundles that might fail to satisfy $\mathscr{P}$. To extract effective bounds for the potential counter-examples, we need to make this process tractable by finding an explicit  upper bound for each $\bar{\bl{b}}$. 

Obviously, the way to write $\mathfrak{L}^J_{\bl{b}}\backslash \mathfrak{L}^J_{\bl{b}'}$ as union of linear subsets with lower dimensions is not unique.
For our purpose we will require each component on the right of (\ref{rn}) is a maximal reduced linear subset. Then the index $\bar{J}$ must be minimal, i.e. there exists no other reduced linear subset $\mathfrak{L}^{\bar{J}'}_{\bar{\bl{b}}'}\subset \mathfrak{L}^J_{\bl{b}}\backslash \mathfrak{L}^J_{\bl{b}'}$ such that $\bar{J}'\subsetneq\bar{J}$ and $\mathfrak{L}^{\bar{J}}_{\bar{\bl{b}}}\subsetneq\mathfrak{L}^{\bar{J}'}_{\bar{\bl{b}}'}$. 

\begin{lem}\label{bbb}For each $\bar{J}\subset I$ in (\ref{rn})
there exists $\bar{\bl{b}}^{\bar{J}}\in\mathds{Z}^I_{\geq0}$ such that for  any maximal reduced linear subset $\mathfrak{L}^{\bar{J}}_{\bar{\bl{b}}}=\bar{L}+\mathfrak{L}^{\bar{J}}_{\bl{0}}$ contained in $\mathfrak{L}^J_{\bl{b}}\backslash \mathfrak{L}^J_{\bl{b}'}$ we have
$\bar{L}.C_{\tau}\leq \bar{b}^{\bar{J}}_{\tau}$ for all $\tau\in I$. 
\end{lem}
\begin{proof} Let $\mathfrak{B}=\{B_j\}_{1\leq j\leq m}$ be a Hilbert basis of $\nef(X)$, for any $S\subset I$ we define
\begin{align}
\mathfrak{B}^0_S&=\{B_j\in\mathfrak{B}\:|\: B_j.C_{\tau}=0\;\text{for all}\;\tau\in S\},\notag\\
\mathfrak{B}^+_S&=\{B_j\in\mathfrak{B}\:|\: B_j.C_{\tau}\neq 0\;\text{for some}\;\tau\in S\}.\notag\
\end{align}
Given $\bar{J}\subset I$ as on the right in (\ref{rn}), we label the line bundles in $\mathfrak{B}$ as follows
\begin{align}
\mathfrak{B}^+_J&=\{B_j\}_{1\leq j\leq k},\label{itau}\\
\mathfrak{B}^0_J\cap\mathfrak{B}^+_{\bar{J}\backslash J}&=\{B_j\}_{k+1\leq j\leq l},\label{jtau}\\
\mathfrak{B}^0_{\bar{J}}=\mathfrak{B}^0_{\bar{J}}\cap\mathfrak{B}^+_{I\backslash\bar{J}}&=\{B_j\}_{l+1\leq j\leq m}\label{ltau},
\end{align}
from which one sees easily
\begin{equation}\label{dxg}J=\{\tau\in I\;|\; B_j.C_{\tau}=0 \;\text{for any}\; k+1\leq j\leq m\}. 
\end{equation}
Let
\[(I\backslash J)_{\bar{J}, \mathfrak{B}}=\{\tau\in I\backslash J\;|\; B_j.C_{\tau}>0\;\text{for some}\; k+1\leq j\leq l\;\text{and}\; B_j.C_{\tau}=0 \;\text{for any}\; l+1\leq j\leq m\}, 
\]
then we claim \begin{equation}\label{dyg}(I\backslash J)_{\bar{J}, \mathfrak{B}}=\bar{J}\backslash J.\end{equation}
Firstly, for $\tau\in (I\backslash J)_{\bar{J}, \mathfrak{B}}$ we have $B_j.C_{\tau}$ for any $l+1\leq j\leq m$ hence by (\ref{ltau}) $\tau\in\bar{J}$ and $\bar{J}\backslash J\supseteq (I\backslash J)_{\bar{J}, \mathfrak{B}}$. Next we show $\bar{J}\backslash J\subseteq (I\backslash J)_{\bar{J}, \mathfrak{B}}$. Suppose $\tau\in\bar{J}\backslash J$ satisfies \begin{equation}\label{asmp}B_j.C_{\tau}=0, \;\;k+1\leq j\leq l,\end{equation} 
then there must be $1\leq j\leq k$ such that $B_j.C_{\tau}> 0$. On the other hand,  take $\hat{L}\in\mathfrak{L}^J_{\bl{b}}=L+\mathfrak{L}^J_{\bl{0}}$ and let
\begin{equation}\label{lhat}\hat{L}=L\otimes \bigotimes_{1\leq j\leq m}B_j^{n_j}, \;\;n_j\geq 0,
\end{equation}
by the definition of a linear subset and (\ref{dxg}) we get
\begin{equation}\label{ljbf}\sum_{1\leq j\leq k}B^{n_j}_j.C_{\tau'}=b_{\tau'}-L.C_{\tau'}=0,\;\;\tau'\in J.
\end{equation}
Then we can deduce $n_j=0$ for any $1\leq j\leq k$ since $B_j.C_{\tau'}\neq0$ for some $\tau'\in J$ by our definition (\ref{itau}).  However, combining this with (\ref{asmp}) one would get $\hat{L}.C_{\tau}=L.C_{\tau}$ for all $\hat{L}\in\mathfrak{L}^J_{\bl{b}}$ hence $\tau\in J$ as $\mathfrak{L}^J_{\bl{b}}$ is reduced, which contradicts with our assumption that $\tau\in \bar{J}\backslash J$. 

Next we will find minimal elements of the following set with respect to $\prec_{\bar{J}}$, which is defined by the subcone $\mathfrak{L}^{\bar{J}}_{\bl{0}}$ of $\nef(X)$ (see Section 2.1 for  the definition)
\[\mathfrak{S}^{\bar{J}}_{\bl{b}\backslash\bl{b}'}=\{\hat{L}\in\mathfrak{L}^J_{\bl{b}}\;|\; (\hat{L}+\mathfrak{L}^{\bar{J}}_{\bl{0}})\cap\mathfrak{L}^J_{\bl{b}'}=\emptyset\}. 
\]
Following the notations above (\ref{lhat}), since $\hat{L}\in\mathfrak{L}^J_{\bl{b}}$ by the argument above we have $n_j=0$ for all $1\leq j\leq k$. On the other hand, let $\hat{L}\in\mathfrak{S}^{\bar{J}}_{\bl{b}\backslash\bl{b}'}$, if $n_j>0$ for some $l+1\leq j\leq m$, then it is easy to see $\hat{L}\otimes B_j^{-n_j}\in\mathfrak{S}^{\bar{J}}_{\bl{b}\backslash\bl{b}'}$. Hence by the minimality requirement we can just assume $n_j=0$ for all $l+1\leq j\leq m$ in (\ref{lhat}).  Then for $\hat{L}\in\mathfrak{S}^{\bar{J}}_{\bl{b}\backslash\bl{b}'}$ (\ref{lhat}) can be further written as
\begin{equation}\label{lhatt}\hat{L}=L\otimes\bigotimes_{k+1\leq j\leq l}B_j^{n_j}.
\end{equation}


Note that the property $\mathscr{P}$ is satisfied by all line bundles in $\mathfrak{L}^J_{\bl{b}'}$, hence for $\hat{L}$ in the form of (\ref{lhatt}) to stay outside $\mathfrak{L}^J_{\bl{b}'}$, the integer vector $(n_j)_{ k+1\leq j\leq l}$ needs to satisfy the following inequality (recall that $\bl{b}=(b_{\tau})_{\tau\in I}$ and $\bl{b}'=(b'_{\tau})_{\tau\in I}$) for at least one $\tau_1\in \bar{J}\backslash J$ (see also (\ref{ddp}))
\begin{equation}\label{ine3}\sum_{k+1\leq j\leq l}n_jB_j.C_{\tau_1}\leq b'_{\tau_1}-b_{\tau_1}.
\end{equation}

Take
\[\{\tau_1\}^{(1)}=\{\tau\in I\backslash J\:|\: B_j.C_{\tau}=0\;\text{for all}\; B_j\in \mathfrak{B}^0_{\{\tau_1\}}\},\]
and 
\begin{equation}\{\tau_1\}^{(k+1)}=\{\tau\in I\backslash J\:|\: B_j.C_{\tau}=0\;\text{for all}\; B_j\in \mathfrak{B}^0_{\{\tau_1\}^{(k)}}\},\;\; k\geq 1.\end{equation}
One sees easily there are only finitely many different terms in the non-decreasing sequence $\{\{\tau_1\}^{(k)}\}_{k\geq1}$ and we will denote by $\overline{\{\tau_1\}}$ the largest one, then by definition we have 
\begin{equation}\label{tyb}\overline{\{\tau_1\}}=\{\tau\in I\backslash J\:|\: B_j.C_{\tau}=0\;\text{for all}\; B_j\in \mathfrak{B}^0_{\overline{\{\tau_1\}}}\}.
\end{equation}We claim $\overline{\{\tau_1\}}=\bar{J}\backslash J$. Firstly, by definition for any $l+1\leq j\leq m$ we have $B_j.C_{\tau_1}=0$ hence also $B_j.C_{\tau}=0$  for any $\tau\in\overline{\{\tau_1\}}$. Let $\tau\in\overline{\{\tau_1\}}$, then since $\tau\notin J$ by (\ref{dxg}) we can find $k+1\leq j\leq l$ such that $B_j.C_{\tau}>0$, then by (\ref{dyg}) we get $\tau\in\bar{J}$ hence $\overline{\{\tau_1\}}\subseteq\bar{J}\backslash J$. On the other hand, by (\ref{tyb}) for any $\tau_2\in (I\backslash J)\backslash\overline{\{\tau_1\}}$, there exists $B_j\in {B}^0_{\overline{\{\tau_1\}}}$ such that $B_j.C_{\tau_2}>0$. Therefore, by definition \ref{lss} the index set $J\cup\overline{\{\tau_1\}}$ is reduced, hence $\bar{J}=J\cup\overline{\{\tau_1\}}$ as $\bar{J}$ is minimal by our assumption.

Therefore, by our definition (\ref{jtau}) for any $k+1\leq j\leq l$ we have 
\[B_j.C_{\tau}>0\] 
for some $\tau\in \overline{\{\tau_1\}}$.
Then similar as (\ref{ine3}) for each $k+1\leq j\leq l$ we have
\begin{equation}\label{ine4}n_j\leq b'_{\tau}-b_{\tau}
\end{equation}
for some $\tau\in\overline{\{\tau_1\}}$.
Now by (\ref{lhatt}) and (\ref{ine4}) for a minimal $\hat{L}$ and $\tau\in I\backslash J$ we get
\begin{align}
\hat{L}.C_{\tau}&=(L\otimes \bigotimes_{k+1\leq j\leq l}B_j^{n_j}).C_{\tau}\notag\\
&=b_{\tau}+\sum_{k+1\leq j\leq l}n_jB_j.C_{\tau}\notag\\
&\leq b_{\tau}+a_{\tau}\sum_{\tau\in\bar{J}\backslash J}b'_{\tau}-b_{\tau},\notag
\end{align}
where $a_{\tau}=\max_{B_j\in\mathfrak{B}}B_j.C_{\tau}$. Then it is easy to see the vector $\bar{\bl{b}}_{\bar{J}}$ given by the formula below satisfies the condition in Claim \ref{bbb}
\begin{equation}\label{rf1}\bar{b}^{\bar{J}}_{\tau}=\begin{cases}b_{\tau}&\text{if}\;\tau\in J;\\
b_{\tau}+a_{\tau}\sum_{\tau\in\bar{J}\backslash J}b'_{\tau}-b_{\tau}&\text{otherwise}.\end{cases}
\end{equation}

\end{proof}
Now let 
\[\mathscr{J}(X)=\{J\subseteq I\:|\: \mathfrak{L}^J_{\bl{0}}\;\text{is a reduced subcone of}\; \nef(X)\},\]
and 
$\emptyset=J_0\subsetneq J_1\subsetneq J_2\subsetneq\cdots\subsetneq J_s\subsetneq J_{s+1}=I$ be a maximal chain consisting of elements from $\mathscr{J}$ in the dimension lowering process for detecting line bundles satisfying property $\mathscr{P}$. That is, there exist sequences of vectors $\{\bl{b}^{J_i}\}_{1\leq i\leq s+1}$ and $\{\bl{b}'^{J_i}\}_{1\leq i\leq s+1}$ such that all line bundles in $\mathfrak{L}^{J_i}_{\bl{b}'^{J_i}}$ satisfy property $\mathscr{P}$ and for each $0\leq i\leq s$, there are linear subsets of the form $\hat{L}+\mathfrak{L}^{J_{i+1}}_{\bl{0}}$ contained in $\mathfrak{L}^{J_i}_{\bl{b}^{J_i}}\backslash\mathfrak{L}^{J_i}_{\bl{b}'^{J_i}}$. Then by (\ref{rf1}) we have
\begin{equation}\label{rf2}\hat{L}.C_{\tau}\leq b^{J_{i+1}}_{\tau}=\begin{cases}b^{J_i}_{\tau}&\text{if}\;\tau\in J_i;\\
b^{J_i}_{\tau}+a_{\tau}\sum_{\tau\in J_{i+1}\backslash J_i}b'^{J_i}_{\tau}-b^{J_i}_{\tau}&\text{otherwise}.\end{cases}
\end{equation}

With these recursive relations one can derive the inequalities that the line bundles left  in the zero-dimensional linear subsets need to satisfy.

\subsection{The Estimation}\hfill

For each $J\in\mathscr{J}(X)$ let $\mathfrak{A}_J$ be a Hilbert basis of $\mathfrak{L}^J_{\bl{0}}\cap\nef(X)$ and take
\[c=\max_{\tau\in I}\max_{\substack{A\in\mathfrak{A}_J,\\J\in\mathscr{J}(X)}}A.C_{\tau}. 
\]
When $J=\emptyset$, we will just take $\mathfrak{A}_{\emptyset}=\mathfrak{B}$, then one sees easily $a_{\tau}\leq c$ for all $\tau\in I$. Moreover, by \cite[Theorem 2.12]{bg2009} we can find a finite set of ample line bundles $\{L^{(i)}\}_i$ such that
\begin{equation}\label{dac}\amp(X)=\bigcup_i (L^{(i)}+\nef(X))\cap\nef(X).
\end{equation}
For each $\tau\in I$ we will take 
\begin{equation}\label{ctau}c_{\tau}=\max_iL^{(i)}.C_{\tau}.
\end{equation}

\begin{pro}\label{loqr} The map (\ref{sp}) is surjective for any ample line bundle $L_1$ and nef line bundle $L_2$ on a projeictive toric variety with dimension $d$ if so is true when $L_1, L_2$ satisfy $L_i.C_{\tau}\leq\#Idc^2$ for any $\tau\in I$, $i=1, 2$.

\end{pro}
\begin{proof}
For any nef line bundle $L$, let $\mathscr{P}$ be the property that $L\prec\tilde{L}\Rightarrow L\prec_c\tilde{L}$. Then by our proofs of Proposition \ref{p1} and Proposition \ref{p3}, for any reduced linear subset $\mathfrak{L}^J_{\bl{b}}$, the linear subset $\mathfrak{L}^J_{\bl{b}'}$ will satisfy the condition of Lemma \ref{pfr} if 

\begin{equation}b'_{\tau}=\begin{cases}b_{\tau}&\text{if}\;\tau\in J;\\
b_{\tau}+\tfrac{d+\epsilon}{1-\epsilon}c&\text{otherwise},\end{cases}
\end{equation}
where $\epsilon$ is any positive number. Then by our argument in (\ref{rf2})
(note that $\mathfrak{L}^{J_0}_{\bl{b}_0}=\mathfrak{L}^{\emptyset}_{\bl{0}}$) for any reduced linear subset $\mathfrak{L}^J_{\bl{b}}=L+\mathfrak{L}^J_{\bl{b}}$ turning up in this dimension lowering process we have 
 \[L.C_{\tau}=b_{\tau}\leq \#Idc^2.
\]

\end{proof}

For now on till the end of this section, we denote by $X$ a smooth projective toric threefold. For $\rho\in\Sigma_X(1)$
we take
\[r_{\rho}=\min_{\substack{\tau\in\Sigma_X(2),\\\langle\rho,u_{\tau}\rangle\neq0}}|\langle u_{\tau},\rho\rangle|.
\]

\begin{pro}\label{lopr} Let $L_1$ be an ample and $L_2$ a nef line bundle on $X$, if
 \[L_1.C_{\tau}>c_{\tau}+4\#I\tfrac{c^2w_{\rho}(L_2)}{r_{\rho}}\] for some $\tau\in I$, where $w_{\rho}(L_2)=\max_{x, y\in P_{L_2}\cap M}|\langle \overrightarrow{xy},\rho\rangle|$, then we have $\coker\:\Phi_{L_1, L_2}=0$.
 
\end{pro}

\begin{proof} For given $L_2$ let $\mathscr{P}$ be the property that $\coker\:\Phi_{L_1, L_2}=0$ for $L_1\in\amp(X)$. Similar as the proof of Proposition \ref{loqr}, it suffices to find an explicit $\bl{b}'$ for a reduced linear subset $\mathfrak{L}^J_{\bl{b}}$ in $\amp(X)$. Following the notations in the proof of Proposition \ref{fe} for any $\rho\in\Sigma_X(1)\backslash\Sigma_{\ub{\sigma}}(1)$ we can find $\tau\in I\backslash J$ such that $\langle u_{\tau},\rho\rangle\neq 0$. By our proof there for (\ref{shna}) to be valid we only need
\[(1-\delta)\ul{L}.C_{\tau}\tfrac{r_{\rho}}{\|\rho\|}-\tfrac{w_{\rho}(L_2)}{\|\rho\|}\geq0,
\]
where $\delta\in(\tfrac{3}{4}, 1)$. Then we can take
\begin{equation}\notag b'_{\tau}=\begin{cases}b_{\tau}&\text{if}\;\tau\in J;\\
b_{\tau}+\tfrac{1}{1-\delta}\tfrac{cw_{\rho}(L_2)}{r_{\rho}}&\text{otherwise}.\end{cases}
\end{equation}
To apply (\ref{rf2}) we need to start from a reduced subset $L^{(i)}+\nef(X)$ (see (\ref{dac})) hence accordingly \[(b_{J_0})_{\tau}=L^{(i)}.C_{\tau}\leq c_{\tau}.\]
Note that the property $\mathscr{P}$ is not satisfied by an ample line bundle $L_1$ only if $L_1$ is contained in the linear subsets with zero dimension. That is, \[L_1.C_{\tau}\leq c_{\tau}+\tfrac{1}{1-\delta}\#I\tfrac{c^2w_{\rho}(L_2)}{r_{\rho}}\]
for any $\tau\in I$ and $\delta\in(\tfrac{3}{4},1)$, from which we can deduce the conclusion. 

\end{proof}

\begin{rem} An estimation for the size the finite subset of $\amp(X)$ in Theorem \ref{t2} can also be obtained by combining the results of the propositions above. As for Proposition \ref{b2}, to explicite a uniform bound of $\dim\coker\:\Phi_{L_1, L_2}$ would be challenging even for those $L_2=L_3^k$ with $L_3$ a given nef line bundle and $k\geq1$. On the other hand, in this setting the surjecitivity of $\Phi_{L_1, L_3^k}$ for all $k\geq 1$ would imply the surjectivity of the following map
\[P_L\cap M+\ub{\sigma}\cap M\rightarrow (P_L+\ub{\sigma})\cap M
\] 
for any full-dimensional cone $\ub{\sigma}$ of $\Sigma_{L_3}$. The surjectivity of  this map can be verified by using formal sum of characters, however a moderate bound for the lattices that might fall out of the image seems still desirable allowing for the amount of computation.
\end{rem}

\section{Covering of Smooth Lattice Polygons}
\begin{thm}\label{sfhn} Let $L_1\prec L_2$ be two ample line bundles on a smooth projective toric surface, then the following map is surjective as long as $L_1$ is not $\mathcal{O}_{\mathbb{P}^2}(1)$
\begin{equation}\label{f12}\phi_{L_1, L_2}: P_{L_1}\times P_{L_1^{-1}\otimes L_2}\cap M\rightarrow P_{L_2}.
\end{equation}

\end{thm}

\subsection{Some Preliminaries}\hfill

Let $\bl{u}$ be a given nonzero vector, recall that $l_{\bl{u}}(x)$ is the line through $x$ with direction $\bl{u}$. For a convex polygon or an infinite convex polygonal region $P$ such that $P\cap l_{\bl{u}}(x)$ is a segment for all $x\in P$, we will denote by $\mathcal{L}^P_{\bl{u}}$ or $\mathcal{L}_{\bl{u}}$ the following function
\[P\rightarrow \mathds{R}_{\geq0},\;\;\; x\rightarrow \|P\cap l_{\bl{u}}(x)\|
\]
and by $\mathcal{L}_{\bl{u}}(P)$ the number $\max_{x\in P}\mathcal{L}_{\bl{u}}(x)$.

Let $P$ be a convex polygon, $\pi_{\bl{u}}: P\rightarrow \mathds{R}$ be the projection of $P$ along the direction of $\bl{u}$. Then  $\pi_{\bl{u}}(P)=[a,b]\subset\mathds{R}$ and $\mathcal{L}_{\bl{u}}^P$ can be regarded as a function defined on $[a, b]$ by letting
$\mathcal{L}_{\bl{u}}^P(t)=\mathcal{L}_{\bl{u}}^P(x)$, where $t\in[a, b]$ and $t=\pi_{\bl{u}}(x)$ for some $x\in P$. 

By investigating the Steiner symmetrization  \cite[9.1]{gruber2007convex} of $P$ along the direction of $\bl{u}$, one sees easily $\mathcal{L}^P_{\bl{u}}(t)$ is a convex unimodal function, which means this function attains a unique maximum in the interval where it is defined.

 If $\mathcal{L}_{\bl{u}}(P)\geq\|\bl{u}\|$, then we have a nonempty interval 
$[c, d]\subseteq [a,b]$ such that $\mathcal{L}_{\bl{u}}(x)\geq\|\bl{u}\|$ iff $\pi_{\bl{u}}(x)\in[c,d]$. Moreover,
if $\mathcal{L}_{\bl{u}}(P)=\|\bl{u}\|$ and there exists exactly one chord of $P$ paralleling with $\bl{u}$ and with length $\|\bl{u}\|$, then $c=d$. Otherwise, we have $c<d$.

Following the notations above, we have the following
\begin{definition}\label{ctp}
A contact point of $P$ with respect to $\bl{u}$  (or of $\bl{u}$ for brevity) is a point $x\in \partial P\cap\partial(\bl{u}+P)$ such that $\pi_{\bl{u}}(x)=c$ or $\pi_{\bl{u}}(x)=d$ and $x-\epsilon\bl{u}\notin\bl{u}+P$ for any $\epsilon>0$.
\end{definition}
It is easy to see from our definition $P$ has one or two contact points in the direction of $\bl{u}$ depending on  $c=d$ or $c<d$. The lemma below should be clear from our definition and we omit its proof.

\begin{lem}\label{nc}
Let $x, y$ be the two contact points $P$ in the direction of $\bl{u}$, then $P\cap (\bl{u}+P)$ is contained in the closed region  bounded by the lines $l_{\bl{u}}(x)$ and $l_{\bl{u}}(y)$. 
\end{lem}


For a parallelogram $P$, each pair of its opposite sides defines an infinite region (including the boundary) bounded by the pair of lines containing these sides. A point set is said to be \emph{bounded by one direction} of $P$ if it is contained in such a region. More specifically, a point set is said to be bounded by direction $AB$ (or $CD$) of a parallelogram $ABCD$ if it is contained in the region bounded by the lines $AB$ and $CD$. A point set is said to be \emph{bounded by two directions} of $P$ if it is contained in neither of these infinite regions but their union. The following lemma is easy to verify and we omit its proof.

\begin{lem}\label{ae}Given two parallelograms all of whose vertices stay on the boundary of their convex hull, then each of them is bounded by one or two directions of the other. If furthermore they have no common points, their vertices are contained in a pair of paralleling lines.

\end{lem}


To determine the relative positions of two chords of $P$ will be essential in the proof of Theorem \ref{sfhn}. Next we show this question has a simple answer under suitable conditions and that is enough for our purpose. Before stating the result we make some preparations. 

Given a convex polygon $P$ in the plane and a vertex $x_0$ of it we choose a coordinate system such that $x_0$ is the unique lowest point of  $P$. In the following the slopes of vectors or segments are measured in this chosen coordinate system.

Let $s_1, s_{-1}$ be the two edges of $P$ which share $x_0$ as their common endpoint. The consecutive edges of $P$ starting from $s_1$ (resp. $s_{-1}$) will be denoted by $s_2, s_3,\cdots$ (resp.  $s_{-2}, s_{-3},\cdots$). Let $\mu_i>0$ be the slope of $s_i$, we assume $\mu_1>\mu_{-1}$, $\mu_i>\mu_{i+1}$ when $i>0$ and $\mu_i<\mu_{i-1}$ when $i<0$.

Recall that $P$ can be represented as 
\[P=\bigcap_{s_i} H^+_{s_i},
\]
where $H^+_{s_i}$ is the half plane which is cut out by the line containing the edge $s_i$ and contains $P$. For each pair of positive integers $m, n$ such that $\mu_m>\mu_{-n}$, we define a polygonal region $\mathscr{P}_{m, -n}$ containing $P$ as follows
\[\mathscr{P}_{m, -n}=\bigcap_{\substack{1\leq i\leq m,\\ -n\leq i\leq -1}} H^+_{s_i}.
\]
Obviously we have $P\subset\mathscr{P}_{m,-n}\subseteq \mathscr{P}_{m', -n'}$ for $m\geq m', n\geq n'$. 
By abuse of notation, we use $s_m$ (resp. $s_{-n}$) for the ray, which is part of the boundary of $\mathscr{P}_{m, -n}$, with initial at $s_m\cap s_{m-1}$ (resp. $s_{-n}\cap s_{-n+1}$) and containing the edge $s_m$ (resp. $s_{-n}$) of $P$. 

Next we consider the chords of the region $\mathscr{P}_{m, -n}$ and $P$ with the same directions and lengths. Firstly let $\bl{u}$ be a vector such that $\mathcal{L}_{\bl{u}}^{P}(x_0)=0$, one sees easily the following function
\[\mathcal{L}^{\mathscr{P}_{m, -n}}_{\bl{u}}(x)
\]
is strictly monotonically increasing when $x$ moves upwards along $\bigcup_{i\geq1}s_i$ or $\bigcup_{i\leq-1}s_i$. 

Now
suppose $P$ has a chord with direction $\bl{u}$ and length $a>0$. Then it is easy to see we can find minimal $m$ and $n$ such that a chord of $\mathscr{P}_{m, -n}$ has direction $\bl{u}$, length $a$ and at the same time also a chord of $P$, i.e. the endpoints of this chord are contained in $P$.
Furthermore, for any  $m', n'$ such that $m'\geq m$ and $n'\geq n$ this chord is also contained in $\mathscr{P}_{m', -n'}$. On the other hand, by the minimality of $m$ and $n$, for $\mathscr{P}_{m',-n'}$ with $m'<m$ or $n'<n$ the chord of $\mathscr{P}_{m', -n'}$ with direction $\bl{u}$ and length $a$ would be strictly on the lower side of the common chord of $P$ and $\mathscr{P}_{m, -n}$ as mentioned above.

Now we can state our 
\begin{lem}\label{ip}

Let $\bl{u}, \bl{v}$ be two nonzero vectors, $x_{\bl{u}}, x_{\bl{v}}$ (resp. $y_{\bl{u}}, y_{\bl{v}}$) be points on $s_{-1}$ (resp. $s_1$), which are contained in $\mathscr{P}_{1,-1}$. Suppose that

1) the directions of the chords $\overline{x_{\bl{u}}y_{\bl{u}}}$ and $\overline{x_{\bl{v}}y_{\bl{v}}}$ are $\bl{u}$ and $\bl{v}$ respectively and $\overline{x_{\bl{u}}y_{\bl{u}}}$ is lying on the upper right side of $\overline{x_{\bl{v}}y_{\bl{v}}}$;

2) $\mathcal{L}^P_{\bl{u}}(x_0)=\mathcal{L}^P_{\bl{v}}(x_0)=0$;

3) $\mathcal{L}_{\bl{u}}(P)\geq \|\overline{x_{\bl{u}}y_{\bl{u}}}\|$, $\mathcal{L}_{\bl{v}}(P)\geq \|\overline{x_{\bl{v}}y_{\bl{v}}}\|$.

\noindent Then we can find $m, n$ such that $\mathscr{P}_{m, -n}$ and $P$ share two chords with directions  $\bl{u}, \bl{v},$ and lengths $\|\overline{x_{\bl{u}}y_{\bl{u}}}\|$, $\|\overline{x_{\bl{v}}y_{\bl{v}}}\|$ respectively. Moreover, the chord with direction $\bl{u}$ still lies above the chord with direction $\bl{v}$. 

\end{lem}

\begin{proof} By our argument above, we can find a pair of positive integers $m, n$ such that $\mathscr{P}_{m,-n}$ shares with $P$ a chord whose direction is $\bl{v}$ and length $\|\overline{x_{\bl{v}}y_{\bl{v}}}\|$. Moreover, we assume $m, n$ are minimal with this proerpty. To prove the lemma, we only need to show the chord of $\mathscr{P}_{m,-n}$ with direction $\bl{v}$ and length $\|\overline{x_{\bl{v}}y_{\bl{v}}}\|$ lies on the lower left of the chord with direction $\bl{u}$ and length $\|\overline{x_{\bl{u}}y_{\bl{u}}}\|$. We will prove this claim by induction, i.e. for any $k\geq1, l\geq1$,  the chord of $\mathscr{P}_{k,-l}$ with direction $\bl{v}$ and length $\|\overline{x_{\bl{v}}y_{\bl{v}}}\|$ lies on the lower left of the chord with direction $\bl{u}$ and length $\|\overline{x_{\bl{u}}y_{\bl{u}}}\|$.  

By condition the claim is true for $k=l=1$. Assuming now it has been proved for $\mathscr{P}_{k, -l}$, next we will prove it for $\mathscr{P}_{k, -l-1}$ and $\mathscr{P}_{k+1, -l}$. 
 Let $x_{\bl{u}}^{k, -l}\in\bigcup_{-l\leq i\leq -1}s_i$, $y_{\bl{u}}^{k, -l}\in\bigcup_{1\leq i\leq k}s_i$ (resp. $x_{\bl{v}}^{k, -l}\in\bigcup_{-l\leq i\leq -1}s_i$, $y_{\bl{v}}^{k, -l}\in\bigcup_{1\leq i\leq k}s_i$) be the endpoints of the two chords of $\mathscr{P}_{k, -l}$ as described in the statement of the claim. We will only consider the special case when $x_{\bl{u}}^{k, -l}=x_{\bl{v}}^{k, -l}$, more general case can be proved similarly. As is shown by the figure on the left below, for $\mathscr{P}_{k, -l-1}$, we have $x_{\bl{u}}^{k, -l-1}=x_{\bl{v}}^{k, -l-1}$ and the two chords aforementioned have the expected relative positions. For $\mathscr{P}_{k+1, -l}$ one reads from the figure on the right easily $\mathcal{L}_{\bl{u}}^{\mathscr{P}_{k+1, -l}}(x_{\bl{v}}^{k+1, -l})<\|\overline{x^{k+1,-l}_{\bl{u}}y^{k+1,-l}_{\bl{u}}}\|=\|\overline{x^{k,-l}_{\bl{u}}y^{k,-l}_{\bl{u}}}\|=\|\overline{x_{\bl{u}}y_{\bl{u}}}\|$. Then $x_{\bl{u}}^{k+1, -l}$ hence the chord of $\mathscr{P}_{k+1, -l}$ with direction $\bl{u}$ and length  $\|\overline{x_{\bl{u}}y_{\bl{u}}}\|$ is on the upper right side of the chord with direction $\bl{v}$ and length $\|\overline{x_{\bl{v}}y_{\bl{v}}}\|$.

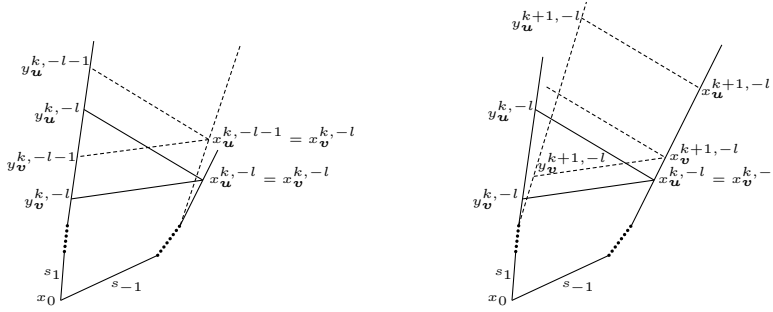
\begin{tikz1}

\draw(0,0)--(0.35, 2.45);
\draw[dash pattern=on 1.5pt off 1pt](1.5,0)--(2.3, 2.4); 

\draw[dash pattern=on 1.5pt off 1pt](1.88,1.14)--(0.13,0.91);
\draw[dash pattern=on 1.5pt off 1pt](1.88,1.14)--(0.3,2.1);

\draw(1.5,0)--(2, 1);

\draw(1.8, 0.6)--(0.05,0.35);
\draw(1.8, 0.6)--(0.22,1.54);

\draw(-0.04,-0.35)--(-0.09, -1);
\draw[line width=1.3pt, line cap=round, dash pattern=on 0pt off 1.5\pgflinewidth](-0.04,-0.35)--(0,0);
\draw(-0.09,-1)--(1.2,-0.4);
\draw[line width=1.3pt, line cap=round, dash pattern=on 0pt off 1.5\pgflinewidth](1.2,-0.4)--(1.5,0);

\node at (2.65,0.64) (0) {\relsize{-20} $x_{\bl{u}}^{k,-l}=x_{\bl{v}}^{k,-l}$};
\node at (-0.3,0.35) (0){\relsize{-20} $y_{\bl{v}}^{k,-l}$};
\node at (2.84,1.18)(0) {\relsize{-20} $x_{\bl{u}}^{k,-l-1}=x_{\bl{v}}^{k,-l}$};
\node at (-0.17,1.5) (0){\relsize{-20} $y_{\bl{u}}^{k,-l}$};
\node at (-0.38,0.85) (0){\relsize{-20} $y_{\bl{v}}^{k,-l-1}$};
\node at (-0.19,2.12)(0) {\relsize{-20} $y_{\bl{u}}^{k,-l-1}$};
\node at (-0.32,-1)(0) {\relsize{-20} $x_0$};
\node at (-0.24,-0.65)(0) {\relsize{-20} $s_1$};
\node at (0.75,-0.83)(0) {\relsize{-20} $s_{-1}$};

-----------------------------------------------------------------------------------
\draw(6,0)--(6.32, 2.24);
\draw[dash pattern=on 1.5pt off 1pt](6,0)--(6.9, 2.97); 
\draw(7.5,0)--(8.7, 2.4);
[dash pattern=on 5pt off 1pt]

\draw(6.22,1.54)--(7.8,0.6);
\draw(7.8,0.6)--(6.05,0.35);

\draw[dash pattern=on 1.5pt off 1pt](6.2,0.65)--(7.95,0.9);
\draw[dash pattern=on 1.5pt off 1pt](6.37,1.84)--(7.95,0.9);

\draw[dash pattern=on 1.5pt off 1pt](6.83,2.76)--(8.41,1.82);

\draw(5.96,-0.35)--(5.91, -1);
\draw[line width=1.3pt, line cap=round, dash pattern=on 0pt off 1.5\pgflinewidth](5.96,-0.35)--(6,0);
\draw(5.91,-1)--(7.2,-0.4);
\draw[line width=1.3pt, line cap=round, dash pattern=on 0pt off 1.5\pgflinewidth](7.2,-0.4)--(7.5,0);

\node at (6.33,2.72) (0) {\relsize{-20} $y_{\bl{u}}^{k+1,-l}$};
\node at (5.86,1.52) (0){\relsize{-20} $y_{\bl{u}}^{k,-l}$};
\node at (8.85,1.82)(0) {\relsize{-20} $x_{\bl{u}}^{k+1,-l}$};
\node at (8.62,0.65) (0){\relsize{-20} $x_{\bl{u}}^{k,-l}=x_{\bl{v}}^{k,-l}$};
\node at (5.66,0.35) (0){\relsize{-20} $y_{\bl{v}}^{k,-l}$};
\node at (6.67,0.83)(0) {\relsize{-20} $y_{\bl{v}}^{k+1,-l}$};
\node at (8.43,0.95)(0){\relsize{-20} $x_{\bl{v}}^{k+1,-l}$};

\node at (5.68,-1)(0) {\relsize{-20} $x_0$};
\node at (5.76,-0.65)(0) {\relsize{-20} $s_1$};
\node at (6.75,-0.83)(0) {\relsize{-20} $s_{-1}$};
\end{tikz1}

\end{proof}

\subsection{The Proof}

\begin{proof}[Proof of Theorem \ref{sfhn}] The proof is given in three steps. 

\subsubsection*{Step 1.  Reduction to Two Cases.}\hfill\medskip

Recall that all smooth projective toric surfaces can be obtained by successively blowing up the fixed points under the torus action. We will prove the proposition by induction on the Picard number $r$ of $X$. Firstly we observe under the condition of Theorem \ref{sfhn},  the surjectivity of (\ref{f12}) can be reformulated as 
\begin{center}
  the lattice translations of $P_{L_1}$ which are contained in $P_{L_2}$ can cover $P_{L_2}$.
\end{center}
\noindent The statement above is easy to verify when $r=2$. Now assume it is proved for all toric surfaces with Picard number $r=k$, next we will prove it for  $r=k+1$. 

Let $X$ be a toric surface obtained from $Y$ with Picard number $k$ by blowing up a $T$-fixed point corresponding to a two-dimensional cone $\sigma\in\Sigma_Y$.  We denote the exceptional curve of this blow-up by $D_{\sigma}$. Then for any ample line bundle $L_1$ on $X$ we have $L_1.D_{\sigma}\geq 1$. Next we will prove the morphism $\phi_{L_1, L_2}$ is surjective for all pairs of ample line bundles $L_1\prec L_2$ by induction on $(s_1, s_2)$, where $s_1=L_1.D_{\sigma}$ and $s_2=L_2.D_{\sigma}$.  When $s_1\geq 2$, the line bundles $L_i(D_{\sigma}), i=1, 2,$ are still ample on $X$ and $L_i(D_{\sigma}).D_{\sigma}=s_i-1$. When $s_1=s_2=1$, $L_i(D_{\sigma})=\pi^*\tilde{L}_i$ for some ample line bundle $\tilde{L}_i,  i=1, 2,$ on $Y$, where $\pi:X\rightarrow Y$ is the projection.

Notice that $s_2\geq s_1\geq1$ and it suffices to prove the following implications

(I1) $(s_1-1, s_2-1)\Rightarrow (s_1, s_2)$ for $s_1\geq 2$ and $s_1=s_2=1$, 

(I2) $(1, s_2)\Rightarrow (1, s_2+1)$ for $s_2\geq 1$.

The proofs of the of (I1) and (I2) will be given in the next two steps.

Before going further, we fix some notations. Let $\alpha,\beta\in\Sigma_Y(1)$ such that $\sigma(1)=\{\alpha, \beta\}$ and $\alpha_1, \alpha_2, \cdots, \alpha_k$ (resp. $\beta_1, \beta_2, \cdots, \beta_l$) be the primitive vectors corresponding to one-dimensional cones of $\Sigma_Y$ lying between $\alpha$ and $-\beta$ (resp. $-\alpha$ and $\beta$). For convenience we occasionally use $\alpha_0$ and $\beta_0$ to denote $\alpha$ and $\beta$ respectively. Moreover, we do not exclude the possibilities that $\alpha_k=-\beta$ and $\beta_l=-\alpha$. 
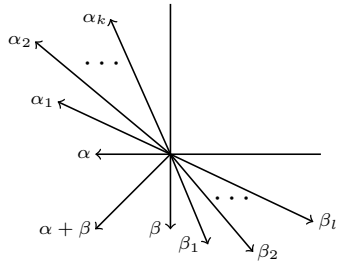
\[
\begin{tikzpicture}
\draw[->, line width=0.6pt] (2,0)--(-1,0);
\draw[->, line width=0.6pt] (0,2)--(0,-1);
\draw[->, line width=0.6pt] (0,0)--(-1,-1);

\draw[->, line width=0.6pt] (0,0)--(-1.5,0.7);
\draw[->, line width=0.6pt] (0,0)--(-1.8,1.5);
\draw[->, line width=0.6pt] (0,0)--(-0.8,1.8);

\draw[->, line width=0.6pt] (0,0)--(0.5,-1.2);
\draw[->, line width=0.6pt] (0,0)--(1.1, -1.3);
\draw[->, line width=0.6pt] (0,0)--(1.9, -0.9);

\node at (-1.15,0)  (0) {\relsize{-2}$\alpha$}; 
\node at (-0.2,-1)  (0) {\relsize{-2}$\beta$}; 
\node at (-1.4,-1)  (0) {\relsize{-2}$\alpha+\beta$};

\node at (-1.7,0.7) (0) {\relsize{-2}$\alpha_1$}; 
\node at(-2,1.5) (0) {\relsize{-2}$\alpha_2$}; 
\node at(-1,1.8) (0) {\relsize{-2}$\alpha_k$}; 
\node at(-0.88,1.2) (0) {\relsize{1.2}$\cdots$};

\node at (0.25,-1.2) (0) {\relsize{-2}$\beta_1$}; 
\node at(1.3,-1.3) (0) {\relsize{-2}$\beta_2$}; 
\node at(2.1,-0.9) (0) {\relsize{-2}$\beta_l$}; 
\node at(0.85,-0.6) (0) {\relsize{1.2}$\cdots$}; 

\end{tikzpicture}
\]

Moreover, we will choose a basis $\bl{u}_h, \bl{u}_v$ of $M_{\mathds{R}}$ defined by $\langle\bl{u}_h,\alpha\rangle=\langle\bl{u}_v,\beta\rangle=-1, \langle\bl{u}_h,\beta\rangle=\langle\bl{u}_v,\alpha\rangle=0$. For each $0\leq i\leq k$ and $ 0\leq j\leq l$, we take $\bl{u}_{\alpha_i}$ and $\bl{u}_{\beta_j}$ to be the primitive integral vectors with non-negative coefficients under the basis $\{\bl{u}_h, \bl{u}_v\}$ such that 
$\langle\bl{u}_{\alpha_i},\alpha_i\rangle=0$ and $\langle\bl{u}_{\beta_j},\beta_j\rangle=0$ respectively, particularly $\bl{u}_{\alpha_0}=\bl{u}_v, \bl{u}_{\beta_0}=\bl{u}_h$. 

The following lemmas will be used later and proof of the first is omitted.

\begin{lem}\label{1q1c}Let $-\sigma$ be the two-dimensional cone spanned by $-\alpha, -\beta$, then $-\sigma\cap\Sigma_X(1)\neq\emptyset$ iff $-\alpha$ or $-\beta$ is contained in $\Sigma_X(1)$. Moreover,  $(-\sigma\backslash\{-\alpha, -\beta\})\cap\Sigma_X(1)\neq\emptyset$ implies  $-\alpha, -\beta\in\Sigma_X(1)$.

\end{lem}

\begin{lem}\label{2ctp}Following notations above, $P_{L_1}$ has two contact points in the directions of $\bl{u}_{\alpha_i}$ and $\bl{u}_{\beta_j}$ for each $0\leq i\leq k$ and $0\leq j\leq l$. 

\end{lem}

\begin{proof}
For $0\leq i\leq k-1$ and $0\leq j\leq l-1$ the conclusion is obviously true as one sees easily
\[\mathcal{L}_{\bl{u}_{\alpha_i}}(P_{L_1})>\|\bl{u}_{\alpha_i}\|,\;\;\mathcal{L}_{\bl{u}_{\beta_j}}(P_{L_1})>\|\bl{u}_{\beta_j}\|.
\]
If $\mathcal{L}_{\bl{u}_{\alpha_k}}(P_{L_1})=\|\bl{u}_{\alpha_k}\|$ or $\mathcal{L}_{\bl{u}_{\beta_l}}(P_{L_1})=\|\bl{u}_{\beta_l}\|$, then we have 
$-\sigma(1)\subseteq\{-\alpha, -\beta\}$. In particular, $e_{\alpha_k}(P_{L_1})$ and $e_{\beta_l}(P_{L_1})$ share a common endpoint, say $x$. Then by Lemma \ref{upp}, we have $\bl{u}_{\alpha_k}+\bl{u}_{\beta_l}+x\in P_{L_1}$, hence the conclusion is also true in this case. 

\end{proof}

\subsubsection*{Step 2. Proof of (I1).}\hfill\medskip

Assuming from now on the morphism $\phi_{L_1(D_{\sigma}), L_2(D_{\sigma})}$ is surjective, we will show so is true for $\phi_{L_1, L_2}$.

\smallskip
\emph{Step 2.1. The Set of Translation Vectors.}\smallskip

By our choice of the integral basis of $M$, the polygons corresponding to $L_i, i=1,2,$ can be obtained from those corresponding to $L_i(D_{\sigma})$ by cutting off the shading region, as is shown in the figure below. There are two possibilities for the relative positions of the polygons corresponding to  $L_1$ and $L_2$.
 
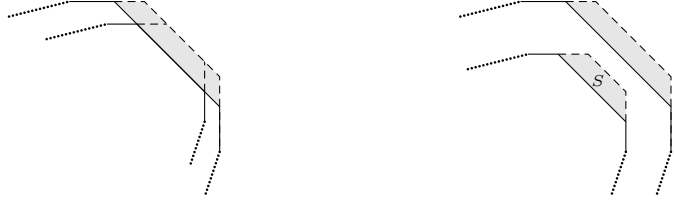
\begin{figure}[H]

\begin{tikzpicture}

\draw (-6,0)--(-5.4,0)--(-4,-1.4)--(-4,-2);
\draw [densely dashed] (-5.4,0)--(-5,0)--(-4,-1)--(-4,-2);
\draw (-5.5, -0.3)--(-5.1, -0.3)--(-4.2,-1.2)--(-4.2, -1.6);

\draw[densely dashed] (-5.1, -0.3)--(-4.7, -0.3);
\draw[densely dashed] (-4.2, -1.2)--(-4.2, -0.8);

\fill [fill opacity=0.1] (-5.4, 0)--(-5,0)--(-4,-1)--(-4,-1.4)--(-5.4, 0);
\draw[line width=1pt, line cap=round, dash pattern=on 0pt off 1.5\pgflinewidth](-6.8, -0.2)--(-6,0);
\draw[line width=1pt, line cap=round, dash pattern=on 0pt off 1.5\pgflinewidth](-6.3, -0.5)--(-5.5,-0.3);
\draw[line width=1pt, line cap=round, dash pattern=on 0pt off 1.5\pgflinewidth](-4, -2)--(-4.2,-2.6);
\draw[line width=1pt, line cap=round, dash pattern=on 0pt off 1.5\pgflinewidth](-4.2, -1.6)--(-4.4,-2.2);
--------------------------------------------------------------------------------------------------------------------------------------------------------------------------------

\draw (0,0)--(0.6,0)--(2,-1.4)--(2,-2);
\draw [densely dashed] (0.6,0)--(1,0)--(2,-1)--(2,-2);
\fill [fill opacity=0.1] (0.6, 0)--(1,0)--(2,-1)--(2,-1.4)--(0.6, 0);
\draw[line width=1pt, line cap=round, dash pattern=on 0pt off 1.5\pgflinewidth](-0.8, -0.2)--(0,0);
\draw[line width=1pt, line cap=round, dash pattern=on 0pt off 1.5\pgflinewidth](-0.7, -0.9)--(0.1,-0.7);

\draw (0.1, -0.7)--(0.5, -0.7)--(1.4,-1.6)--(1.4, -2);
\draw[densely dashed] (0.5, -0.7)--(0.9, -0.7)--(1.4,-1.2)--(1.4,-1.6);
\fill[fill opacity=0.1] (0.5, -0.7)--(0.9,-0.7)--(1.4,-1.2)--(1.4,-1.6)--(0.5, -0.7);
\draw[line width=1pt, line cap=round, dash pattern=on 0pt off 1.5\pgflinewidth](2, -2)--(1.8,-2.6);
\draw[line width=1pt, line cap=round, dash pattern=on 0pt off 1.5\pgflinewidth](1.4, -2)--(1.2,-2.6);
\node at (1.03,-1.05) (0){\relsize{-3}$S$};
        \end{tikzpicture}
\end{figure}    
     
\noindent For the case on the left, there is nothing to prove. For the case on the right, we need to show the shading region $S$ of $P_{L_1(D_{\sigma})}$ is contained in the union of some lattice translations of $P_{L_1}$, each of which is contained in $P_{L_2}$. 

Let $A_0, A'_0$ be the endpoints of $e_{\alpha+\beta}(P_{L_1})$, if $\mathcal{L}_{\bl{u}_{\alpha_i}}(A_0)\geq\|\bl{u}_{\alpha_i}\|$,  $\mathcal{L}_{\bl{u}_{\alpha_i}}(A'_0)\geq\|\bl{u}_{\alpha_i}\|$ for some $0\leq i\leq k$, then we can deduce $S\subset \bl{u}_{\alpha_i}+P_{L_1}$. We have similar result when $\alpha_i$ is replaced by $\beta_j$ for some $0\leq j\leq l$. In the following we will work under the 
assumption below
\begin{equation}\label{dgn} S\;\text{cannot be coverved by either a single}\;\bl{u}_{\alpha_i}+P_{L_1}\;\text{or}\;\bl{u}_{\beta_j}+P_{L_1}.
\end{equation}

Let $p_0\in P_{L_1^{-1}\otimes L_2}\cap M$ be the lattice point corresponding to $P_{L_1}$. In order to find the translations of $P_{L_1}$ which can cover $S$, we will first prove the following
\begin{claim}\label{cl2.0}There exists lattice points $q, q'\in P_{L_1^{-1}\otimes L_2}\cap M$ such that for any point $x$ on the segment $\overline{qq'}$  
\begin{equation}\label{2cps}\mathcal{L}_{\overrightarrow{p_0x}}(P_{L_1})>\|\overrightarrow{p_0x}\|.
\end{equation}
\end{claim}
Then we will show
\begin{equation}\label{1i1}S\subseteq \bigcup _{p\in p_0qq'\cap M}\overrightarrow{p_0p}+P_{L_1},
\end{equation}
where $p_0qq'$ is the triangle with vertices $p_0, q$ and $q'$.  





\smallskip

We first prove this claim under the condition $p_0\notin\partial P_{L_1^{-1}\otimes L_2}$. The proof for the case when $p_0\in\partial P_{L_1^{-1}\otimes L_2}$ will be given at the end of Step 2.1. 
\smallskip

\begin{proof}[Proof of Claim \ref{cl2.0} for $p_0\notin\partial P_{L_1^{-1}\otimes L_2}$.]   For each $0\leq i\leq k$ and $0\leq j\leq l$, let 
\begin{equation}\label{lmbd}\lambda_{\alpha_i}=\max_{\lambda\in\mathds{R}_{\geq0}}\{\lambda\;|\;\lambda\bl{u}_{\alpha_i}+P_{L_1}\subset P_{L_2}\},\;\; \;\lambda_{\beta_j}=\max_{\lambda\in\mathds{R}_{\geq0}}\{\lambda\;|\; \lambda\bl{u}_{\beta_j}+P_{L_1}\subset P_{L_2}\}.\end{equation}

Since $p_0\notin\partial P_{L_0^{-1}\otimes L_2}$ we have \begin{equation}\label{cdz}\lambda_{\alpha_i}>0,\: \lambda_{\beta_j}>0\;\;\text{for all}\;\; 0\leq i\leq k, \:0\leq j\leq l.\end{equation}
Then there exists a $\rho\in\Sigma_X(1)$ such that
\[e_{\rho}(\lambda_{\alpha_i}\bl{u}_{\alpha_i}+P_{L_1})\subseteq e_{\rho}(P_{L_2}).\]
If $\lambda_{\alpha_i}\notin\mathds{Z}$ and particularly when $\lambda_{\alpha_i}<1$, the translation of any lattice point of $P_{L_1}$ under $\lambda_{\alpha_i}\bl{u}_{\alpha_i}$ is not a lattice point. Then $\rho\in\Sigma_X(1)$ satisfying the relation above is unique, we will denote it by $\iota(\alpha_i)$. In the same way we can define $\iota(\beta_j)$.

Since the Picard number of the base toric surface is  at least three, there exist $0\leq k_0\leq k,  0\leq l_0\leq l$ such that $\alpha_{k_0}+\beta_{l_0}=0$ hence $\bl{u}_{\alpha_{k_0}}=\bl{u}_{\beta_{l_0}}$. If $\lambda_{\alpha_{k_0}}=\lambda_{\beta_{l_0}}\geq1$, then \[\mathfrak{L}_{\bl{u}_{\alpha_{k_0}}}(A_0)\geq\|\bl{u}_{\alpha_{k_0}}\|,\;\;\mathfrak{L}_{\bl{u}_{\alpha_{k_0}}}(A'_0)\geq\|\bl{u}_{\alpha_{k_0}}\|,\] hence $S\subset \bl{u}_{\alpha_{k_0}}+P_{L_1}$, which is excluded by our assumption at the beginning (\ref{dgn}). Then the  condition (\ref{cdz}) specializes to \begin{equation}\label{01abt}0< \lambda_{\alpha_{k_0}}=\lambda_{\beta_{l_0}}<1.\end{equation}


Next we will show there exists some $1\leq i\leq k$ and $1\leq j\leq l$ such that
\begin{equation}\label{k1} \lambda_{\alpha_i}\geq1;
\end{equation}
\begin{equation} \label{l1}\lambda_{\beta_j}\geq1.
\end{equation}

Next we will show (\ref{k1}) and (\ref{l1}) under the condition
$\iota(\alpha_{k_0})=\iota(\beta_{l_0})=\beta_s$ for some $0\leq s\leq l$. The cases when
$\iota(\alpha_{k_0})=\iota(\beta_{l_0})=\alpha_s$ for some $0\leq s\leq k$ or $\alpha+\beta$ can be proved similarly and we omit them. 

First we prove (\ref{l1}). By assumption it is easy to see $s<l_0$. Besides, one sees easily $\iota(\beta_0)=\iota(\beta)\neq \beta_j$ for any $j$. In addition, by (\ref{dgn}) $\iota(\beta_0)\neq \alpha+\beta$ since otherwise we will have $\lambda_{\beta_0}\geq 1$ and $S\subset \bl{u}_{\beta_0}+P_{L_1}$. Then $\iota(\beta_0)=\alpha_i$ for some $0\leq i\leq k$. 

If the contratry of (\ref{l1}) is true, i.e. $\lambda_{\beta_j}<1$ for all $0\leq j\leq l$, we will take
\[l_1=\min_{1\leq j\leq l_0}\{j\:|\: \iota(\beta_j)=\beta_{s_1}\;\text{for some}\; 0\leq s_1\leq l\}.
\]
By abuse of notaition, let $\iota(\beta_{l_1})=\beta_{s_1}$, then $s_1\leq l_1-1$ and $\iota(\beta_{l_1-1})\in\{\alpha_0, \alpha_1, \cdots,\alpha_k, \alpha+\beta\}$. Next we assume $\iota(\beta_{l_1-1})=\alpha_t$ for some $0\leq t\leq k$ and the case $\iota(\beta_{l_1-1})=\alpha+\beta$ is omitted as it can be treated similarly. As is shown in the figure on the left below, for each $\alpha_i$ with $0\leq i\leq t$,  $\beta_j$ with $0\leq j\leq s_1$ and $\alpha+\beta$ the corresponding edge of $P_{L_1^{-1}\otimes L_2}$ (if its length is nonzero)  has at least one of its endpoints contained in the angular region with vertex at $p_0$ and bounded by the rays $r_{\bl{u}_{\beta_{l_1}}}(p_0)$ and $r_{\bl{u}_{\beta_{l_1-1}}}(p_0)$. In particular, this angular region contains at least one lattice point of $P_{L_1^{-1}\otimes L_2}$.  However, this is impossible if one notices $\bl{u}_{\beta_{l_1}}, \bl{u}_{\beta_{l_1-1}}$ constitute a basis of $M$, $\lambda_{\beta_{l_1}}<1, \lambda_{\beta_{l_1-1}}<1$ and moreover the slope of $e_{\beta_{s_1}}=e_{\beta_{s_1}}(P_{L_1^{-1}\otimes L_2})$ is no larger than that of $\bl{u}_{\beta_{l_1-1}}$ (recall that $s_1\leq l_1-1$). As a consequence, we have $\lambda_{\beta_{l_1}}\geq1$ for some $1\leq l_1\leq l_0$.

 
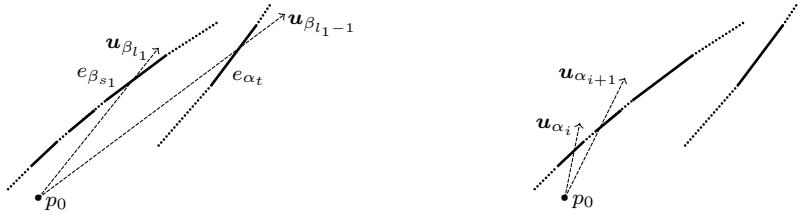
\begin{figure}[H]
\begin{tikzpicture}
\draw[->,densely dashed,dash pattern=on 1.5pt off 0.5pt] (-7,0)--(-5.4,2);

\draw[->,densely dashed,dash pattern=on 1.5pt off 0.5pt] (-7,0)--(-3.72,2.45);

\draw [line width=1pt](-6.1,1.3)--(-5.3,1.9);
\draw [line width=1pt, line cap=round, dash pattern=on 0pt off 1.5\pgflinewidth](-5.3,1.9)--(-4.6,2.35);
\draw [line width=1pt, line cap=round, dash pattern=on 0pt off 1.5\pgflinewidth](-6.4,1.05)--(-6.1,1.3);
\draw [line width=1pt] (-6.25,1.17)--(-6.6,0.88);

\draw [line width=1pt, line cap=round, dash pattern=on 0pt off 1.5\pgflinewidth] (-6.75,0.75)--(-6.6,0.88);
\draw [line width=1pt](-7.1,0.42)--(-6.75,0.75);

\draw [line width=1pt, line cap=round, dash pattern=on 0pt off 1.5\pgflinewidth](-7.4,0.12)--(-7.1,0.42);
------------------------------------------------------------------------------------
\draw [line width=1pt](-4.7,1.5)--(-4.1,2.3);

\draw[line width=1pt, line cap=round, dash pattern=on 0pt off 1.5\pgflinewidth]((-4.7,1.5)--(-5.4,0.7);
\draw[line width=1pt, line cap=round, dash pattern=on 0pt off 1.5\pgflinewidth]((-4.1,2.3)--(-3.9,2.6);

\node at(-6.755,-0.07) (0) {\relsize{-1}$p_0$}; \node at(-7,0) (0) {\relsize{-7}$\bullet$};  
\node at(-5.8,1.98) (0) {\relsize{-1}$\bl{u}_{\beta_{l_1}}$}; 
\node at(-3.22,2.3) (0) {\relsize{-1}$\bl{u}_{\beta_{l_1-1}}$}; 
\node at(-6.2,1.6) (0) {\relsize{-1}$e_{\beta_{s_1}}$}; 
\node at(-4.2,1.6) (0) {\relsize{-1}$e_{\alpha_t}$}; 
===========================================

\draw[->,densely dashed,dash pattern=on 1.5pt off 0.5pt] (0,0)--(0.2,1);

\draw[->,densely dashed,dash pattern=on 1.5pt off 0.5pt] (0,0)--(0.8,1.6);

\draw [line width=1pt](0.9,1.3)--(1.7,1.9);
\draw [line width=1pt, line cap=round, dash pattern=on 0pt off 1.5\pgflinewidth](1.7,1.9)--(2.4,2.35);
\draw [line width=1pt, line cap=round, dash pattern=on 0pt off 1.5\pgflinewidth](0.6,1.05)--(0.9,1.3);
\draw [line width=1pt] (0.75,1.17)--(0.4,0.88);

\draw [line width=1pt, line cap=round, dash pattern=on 0pt off 1.5\pgflinewidth] (0.25,0.75)--(0.4,0.88);
\draw [line width=1pt](-0.1,0.42)--(0.25,0.75);

\draw [line width=1pt, line cap=round, dash pattern=on 0pt off 1.5\pgflinewidth](-0.4,0.12)--(-0.1,0.42);
------------------------------------------------------------------------------------
\draw [line width=1pt](2.3,1.5)--(2.9,2.3);

\draw[line width=1pt, line cap=round, dash pattern=on 0pt off 1.5\pgflinewidth]((2.3,1.5)--(1.6,0.7);
\draw[line width=1pt, line cap=round, dash pattern=on 0pt off 1.5\pgflinewidth]((2.9,2.3)--(3.1,2.6);

\node at(0.25,-0.05) (0) {\relsize{-1}$p_0$}; \node at(0,0) (0) {\relsize{-7}$\bullet$};  
\node at(-0.13,0.95) (0) {\relsize{-1}$\bl{u}_{\alpha_i}$}; 
\node at(0.32,1.61) (0) {\relsize{-1}$\bl{u}_{\alpha_{i+1}}$}; 
\end{tikzpicture}
\end{figure}
Next we prove (\ref{k1}) by contradiction, thus we assume $0< \lambda_{\alpha_i}<1$ for each $0\leq i\leq k$. If $\iota(\alpha_{t_2})=\alpha_{t_1}$ for some $0\leq t_1< t_2\leq k$, then by an argument similar as above, we can find $1\leq k_1\leq k$ satisfying (\ref{k1}).  Otherwise $\iota(\alpha_i)\in\{\alpha+\beta, \beta_0, \beta_1,\cdots, \beta_l\}$ for all $0\leq i\leq k$. It can be seen easily from the figure above on the right there exist no lattice points of $P_{L_1^{-1}\otimes L_2}$ lying in the intersection of $P_{L_1^{-1}\otimes L_2}$ with the angular region bounded by $r_{\bl{u}_{\alpha_i}}(p_0)$ and $r_{\bl{u}_{\alpha_{i+1}}}(p_0)$, hence the intersection of $P_{L_1^{-1}\otimes L_2}$ with the angular region bounded by $r_{\bl{u}_{\alpha_0}}(p_0)$ and $r_{\bl{u}_{\alpha_k}}(p_0)$ either. 

On the other hand, we have already found $1\leq l_1\leq l$ such that $\beta_{l_1}$ satisfies (\ref{l1}). As a consequence, the slope of $\bl{u}_{\beta_{l_1}}$ is smaller than that of $\bl{u}_{\alpha_k}$. Then we can deduce the slope of $\bl{u}_k$ is positive hence $P_{L_1}$ has a unique lowest point. As the slope of $\bl{u}_{\beta_l}$ is strictly larger than that of $\bl{u}_{\alpha_k}$ we must have $l_1<l$. Therefore, we can deduce 
\[\mathcal{L}_{\bl{u}_{\beta_{l_1}}}(P_{L_1})>\|\bl{u}_{\beta_{l_1}}\|.\] 

Then $P_{L_1}$ has two contact points in the direction of $\bl{u}_{\beta_{l_1}}$, particularly we can find a unique pair of points $y_1, y_2\in\partial P_{L_1}\backslash e_{\beta_{l_1}}(P_{L_1})$ such that $\overrightarrow{y_1y_2}=\bl{u}_{\beta_{l_1}}$
and $y_2$ is a contact point of $P_{L_1}$ with respect to $\bl{u}_{\beta_{l_1}}$.

Next we locate the points $y_1$ and $y_2$ on $\partial P_{L_1}$. Firslty by (\ref{dgn}) we have $S\not\subset\bl{u}_{\beta_{l_1}}+P_{L_1}$, hence $y_2\notin\bigcup_{0\leq i\leq k}e_{\alpha_i}(P_{L_1})$. By our argument above, the slope of $\bl{u}_{\beta_{l_1}}$ is small than that of  $\bl{u}_{\alpha_i}, 0\leq i\leq k$, then $y_1\notin\bigcup_{0\leq i\leq k}e_{\alpha_i}(P_{L_1})$. For slope reasons $y_2$ is contained in  $\bigcup_{0\leq j\leq l_1-1}e_{\beta_j}(P_{L_1})$ or $e_{\alpha+\beta}(P_{L_1})$ and $y_1$ must be contained in $\bigcup_{l_1+1\leq j\leq l}e_{\beta_j}(P_{L_1})$, which would imply 
$\|\overrightarrow{y_1y_2}\|>\|e_{\beta_{l_1}}(P_{L_1})\|\geq\|\bl{u}_{\beta_{l_1}}\|$, a contradiction. As a consequence, (\ref{k1}) is proved.





 Next we prove (\ref{2cps}). For saving notations, we just assume $\lambda_{\alpha_i}>1, \lambda_{\beta_j}>1$ from now on. By Lemma \ref{2ctp} $\bl{u}_{\alpha_i}$ and $\bl{u}_{\beta_j}$ satisfy one of the following conditions. 
 
1) $\mathcal{L}_{\bl{u}_{\alpha_i}}(P_{L_1})>\|\bl{u}_{\alpha_i}\|$ or $\mathcal{L}_{\bl{u}_{\beta_j}}(P_{L_1})>\|\bl{u}_{\beta_j}\|$;

2) $\mathcal{L}_{\bl{u}_{\alpha_i}}(P_{L_1})=\|\bl{u}_{\alpha_i}\|$, $\mathcal{L}_{\bl{u}_{\beta_j}}(P_{L_1})=\|\bl{u}_{\beta_j}\|$ and $P_{L_1}$ has two contact points with respect to any of $\bl{u}_{\alpha_i}$ and $\bl{u}_{\beta_j}$. 

\noindent For the latter, we have $i=k, j=l$, $\bl{u}_{\alpha_{k-1}}$ parallels with $\bl{u}_{\beta_l}$ and $\bl{u}_{\alpha_k}$ parallels with $\bl{u}_{\beta_{l-1}}$. Then by (\ref{k1}) and (\ref{l1}) we have $S\subset\bl{u}_{\alpha_k}+P_{L_1}$ and $S\subset\bl{u}_{\beta_l}+P_{L_1}$. Next we assume $\mathcal{L}_{\bl{u}_{\beta_j}}(P_{L_1})>\|\bl{u}_{\beta_j}\|$ and the case when $\mathcal{L}_{\bl{u}_{\alpha_i}}(P_{L_1})>\|\bl{u}_{\alpha_i}\|$ can be treated similarly. Take $A, B\in e_{\alpha_i}(P_{L_1})$ such that $\|\overline{AB}\|=\|\bl{u}_{\alpha_i}\|$. If $\|e_{\beta_j}(P_{L_1})\|>\|\bl{u}_{\beta_j}\|$, then we can find a segment on $e_{\beta_{l_1}}(P_{L_1})$ longer than $\bl{u}_{\beta_j}$. Otherwise, by assumption we can find a chord of $P_{L_1}$ near $e_{\beta_j}(P_{L_1})$ such that it parallels with but is longer than $\bl{u}_{\beta_j}$ and has empty intersection with $\overline{AB}$. In any of the two possibilities we will denote this longer segment by $\overline{CD}$.

The figures below exhibit these segments and their convex hull under the conditions that $\|\overline{AB}\|<\|\overline{CD}\|$ together with the slope of  $\overline{CD}$ is lareger and smaller than that of $\overline{AB}$ respectively. The case when  $\|\overline{AB}\|\geq\|\overline{CD}\|$ is omitted as the proof is similar. In these two cases we take $D'$ and $C'$ on $\overline{CD}$ such that $\|\overline{CD'}\|=\|\overline{C'D}\|=\|\bl{u}_{\beta_j}\|$. From the figures below one sees easily for any point $x_1 \in\overline{D'E}$  and $x_2\in\overline{C'E}$ we have
\begin{equation}\label{fro} \mathcal{L}_{\overrightarrow{Cx_1}}(C)>\|\overrightarrow{Cx_1}\|, \;\;\mathcal{L}_{\overrightarrow{Dx_2}}(D)>\|\overrightarrow{Dx_2}\|.
\end{equation} 
Then  (\ref{2cps}) can be proved by taking $q=p_0+\bl{u}_{\alpha_i}$ and $q'=p_0+\bl{u}_{\beta_j}$. 
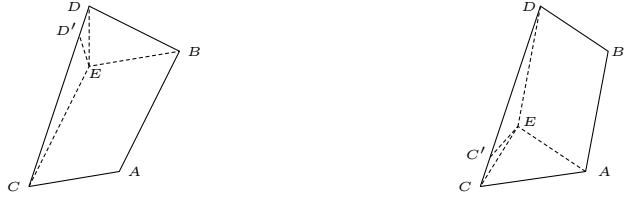
\begin{figure}[H]

\begin{tikzpicture}
\draw (-6,0)--(-5.2,2.4)--(-4,1.8)--(-4.8,0.2)--(-6,0);
\node at(-6.2,0) (0) {\relsize{-7}$C$};
\node at(-5.4,2.4) (0) {\relsize{-7}$D$};
\node at(-5.5,2.1) (0) {\relsize{-7}$D'$};
\node at(-4.6,0.2) (0) {\relsize{-7}$A$};
\node at(-3.8,1.8) (0) {\relsize{-7}$B$};
\node at(-5.12,1.5) (0) {\relsize{-7}$E$};
\draw [dash pattern=on 1.5pt off 1pt]  (-6,0)--(-5.2,1.6)--(-4,1.8);
\draw [dash pattern=on 1.5pt off 1pt]  (-5.2,1.6)--(-5.2,2.4);
\draw [dash pattern=on 1.5pt off 1pt]  (-5.2,1.6)--(-5.33,2.01);
----------------------------------------------------------------------------
\draw (0,0)--(0.8,2.4)--(1.7,1.8)--(1.4,0.2)--(0,0);
\node at(-0.2,0) (0) {\relsize{-7}$C$};
\node at(-0.05,0.45) (0) {\relsize{-7}$C'$};
\node at(0.65,2.4) (0) {\relsize{-7}$D$};
\node at(1.65,0.2) (0) {\relsize{-7}$A$};
\node at(1.83,1.8) (0) {\relsize{-7}$B$};
\node at(0.66,0.86) (0) {\relsize{-7}$E$};
\draw [dash pattern=on 1.5pt off 1pt]  (1.4,0.2)--(0.5,0.8)--(0,0);
\draw [dash pattern=on 1.5pt off 1pt]  (0.8,2.4)--(0.5,0.8);
\draw [dash pattern=on 1.5pt off 1pt]  (0.13,0.39)--(0.5,0.8);
\end{tikzpicture}

\end{figure}
\end{proof}
 Next we introduce some notations, which will be used in the proofs of 
Claim \ref{cl2.0} for $p_0\in\partial P_{L_1^{-1}\otimes L_2}$ and (\ref{1i1}). First of all, we label the primitive vectors with initial at $p_0$ and end at a lattice point in the triangle $p_0qq'$ (including the boundary) clockwisely by $\overrightarrow{p_0p_i}, 1\leq i\leq m$. Then it is easy to see each triangle $p_0p_ip_{i+1}, 1\leq i\leq m-1$, is unimodular. 

Let $A_0, A_0'$
be the endpoints of $e_{\alpha+\beta}(P_{L_1})$ with $A_0$ on the upper left side of $A'_0$ and
\[\mathcal{L}_{\overrightarrow{p_0p_i}}(\overline{A_0A'_0})=\max_{x\in\overline{A_0A'_0}}\mathcal{L}_{\overrightarrow{p_0p_i}}(x).
\]
By combining the unimodality of the function $\mathcal{L}_{\overrightarrow{p_0p_i}}(x), x\in P_{L_1}$ and the relation between $\mathcal{L}_{\overrightarrow{p_0p_i}}(\overline{A_0A'_0})$ and $\|\overrightarrow{p_0p_i}\|$, one checks easily each translation vector belongs to one of the following types. 

a)  $\mathcal{L}_{\overrightarrow{p_0p_i}}(\overline{A_0A'_0})<\|\overrightarrow{p_0p_i}\|$, $\mathcal{L}_{\overrightarrow{p_0p_i}}(A_0)>\mathcal{L}_{\overrightarrow{p_0p_i}}(A'_0)$;

b)  $\mathcal{L}_{\overrightarrow{p_0p_i}}(A_0)\geq\|\overrightarrow{p_0p_i}\|>\mathcal{L}_{\overrightarrow{p_0p_i}}(A'_0)$; 

c) $\mathcal{L}_{\overrightarrow{p_0p_i}}(\overline{A_0A'_0})\geq\|\overrightarrow{p_0p_i}\|$, $\mathcal{L}_{\overrightarrow{p_0p_i}}(A_0)\geq\|\overrightarrow{p_0p_i}\|$, $\mathcal{L}_{\overrightarrow{p_0p_i}}(A'_0)\geq\|\overrightarrow{p_0p_i}\|$;

d) $\mathcal{L}_{\overrightarrow{p_0p_i}}(A'_0)\geq\|\overrightarrow{p_0p_i}\|>\mathcal{L}_{\overrightarrow{p_0p_i}}(A_0)$; 

e) $\mathcal{L}_{\overrightarrow{p_0p_i}}(\overline{A_0A'_0})<\|\overrightarrow{p_0p_i}\|$, $\mathcal{L}_{\overrightarrow{p_0p_i}}(A'_0)>\mathcal{L}_{\overrightarrow{p_0p_i}}(A_0)$;

f) $\mathcal{L}_{\overrightarrow{p_0p_i}}(\overline{A_0A'_0})\geq\|\overrightarrow{p_0p_i}\|$, $\mathcal{L}_{\overrightarrow{p_0p_i}}(A_0)<\|\overrightarrow{p_0p_i}\|$, $\mathcal{L}_{\overrightarrow{p_0p_i}}(A'_0)< \|\overrightarrow{p_0p_i}\|$.

\smallskip

Let $\overrightarrow{p_0q}=a\bl{u}_h+b\bl{u}_v, a\geq 0, b\geq 0$ be a nonzero vector, for convenience in our later argument we will say $\overrightarrow{p_0q}$ is of type g) (resp. type h)) if $P_{L_1}$ has one (resp. zero) contact point with respect to $\overrightarrow{p_0q}$. 

Next we will label the contact points for each vector. Let \[\gamma_i: P_{L_1}\rightarrow A_0A'_0,\;\;x\mapsto l_{\overrightarrow{p_0p_i}}(x)\cap A_0A'_0,\] 
then one can define a real-valued function on $\gamma_i(P_{L_1})$ simply by assigning $\mathcal{L}_{\overrightarrow{p_0p_i}}(x)$ to $t$ for $t=\gamma_i(x)$ and this function is still unimodal. Let $A_i, A'_i$ be the two contact points of $P_{L_1}$ in the direction of $\overrightarrow{p_0p_i}$ we assume $\tilde{A}_i=\gamma_i(A_i)$ always lies on the upper left of $ \tilde{A}'_i=\gamma_i(A'_i)$. Then the relative positions of $\tilde{A}_i, \tilde{A}'_i$ and $A_0, A'_0$ on the line $A_0A'_0$ for various types are given as follows




a) $\tilde{A}_i, \tilde{A}'_i, A_0, A'_0$, ($\tilde{A}'_i\neq A_0$);

b) $\tilde{A}_i, A_0, A'_i, A'_0$, ($A'_i\neq A_0, A'_i\neq A'_0,\tilde{A}'_i=A'_i)$;

c) $ \tilde{A}_i, A_0, A'_0,  \tilde{A}'_i$;

d) $A_0, A_i, A'_0, \tilde{A}'_i$, ($A_i\neq A_0, A_i\neq A'_0, \tilde{A}_i=A_i)$;

e) $A_0, A'_0, \tilde{A}_i, \tilde{A}'_i$, ($\tilde{A}_i\neq A'_0)$;

f) $A_0, A_i, A'_i, A'_0$,  $(\tilde{A}_i=A_i, \tilde{A}'_i=A'_i)$.

\begin{proof}[Proof of Claim \ref{cl2.0} for $p_0\in\partial P_{L_1^{-1}\otimes L_2}$]

If $p_0\in e_{\alpha+\beta}(P_{L_1^{-1}\otimes L_2})$, then there is nothing to prove as $S$ is not contained in $P_{L_2}$. If $p_0\in \bigcup_{\rho\in-\sigma\cap\Sigma_X(1)}e_{\rho}(P_{L_1^{-1}\otimes L_2})$ and not an endpoint of this polygonal chain, then one sees easily $\iota(\alpha_i)$ and $\iota(\beta_j)$ are well defined for all $0\leq i\leq k$ and $0\leq\beta\leq l$. Then the arguement above for $p_0\notin\partial P_{L_1^{-1}\otimes L_2}$ can still serve here. Next we will prove the case when $p_0$ is contaiend in $\bigcup_{0\leq j\leq l}e_{\beta_j}(P_{L_1^{-1}\otimes L_2})$. The case when $p_0\in\bigcup_{0\leq i\leq k}e_{\alpha_i}(P_{L_1^{-1}\otimes L_2})$ can be proved similarly and we omit it. If $p_0\in e_{\beta}(P_{L_1^{-1}\otimes L_2})\backslash e_{\alpha+\beta}(P_{L_1^{-1}\otimes L_2})$, then we have $S\subset \bl{u}_h+P_{L_1}$. Besides, if $\|e_{\beta}(P_{L_1^{-1}\otimes L_2})\|=0$ and $p_0$ lie at the top position of the $\bigcup_{1\leq j\leq l}e_{\beta_j}(P_{L_1^{-1}\otimes L_2})$, then there will be nothing to prove since $S$ is not contained in $P_{L_2}$. Furthermore, if $p_0$ lies at the lowest position of $P_{L_1^{-1}\otimes L_2}$, one sees easily (\ref{k1}) and (\ref{l1}) are valid for $p_0$ automatically, hence our previous proof will provide a translation vector with which we can cover $S$ as indicated by (\ref{1i1}).  

By the argument above, we can find some $1\leq j\leq l$ and $p_1\in e_{\beta_j}(P_{L_1^{-1}\otimes L_2})$ such that $\overrightarrow{p_0p_1}=\bl{u}_{\beta_j}$, as is shown in the figure on the left below.  If the open quadrant on the lower right of $p_0$ contains a lattice point of $P_{L_1^{-1}\otimes L_2}$, then one checks easily $p_0+\bl{u}_h\in P_{L_1\\^{-1}\otimes L_2}$ hence we are done. Next we will assume all lattice points of $P_{L_1^{-1}\otimes L_2}$ below $l_{\bl{u}_h}(p_0)$ are contained in the cone bounded by $r_{-\bl{u}_h}(p_0)$ and $r_{-\bl{u}_v}(p_0)$.

Let $R$ be the open angular region bounded by $r_{\overrightarrow{p_0p_1}}(p_0)$ and $r_{\bl{u}_h}(p_0)$, then we claim there must be some lattice point $p$ in $P_{L_1^{-1}\otimes L_2}\cap R$ such that $\overrightarrow{p_0p}$ is of type c) or d) or e) if $\overrightarrow{p_0p_1}$ is not of type c). If the contrary is true, then for all $p\in P_{L_1^{-1}\otimes L_2}\cap R\cap M$ the type of $\overrightarrow{p_0p}$ must be one of the following: a), b), f), g) and h).

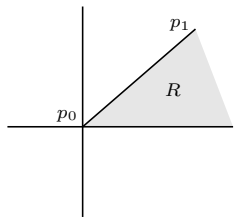
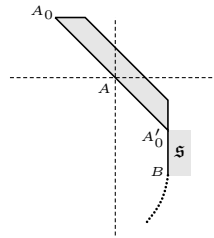
\begin{figure}[H]
\begin{tikzpicture}
\draw[ line width=0.6pt] (-1,0)--(-4,0);
\draw[ line width=0.6pt] (-3,1.6)--(-3,-1.2);
\draw[line width=0.6pt] (-3,0)--(-1.5,1.3);

\fill [fill opacity=0.1](-3,0)--(-1.5,1.3)--(-1,0)--(-3,0);

\node at (-1.7,1.35)  (0) {\relsize{-2}$p_1$}; 
\node at (-3.2,0.15)  (0) {\relsize{-2}$p_0$};

\node at (-1.8,0.5)  (0) {\relsize{-2}$R$}; 

----------------------------------------------------------------------------



\draw[ line width=0.6pt](3.6, 1.6)--(4,1.6)--(5.1,0.5)--(5.1,0.1)--(3.6, 1.6);
\fill [fill opacity=0.1] (3.6, 1.6)--(4,1.6)--(5.1,0.5)--(5.1,0.1)--(3.6, 1.6);


\draw[ line width=0.6pt](5.1,0.1)--(5.1, -0.5);
\fill [fill opacity=0.1] (5.1, 0.1)--(5.1,-0.5)--(5.4,-0.5)--(5.4,0.1)--(5.1, 0.1);


\draw [dash pattern=on 1.5pt off 1pt](4.4, -1.3)--(4.4,1.6);
\draw [dash pattern=on 1.5pt off 1pt](3, 0.8)--(5.7,0.8);


\draw[line width=1pt, line cap=round, dash pattern=on 0pt off 1.5\pgflinewidth](5.1,-0.5) .. controls (5.05,-0.9) and (4.9,-1) .. (4.8,-1.15);


\node at (3.43,1.67)  (0) {\relsize{-5}$A_0$};

\node at (4.9,0)  (0) {\relsize{-5}$A'_0$};

\node at (4.25,0.65)  (0) {\relsize{-5}$A$}; 

\node at (4.96,-0.45)  (0) {\relsize{-5}$B$}; 

\node at (5.26,-0.18)  (0) {\relsize{-0.5}$\mathfrak{s}$}; 

\end{tikzpicture}

\end{figure}

If there exists at least one lattice point $p$ such that $\overrightarrow{p_0p}$ is of type b) or f), then by defintion on the segment $\overline{A_0A'_0}$ there lies some contact points of $P_{L_1}$ with respect to $\overrightarrow{p_0p}$. Since $P_{L_1^{-1}\otimes L_2}\cap R$ contains only finitely many lattice points, we can pick one such contact point with the lowest position on $\overline{A_0A'_0}$. We will denote this point by $A$. Otherwise, we will just let $A=A_0$. Note that in either case we have $A\neq A'_0$. Now we take $p\in P_{L_1^{-1}\otimes L_2}\cap R\cap M$ such that $\overrightarrow{p_0p}$ is not of type g) or h), then by our choice of $A$, the intersection of $\overrightarrow{p_0p}+P_{L_1(D_{\sigma})}$ and the open cone which is bounded by $r_{\bl{u}_h}(A)$ and $r_{-\bl{u}_v}(A)$ is empty. 

Let $B$ be the lower endpoint of $e_{\alpha}(P_{L_1})$ and $\mathfrak{s}$ be the interior of the convex hull of $A'_0, \epsilon\bl{u}_h+A'_0, B$ and $\epsilon\bl{u}_h+B$ for some $\epsilon>0$. Since by condition $p_0$ is not the lowest point on $P_{L_1^{-1}\otimes L_2}$, for all sufficiently small $\epsilon>0$ we have $\mathfrak{s}\subset\epsilon\bl{u}_h+P_{L_1(D_{\sigma})}\subset P_{L_2(D_{\sigma})}$.


Next we will show $\mathfrak{s}$ cannot be covered by integral translations of $P_{L_1(D_{\sigma})}$ that are contained in $P_{L_2(D_{\sigma})}$ when $\epsilon$ is sufficiently small. First of all, by our argument above for all lattice point $p$ such that $p$ lies below the line $l_{\bl{u}_h}(p_0)$ or $\overrightarrow{p_0p}$ of type a) b) and f) 
we always have $(\overrightarrow{p_0p}+P_{L_1(D_{\sigma})})\cap\mathfrak{s}=\emptyset$. The same result is also true for a lattice ponit $p$ such that $\overrightarrow{p_0p}$ is of type h) or type g) with the only contact point not contained in $e_{\alpha}(P_{L_1})$ as long as $\epsilon>0$ is suffiently small as there are only finitely many such lattice points. As for a lattice point $p$ such that $\overrightarrow{p_0p}$ is of type g) and the contact point lies on $e_{\alpha}(P_{L_1})$ one sees easily the intersection $(\overrightarrow{p_0p}+P_{L_1(D_{\sigma})})\cap\mathfrak{s}$ is contained in a neighborhood of the contact point whose size will shrink with $\epsilon$. Again as there are only finitely many such lattice points some points in the two-dimensional region $\mathfrak{s}$ will be missed out by these intersections for suffciently small $\epsilon$. Then our claim is proved and we get a contradiction since by inductive assumption we have 
\[P_{L_2(D_{\sigma})}\subseteq \bigcup_{p\in P_{L_1^{-1}\otimes L_2}\cap M}\overrightarrow{p_0p}+P_{L_1(D_{\sigma})}. 
\]
As a consequence, we can find $p\in P_{L_1^{-1}\otimes L_2}\cap R\cap M$ such that $\overrightarrow{p_0p}$ is of type c) or d) or e). Moreover, note that as $\|e_{\beta_j}(P_{L_1^{-1}\otimes L_2})\|\geq \|\bl{u}_{\beta_j}\|$, we have 
\[\|e_{\beta_j}(P_{L_1})\|\geq 2\|\bl{u}_{\beta_j}\|.
\]
If $\overrightarrow{p_0p}$ is not of type c), then by a similar argument as the proof of (\ref{fro}) we can show
\[\mathcal{L}_{\overrightarrow{p_0x}}(P_{L_1})>\|\overrightarrow{p_0x}\|,
\]
 for any lattice point $x$ in the triangle $p_0p_1p$. In particular, $P_{L_1}$ has two contact points with respect to $\overrightarrow{p_0x}$ and $\overrightarrow{p_0x}$ is not of type f). 

\end{proof}

\medskip\emph{Step. 2.2. Impossible Types of Neighboring Translation Vectors.}

In the sequel we will use  x)$_i$-y)$_{i+1}$ to indicate that the consecutive translation vectors $\overrightarrow{p_0p_i}$ and  $\overrightarrow{p_0p_{i+1}}$ are of type x) and y) respectively. Besides, for $1\leq i\leq m$ we will take $B_i=A_i+\overrightarrow{p_ip_0}$ and $B'_i=A'_i+\overrightarrow{p_ip_0}$, then $B_i, B'_i\in\partial P_{L_1}$ as $A_i, A'_i$ are contact points of $\overrightarrow{p_0p_i}$.  For  computation in our proof, we will assume \[\overrightarrow{p_0p_{i+1}}=a\bl{u}_h+b\bl{u}_v,\;\;\; \overrightarrow{p_0p_i}=c\bl{u}_h+d\bl{u}_v.\] Then it is easy to see we have the following equivalence
\begin{equation}\label{abcd}\|\overrightarrow{p_0p_i}\|>\|\overrightarrow{p_0p_{i+1}}\|\Leftrightarrow c\geq a, d>b,\;\;\; 
\|\overrightarrow{p_0p_i}\|<\|\overrightarrow{p_0p_{i+1}}\|\Leftrightarrow c< a, d\leq b.
\end{equation}

Next we will prove the following 
\begin{claim}\label{itnv}  \normalfont For $1\leq i\leq m-1$ if one of $\overrightarrow{p_0p_i}, \overrightarrow{p_0p_{i+1}}$ is of type d) or e), so is the other.

\end{claim}

Assuming Claim \ref{itnv} is proved, then there must be a translation vector of type c). Firstly, the possiblity of being type f) for these translation vectors are ruled out by (\ref{fro}). 
Then by the proof of (\ref{2cps}) and particularly (\ref{k1}) and (\ref{l1})  one of $\overrightarrow{p_0p_1}$ and $\overrightarrow{p_0p_m}$ is of type d) or e) and the other of type a) or b).  However, by Claim \ref{itnv}, there will be no translation vectors of types other than d) and e), a contradiction. Therefore,  (\ref{1i1}) can be deduced from  our definition of type c).

\bigskip
\emph{Proof of Claim \ref{itnv}.} 

\smallskip
2.2.1.1  The case e)$_i$-a)$_{i+1}$.

Following the notations above, if the contrary is true, $A_iA'_iB'_iB_i$ lies on the open right half plane cut out by $l_{\overrightarrow{p_0p_i}}(A'_0)$  while $A_{i+1}A'_{i+1}B'_{i+1}B_{i+1}$  lies on the left of $l_{\overrightarrow{p_0p_{i+1}}}(A'_0)$. Since the slope of $\overrightarrow{p_0p_i}$ is larger than that of $\overrightarrow{p_0p_{i+1}}$, the two parallelograms are separate. By Lemma \ref{ae}, their vertices are contained in a pair of paralleling lines.
Moreover, each of the two lines contains four points of $\partial P_{L_1}$ hence an edge of $P_{L_1}$. Note that $\overrightarrow{p_0p_i}$, $\overrightarrow{p_0p_{i+1}}$ cannot both parallel with these lines. Therefore, we have at least a pair of points, say $A_i$ and $A'_i$, lying on one of these lines. However, that would imply one of $A_i$ and $A'_i$ is not a contact point as these two parallelograms have empty intersection, a contradiction.

2.2.1.2  The case e)$_i$-b)$_{i+1}$. 

The proof is the same as the case e)$_i$-a)$_{i+1}$ hence omitted.

\smallskip
2.2.1.3 The case a)$_i$-e)$_{i+1}$.

Let $B=A_0+\overrightarrow{p_ip_0}, B'=A'_0+\overrightarrow{p_{i+1}p_0}$, then by combining  $\|\overline{A_0A'_0}\|\geq \sqrt2$ with the fact that $\overrightarrow{p_0p_i}, \overrightarrow{p_0p_{i+1}}$ form a basis of $M$ (these vectors are neither horizontal nor vertical by (\ref{k1}) and (\ref{l1}) one checks easily $\overline{A_0B}, \overline{A'_0B'}$ have no common points. By condition $\mathcal{L}_{\overrightarrow{p_0p_i}}(A_0)<\|\overrightarrow{p_0p_i}\|$, the chord of $P_{L_1}$ through $A_0$ and paralleling with $\overrightarrow{p_0p_i}$ is contained in $\overline{A_0B}$. Similarly, the chord of $P_{L_1}$ through $A'_0$ paralleling with $\overrightarrow{p_0p_{i+1}}$ is contained in $\overline{A'_0B'}$. We will denote these chords by $\overline{A_0B_0}$ and $\overline{A'_0B'_0}$ respectively.

Let $\tilde{A}'_i=A'_iB'_i\cap A_0A'_0, \tilde{A}_{i+1}=A_{i+1}B_{i+1}\cap A_0A'_0$ and $C=A'_iB'_i\cap A_{i+1}B_{i+1}$. Since $\tilde{A}'_i$ (resp. $\tilde{A}_{i+1}$) lies on the upper left (resp lower right) of $A_0$ (resp. $A'_0$), $C$ is on the lower left of $A_0A'_0$. The following figure lists all the cases of the relative positions of the chords $\overline{A'_iB'_i}$ and $\overline{A_{i+1}B_{i+1}}$ on the union $A'_iB'_i\cup A_{i+1}B_{i+1}$. 

For the first one, i.e. when one of $\overline{A_iB_i}$ and $\overline{A'_{i+1}B'_{i+1}}$ containing $C=A_iB_i\cap A'_{i+1}B'_{i+1}$, in order that the four endpoints of these chords are contained on the boundary of their convex hull, $C$ must be contained in both of them. Then $B_0, B'_0$ are contained in the interior of the convex hull of $A_0, A'_0, A_i, A'_{i+1}$ and $C$, which is impossible. The proofs for the other cases are similar and we omit them. Therefore the types of $\overrightarrow{p_0p_i}, \overrightarrow{p_0p_{i+1}}$ cannot be a), e). 

\begin{figure}[H]
\begin{tikz1}

\draw (-6.8,0.8)--(-6,0)--(-5.2,-0.8);

\draw [dash pattern=on 1.5pt off 1pt](-6.55,0.55)--(-7.13,-3.51);
\draw [dash pattern=on 1.5pt off 1pt](-5.45,-0.55)--(-8.04,-2.4);

\draw [line width=1pt] (-6.75,-0.85)--(-6.95,-2.25);

\draw [line width=1pt] (-6.43,-1.25)--(-7.13,-1.75);

\draw (-6.25,0.25)--(-6.35,-0.45);
\draw (-5.75,-0.25)--(-6.1,-0.5);

\node at(-6.1,0.32) (0) {\relsize{-10}$A_0$};
\node at(-6.42,-0.55) (0) {\relsize{-10}$B_0$};
\node at(-5.6,-0.18) (0) {\relsize{-10}$A'_0$};
\node at(-6.03,-0.62) (0) {\relsize{-10}$B'_0$};
\node at(-6.4,0.6) (0) {\relsize{-10}$\tilde{A}'_i$};
\node at(-5.15,-0.48) (0) {\relsize{-10}$\tilde{A}_{i+1}$};
\node at(-6.73,-1.65)(0) {\relsize{-10}$C$};

\node at(-6.92,-0.85) (0) {\relsize{-10}$A'_i$};
\node at(-6.8,-2.3) (0) {\relsize{-10}$B'_i$};
\node at(-6.1,-1.32) (0) {\relsize{-10}$A_{i+1}$};
\node at(-7.45,-1.65)(0) {\relsize{-10}$B_{i+1}$};

-------------------------------------------------------------
\draw (-3.8,0.8)--(-3,0)--(-2.2,-0.8);
\draw [dash pattern=on 1.5pt off 1pt](-3.55,0.55)--(-4.13,-3.51);
\draw [dash pattern=on 1.5pt off 1pt](-2.45,-0.55)--(-5.04,-2.4);

\draw [line width=1pt] (-3.6,0.2)--(-3.8,-1.2);

\draw [line width=1pt] (-2.73,-0.75)--(-3.43,-1.25);

\draw (-3.25,0.25)--(-3.35,-0.45);
\draw (-2.75,-0.25)--(-3.1,-0.5);

\node at(-3.1,0.32) (0) {\relsize{-10}$A_0$};
\node at(-3.42,-0.55) (0) {\relsize{-10}$B_0$};
\node at(-2.6,-0.18) (0) {\relsize{-10}$A'_0$};
\node at(-3.03,-0.62) (0) {\relsize{-10}$B'_0$};
\node at(-3.4,0.6) (0) {\relsize{-10}$\tilde{A}'_i$};
\node at(-2.15,-0.48) (0) (0) {\relsize{-10}$\tilde{A}_{i+1}$};
\node at(-3.73,-1.65)(0) {\relsize{-10}$C$};

\node at(-3.8,0.2) (0) {\relsize{-10}$A'_i$};
\node at(-4,-1.18) (0) {\relsize{-10}$B'_i$};
\node at(-2.46,-0.87) (0) {\relsize{-10}$A_{i+1}$};
\node at(-3.15,-1.4)(0) {\relsize{-10}$B_{i+1}$};
-------------------------------------------------------------

\draw (-0.8,0.8)--(0,0)--(0.8,-0.8);
\draw [dash pattern=on 1.5pt off 1pt](-0.55,0.55)--(-1.13,-3.51);
\draw [dash pattern=on 1.5pt off 1pt](0.55,-0.55)--(-2.04,-2.4);

\draw[line width=1pt] (-0.9,-1.9)--(-1.1,-3.3);
\draw[line width=1pt] (-1.13,-1.75)--(-1.83,-2.25);

\draw (-0.25,0.25)--(-0.35,-0.45);
\draw (0.25,-0.25)--(-0.1,-0.5);

\node at(-0.1,0.32) (0) {\relsize{-10}$A_0$};
\node at(-0.42,-0.55) (0) {\relsize{-10}$B_0$};
\node at(0.4,-0.18) (0) {\relsize{-10}$A'_0$};
\node at(-0.03,-0.62) (0) {\relsize{-10}$B'_0$};
\node at(-0.4,0.6) (0) {\relsize{-10}$\tilde{A}'_i$};
\node at(0.85,-0.48)(0) {\relsize{-10}$\tilde{A}_{i+1}$};
\node at(-0.73,-1.65)(0) {\relsize{-10}$C$};

\node at(-0.73,-1.9) (0) {\relsize{-10}$A'_i$};
\node at(-0.93,-3.35) (0) {\relsize{-10}$B'_i$};

\node at(-1.4,-1.57) (0) {\relsize{-10}$A_{i+1}$};
\node at(-2.16,-2.1) (0) {\relsize{-10}$B_{i+1}$};

-------------------------------------------------------------

\draw (2.2,0.8)--(3,0)--(3.8,-0.8);
\draw [dash pattern=on 1.5pt off 1pt](2.45,0.55)--(1.87,-3.51);
\draw [dash pattern=on 1.5pt off 1pt](3.55,-0.55)--(0.96,-2.4);

\draw [line width=1pt] (2.4,0.2)--(2.2,-1.2);
\draw[line width=1pt] (1.87,-1.75)--(1.17,-2.25);

\draw (2.75,0.25)--(2.65,-0.45);
\draw (3.25,-0.25)--(2.9,-0.5);

\node at(2.9,0.32) (0) {\relsize{-10}$A_0$};
\node at(2.58,-0.55) (0) {\relsize{-10}$B_0$};
\node at(3.4,-0.18) (0) {\relsize{-10}$A'_0$};
\node at(2.97,-0.62) (0) {\relsize{-10}$B'_0$};
\node at(2.6,0.6) (0) {\relsize{-10}$\tilde{A}'_i$};
\node at(3.85,-0.48)(0) {\relsize{-10}$\tilde{A}_{i+1}$};
\node at(2.27,-1.65)(0) {\relsize{-10}$C$};

\node at(2.2,0.2) (0) {\relsize{-10}$A'_i$};
\node at(2,-1.18) (0) {\relsize{-10}$B'_i$};

\node at(1.6,-1.57) (0) {\relsize{-10}$A_{i+1}$};
\node at(0.84,-2.1) (0) {\relsize{-10}$B_{i+1}$};
-------------------------------------------------------------
\draw (5.2,0.8)--(6,0)--(6.8,-0.8);

\draw [dash pattern=on 1.5pt off 1pt](5.45,0.55)--(4.87,-3.51);
\draw [dash pattern=on 1.5pt off 1pt](6.55,-0.55)--(3.96,-2.4);

\draw [line width=1pt] (6.27,-0.75)--(5.57,-1.25);
\draw[line width=1pt] (5.1,-1.9)--(4.9,-3.3);

\draw (5.75,0.25)--(5.65,-0.45);
\draw (6.25,-0.25)--(5.9,-0.5);

\node at(5.9,0.32) (0) {\relsize{-10}$A_0$};
\node at(5.58,-0.55) (0) {\relsize{-10}$B_0$};
\node at(6.4,-0.18) (0) {\relsize{-10}$A'_0$};
\node at(5.97,-0.62) (0) {\relsize{-10}$B'_0$};
\node at(5.6,0.6) (0) {\relsize{-10}$\tilde{A}'_i$};
\node at(6.85,-0.48)(0)(0) {\relsize{-10}$\tilde{A}_{i+1}$};
\node at(5.27,-1.65)(0) {\relsize{-10}$C$};

\node at(5.27,-1.9) (0) {\relsize{-10}$A'_i$};
\node at(5.07,-3.35) (0) {\relsize{-10}$B'_i$};

\node at(6.54,-0.87) (0) {\relsize{-10}$A_{i+1}$};
\node at(5.85,-1.4)(0) {\relsize{-10}$B_{i+1}$};
-------------------------------------------------------------

\end{tikz1}
\caption{The case a)$_i$-e)$_{i+1}$. } \label{fig:M1}

\end{figure}

2.2.1.4 The case b)$_i$-e)$_{i+1}$. 

Firstly we consider the case $\|\overrightarrow{p_0p_i}\|<\|\overrightarrow{p_0p_{i+1}}\|$. Let $C=l_{\overrightarrow{p_0p_{i+1}}}(B'_i)\cap (A_0A'_0\cup l_{\bl{u}_v}(A'_0))$, then we claim $C$ is not lying on $A_0A'_0$. Otherwise, let $D$ be the point on $\partial P_{L_1}$ such that $\overline{A_0'D}$ is the chord through $A_0'$ paralleling with $\overrightarrow{p_0p_{i+1}}$. Since $\overrightarrow{p_0p_{i+1}}$ is of type e), we have $\|\overline{A'_0D}\|=\mathcal{L}_{\overrightarrow{p_0p_{i+1}}}(A'_0)>\mathcal{L}_{\overrightarrow{p_0p_{i+1}}}(C)\geq\|\overline{B_i'C}\|$. Let $E=DB'\cap A_0A_0'$, then $E$ is on the upper left side of $A'_i$. Indeed, if $E$ lies strictly on the lower right of $A'_i$, $B'_i$ would be contained in the interior of the triangle $A'_0A'_iD$, which is impossible. Besides, $E=A_i'$ would imply $A'_i, B'_i, D$ colinear, then $\overline{A'_iB'_i}$ is contained in  an edge of $P_{L_1}$, which is again impossible.  Now one sees easily $\mathcal{L}_{\overrightarrow{p_0p_i}}(A'_i+\epsilon(\bl{u}_h-\bl{u}_v))>\|\overline{A'_iB'_i}\|=\|\overrightarrow{p_0p_i}\|$ for any sufficiently small $\epsilon>0$, which contradicts with the fact that $\overrightarrow{p_0p_i}$ is of type b). 

Let $C=l_{\overrightarrow{p_0p_{i+1}}}(B_i')\cap r_{-\bl{u}_v}(A'_0)$, one checks easily $C$ is contained in the edge $e_{\alpha}(P_{L_1})$ and $\|\overline{B_i'C}\|<\|\overrightarrow{p_0p_{i+1}}\|$. On the other hand, since $\overrightarrow{p_0p_i}$ is of type b), we have $A'_i\neq A_0$, let $F=A_i'+\epsilon(-\bl{u}_h+\bl{u}_v)$ for sufficiently small $\epsilon>0$ and $G=F+\overrightarrow{p_ip_0}$. Since $\mathcal{L}_{\overrightarrow{p_0p_i}}(F)\geq\|\overrightarrow{p_0p_i}\|$, $G$ is contained in the interior of $P_{L_1}$ and when $\epsilon$ is sufficiently small the intersection of $l_{\overrightarrow{p_0p_{i+1}}}(G)$ with $l_{\bl{u}_v}(A'_0)$ lies between $A'_0$ and $C$, as is shown in figure (ii) below.  Then one sees easily $\mathcal{L}_{\overrightarrow{p_0p_{i+1}}}(G)>\mathcal{L}_{\overrightarrow{p_0p_{i+1}}}(B_i')=\|\overline{B_i'C}\|$. Therefore, the contact points $A_{i+1}, A'_{i+1}$ of $P_{L_1}$ in the direction of $\overrightarrow{p_0p_{i+1}}$ are contained in $\overline{A_0'C}$. 

Note that $\|\overline{A'_{i+1}B'_{i+1}}\|>\|\overline{B'_iC}\|$, therefore $\overline{A'_{i+1}B'_{i+1}}\cap \overline{A'_iB'_i}\neq\emptyset$ and the common point is not an endpoint. Let $H_{i+1}=A_0'+\overrightarrow{p_{i+1}p_0}$, $H_i=A_0'+\overrightarrow{p_ip_0}$, $I=H_{i+1}H_i\cap l_{\bl{u}_v}(A_0')$ and $J=B'_iB'_{i+1}\cap l_{\bl{u}_v}(A_0')$, then one sees easily $J$ is lying between $A_0'$ and $I$. Moreover, one computes $\|\overline{A_0'I}\|=\tfrac{1}{a-c}\leq 1$, which implies $\overline{A'_0I}\subseteq e_{\alpha}(P_{L_1})$. Then as $J, B'_{i+1}$ and $B'_i$ are also contained in $\partial P_{L_1}$, $\overline{B'_{i+1}B'_i}$ is contained in some edge of $P_{L_1}$. However, this would imply $I$ lies outside $P_{L_1}$, a contradiction.

\begin{figure}[H]
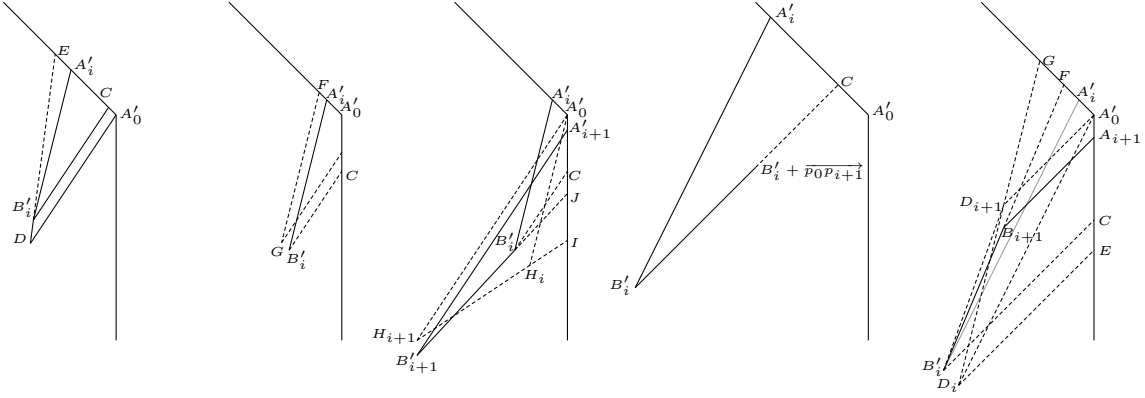

\begin{tikz1}

\draw (-6,-3)--(-6,0)--(-7.5,1.5);

\draw (-6.6,0.6)--(-7.1,-1.4);

\draw(-7.1, -1.4)--(-6.1,0.1);
\draw (-7.14,-1.71)--(-6,0);
\draw (-7.14,-1.71)--(-7.1,-1.4);

\draw [dash pattern=on 1.5pt off 1pt]  (-7.1,-1.4)--(-6.81,0.81);

\node at(-5.8,0) {\relsize{-20}$A'_0$};
\node at(-6.13,0.3) {\relsize{-20}$C$};
\node at(-7.25,-1.23) {\relsize{-20}$B'_i$};
\node at(-6.42,0.65) {\relsize{-20}$A'_i$};
\node at(-6.7,0.86) {\relsize{-20}$E$};
\node at(-7.28,-1.64) {\relsize{-20}$D$};

------------------------------------------------------
\draw (-3,-3)--(-3,0)--(-4.5,1.5);
\draw (-3.2,0.2)--(-3.7,-1.8);
\draw[dash pattern=on 1.5pt off 1pt] (-3.3,0.3)--(-3.8,-1.7);
\draw[dash pattern=on 1.5pt off 1pt] (-3.7,-1.8)--(-3, -0.75);
\draw[dash pattern=on 1.5pt off 1pt] (-3.8,-1.7)--(-3, -0.5);

\node at(-2.88,0.1) {\relsize{-20}$A'_0$};
\node at(-3.08,0.27) {\relsize{-20}$A'_i$};
\node at(-3.25,0.41) {\relsize{-20}$F$};
\node at(-3.6,-1.93) {\relsize{-20}$B'_i$};
\node at(-2.86,-0.8) {\relsize{-20}$C$};
\node at(-3.86,-1.82) {\relsize{-20}$G$};
-----------------------------------------------------
\draw (0,-3)--(0,0)--(-1.5,1.5);
\draw [dash pattern=on 1.5pt off 0.8pt](0,0)--(-2,-3);
\draw[dash pattern=on 1.5pt off 1pt] (0,0)--(-0.5,-2);
\draw[dash pattern=on 1.5pt off 1pt](-2,-3)--(-0.5,-2)--(0,-1.67);

\draw (0,-0.2)--(-2,-3.2);
\draw (-0.2,0.2)--(-0.7,-1.8);
\draw (-2, -3.2)--(-0.7,-1.8);
\draw [dash pattern=on 1.5pt off 0.8pt](-0.7,-1.8)--(0,-1.046);
\draw[dash pattern=on 1.5pt off 1pt] (-0.7,-1.8)--(0, -0.75);

\node at(0.12,0.1) {\relsize{-20}$A'_0$};
\node at(-0.08,0.27) {\relsize{-20}$A'_i$};
\node at(0.1,-1.7) {\relsize{-20}$I$};
\node at(0.1,-1.1) {\relsize{-20}$J$};
\node at(0.1,-0.8) {\relsize{-20}$C$};
\node at(-0.83,-1.69) {\relsize{-20}$B'_i$};
\node at(-0.43,-2.15) {\relsize{-20}$H_i$};
\node at(-2.3,-2.95) {\relsize{-20}$H_{i+1}$};
\node at(0.3,-0.2) {\relsize{-20}$A'_{i+1}$};
\node at(-2,-3.3) {\relsize{-20}$B'_{i+1}$};

------------------------------------------------------------------------
\draw (4,-3)--(4,0)--(2.5,1.5);
 \draw(2.7,1.3)--(0.9,-2.3);
 \draw(0.9,-2.3)--(2.5,-0.7);
 \draw[dash pattern=on 1.5pt off 1pt](2.5,-0.7)--(3.6,0.4);
\node at(2.9,1.35) {\relsize{-20}$A'_i$};
\node at(0.7,-2.27) {\relsize{-20}$B'_i$};
\node at(3.25,-0.75) {\relsize{-20}$B'_i+\overrightarrow{p_0p_{i+1}}$};
\node at(3.72,0.44) {\relsize{-20}$C$};
\node at(4.2,0.05) {\relsize{-20}$A'_0$};
---------------------------------------------------------------------------------------------------------------------------
\draw (7,-3)--(7,0)--(5.5,1.5);
\draw [dash pattern=on 1.5pt off 0.8pt](7,0)--(5.8,-1.2);
\draw[dash pattern=on 1.5pt off 1pt] (7,0)--(5.2,-3.6);

\draw[dash pattern=on 1.5pt off 1pt]  (5.2,-3.6)--(5.8,-1.2)--(6.28,0.72);

\draw (7,-0.3)--(5.8,-1.5);

\draw [opacity=0.4]   (6.8,0.2)--(5,-3.4);
\draw[dash pattern=on 1.5pt off 1pt]  (5,-3.4)--(5.8,-1.2);
\draw (5,-3.4)--(5.8,-1.5);
\draw [dash pattern=on 1.5pt off 1pt](5.8,-1.5)--(6.6,0.4);
\draw [dash pattern=on 1.5pt off 1pt] (5,-3.4)--(7,-1.4);
\draw [dash pattern=on 1.5pt off 1pt] (5.2,-3.6)--(7,-1.8);

\node at(7.2,0) {\relsize{-20}$A'_0$};
\node at(7.35,-0.3) {\relsize{-20}$A_{i+1}$};
\node at(7.15,-1.4) {\relsize{-20}$C$};
\node at(6.05,-1.61){\relsize{-20}$B_{i+1}$};
\node at(4.85,-3.35) {\relsize{-20}$B'_i$};
\node at(6.4,0.73) {\relsize{-20}$G$};

\node at(5.05,-3.6) {\relsize{-20}$D_i$};
\node at(5.5,-1.2) {\relsize{-20}$D_{i+1}$};

\node at(7.15,-1.8) {\relsize{-20}$E$};
\node at(6.6,0.52) {\relsize{-20}$F$};
\node at(6.9,0.3) {\relsize{-20}$A'_i$};

\end{tikz1}

\caption{The case b)$_i$-e)$_{i+1}$. } \label{fig:M1}

\end{figure}

Next we consider the case $\|\overrightarrow{p_0p_i}\|>\|\overrightarrow{p_0p_{i+1}}\|$. Let $C=l_{\overrightarrow{p_0p_{i+1}}}(B'_i)\cap (A_0A'_0\cup l_{\bl{u}_v}(A'_0))$, then $C$ cannot be contained in $\overline{A_0A'_0}$, otherwise $\mathcal{L}_{\overrightarrow{p_0p_{i+1}}}(C)>\|\overrightarrow{p_0p_{i+1}}\|$, which contradicts with the condition that $\overrightarrow{p_0p_{i+1}}$ is of type e). 

By abuse of notation, let $C= l_{\overrightarrow{p_0p_{i+1}}}(B'_i)\cap l_{\bl{u}_v}(A'_0)$, $D_i=A'_0+\overrightarrow{p_ip_0}, D_{i+1}=A_0'+\overrightarrow{p_{i+1}p_0}$ and $E=l_{\overrightarrow{p_0p_{i+1}}}(D_i)\cap l_{\bl{u}_v}(A'_0)$. Then one computes easily $\|\overline{A'_0C}\|<\|\overline{A'_0E}\|=\tfrac{1}{a}\leq 1.$
Therefore $C$ is contained in $e_{\alpha}(P_{L_1})$. From the condition $\|\overrightarrow{p_0p_i}\|>\|\overrightarrow{p_0p_{i+1}}\|$ we can check easily $\|\overline{D_iE}\|\geq\|\overrightarrow{p_0p_{i+1}}\|$. Since $\overrightarrow{p_0p_i}$ is of type b), we have $A'_i\neq A'_0$, then $\mathcal{L}_{\overrightarrow{p_0p_{i+1}}}(C)=\|\overline{B'_iC}\|>\|\overline{D_iE}\|\geq\|\overrightarrow{p_0p_{i+1}}\|$. Then we can deduce $A_{i+1}$ lie between $A_0'$ and $C$ as $\overrightarrow{p_0p_{i+1}}$ is of type e). We claim $\overline{A_{i+1}B_{i+1}}$ and $\overline{A'_iB'_i}$ must intersect and the intersection points cannot be the endpoints of these segments. Indeed, since $B_{i+1}$ lies on the border of $P_{L_1}$, it cannot stay on the right of the line $A'_iB'_i$. Moreover, if it lies on $\overline{A'_iB'_i}$, since $\|\overline{B'_iC}\|>\|\overline{A_{i+1}B_{i+1}}\|$ one sees easily it must be different from $A'_i$ and $B'_i$. Then $\overline{A'_iB'_i}$ lies in the edge of $P_{L_1}$ containing $A'_i, B'_i$ and $B_{i+1}$, which  is impossible.

Let $F=B'_iB_{i+1}\cap A_0A'_0, G=D_iD_{i+1}\cap A_0A_0'$, then by our argument above $F$ lies on the upper left of $A'_i$ and $G$ lies on the upper left of $F$. Moreover, one computes 
$\|\overline{GA'_0}\|=\tfrac{\sqrt{2}}{c+d-a-b}\leq \sqrt{2},
$
from which we can deduce $F$ and $G$ are contained in $e_{\alpha+\beta}(P_{L_1})$. On the other hand, the line $B'_iB_{i+1}$ contains $F\in\partial P_{L_1}$, hence this line contains an edge of $P_{L_1}$. Then $G$ would lie in the exterior of $P_{L_1}$, a contradiction.

\medskip

2.2.2.1 The case d)$_i$-a)$_{i+1}$.

The proof is similar to that of 2.2.1.1 and is omitted.

2.2.2.2 The case d)$_i$-b)$_{i+1}$. 

First assuming $\|\overrightarrow{p_0p_i}\|<\|\overrightarrow{p_0p_{i+1}}\|$, in the figure below we list all the possibilites of the relative positions of $\overline{A_iB_i}$ and $\overline{A'_{i+1}B'_{i+1}}$. For the first two, we have $\mathcal{L}_{\overrightarrow{p_0p_i}}(A_i+\epsilon(\bl{u}_v-\bl{u}_h))>\|\overrightarrow{p_0p_i}\|$ for all $\epsilon>0$ sufficiently small, which contradicts with the fact that $\overrightarrow{p_0p_i}$ is of type d).  For the third one, we have $\mathcal{L}_{\overrightarrow{p_0p_{i+1}}}(A_0)\geq\|\overrightarrow{p_0p_{i+1}}\|=\|\overline{A'_{i+1}B'_{i+1}}\|$, hence $C=A_0+\overrightarrow{p_{i+1}p_0}$ is contained in $P_{L_1}$. However, that would imply $B_i$ is contained in the interior of the parallelogram $A_0CB'_{i+1}A'_{i+1}$, which  contradicts with the fact that $B_i\in\partial P_{L_1}$. The case when $\|\overrightarrow{p_0p_i}\|>\|\overrightarrow{p_0p_{i+1}}\|$ can be proved similarly and we omit it. 

  \begin{figure}[H]
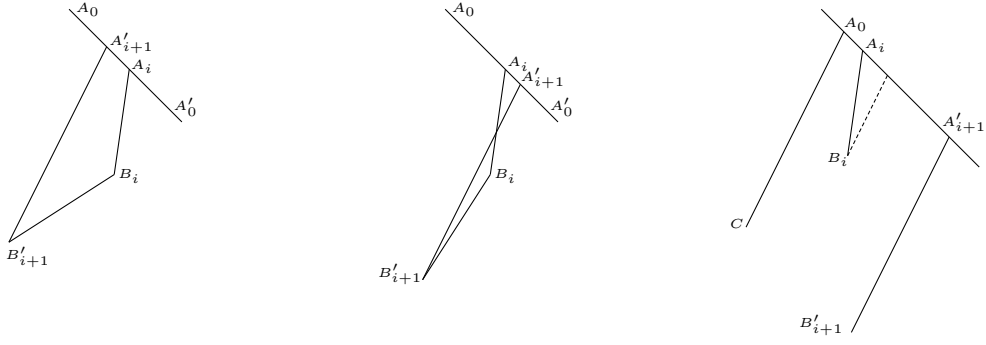


\begin{tikz1}
\draw (-5.1,0.1)--(-3.6,-1.4);

\draw (-4.6,-0.4)--(-5.9,-3);
\draw (-4.3,-0.7)--(-4.5,-2.1);

\draw(-5.9,-3)--(-4.5, -2.1);

\node at (-4.85,0.1) {\relsize{-20}$A_0$};
\node at (-4.27,-0.35) {\relsize{-20}$A'_{i+1}$};
\node at (-3.55,-1.2) {\relsize{-20}$A'_0$};
\node at (-4.13,-0.66) {\relsize{-20}$A_i$};
\node at (-4.3, -2.12) {\relsize{-20}$B_i$};
\node at (-5.65,-3.2) {\relsize{-20}$B'_{i+1}$};












-------------------------------------------------------------------------------------------

\draw (-0.1,0.1)--(1.4,-1.4);
\draw (0.7,-0.7)--(0.5,-2.1);

\draw (0.9,-0.9)--(-0.4,-3.5);
\draw (-0.4,-3.5)--(0.5,-2.1);

\node at (0.15,0.1) {\relsize{-20}$A_0$};
\node at (0.87,-0.63) {\relsize{-20}$A_i$};
\node at (1.2,-0.82) {\relsize{-20}$A'_{i+1}$};
\node at (1.45,-1.2) {\relsize{-20}$A'_0$};
\node at (-0.7,-3.43) {\relsize{-20}$B'_{i+1}$};
\node at (0.7,-2.12) {\relsize{-20}$B_i$};

------------------------------------------------------------------------

\draw (4.9,0.1)--(7,-2);
\draw (5.2,-0.2)--(3.9,-2.8);

\draw (5.45,-0.45)--(5.25,-1.85);

\draw (6.6,-1.6)--(5.3,-4.2);

\draw[dash pattern=on 1.5pt off 1pt](5.25,-1.85)--(5.78,-0.78);

\node at (5.3565,-0.1) {\relsize{-20}$A_0$};
\node at (3.78,-2.75) {\relsize{-20}$C$};
\node at (5.62,-0.38) {\relsize{-20}$A_i$};
\node at (6.8,-1.42) {\relsize{-20}$A'_{i+1}$};
\node at (4.9,-4.1) {\relsize{-20}$B'_{i+1}$};
\node at (5.12,-1.9) {\relsize{-20}$B_i$};

---------------------------------------------------------------------------------------------------

\end{tikz1}

\caption{The case d)$_i$-b)$_{i+1}$, $\|\protect\overrightarrow{p_0p_i}\|<\|\protect\overrightarrow{p_0p_{i+1}}\|$.  } \label{fig:M1}
\end{figure}

\medskip
2.2.2.3 The case a)$_i$-d)$_{i+1}$. 

Firstly we assume $\|\overrightarrow{p_0p_i}\|<\|\overrightarrow{p_0p_{i+1}}\|$. Let $C=l_{\overrightarrow{p_0p_i}}(B_{i+1})\cap(A_0A'_0\cup l_{\bl{u}_h}(A_0))$, then $C\notin\overline{A_0A'_0}$ since otherwise we will have  $\mathcal{L}_{\overrightarrow{p_0p_i}}(A_0)>\|\overrightarrow{p_0p_i}\|$, which contradicts with the condition that $\overrightarrow{p_0p_i}$ is of type a). By abuse of notation, let $C=l_{\bl{u}_h}(A_0)\cap l_{\overrightarrow{p_0p_i}}(B_{i+1})$. Under the condition that $\overrightarrow{p_0p_i}, \overrightarrow{p_0p_{i+1}}$ form a basis of $M$, one checks easily $\|\overline{CA_0}\|<\tfrac{1}{d}\leq1$. Therefore $C$ is contained in the edge $e_{\beta}(P_{L_1})$. From $\|\overrightarrow{p_0p_i}\|<\|\overrightarrow{p_0p_{i+1}}\|$ we can deduce $\|\overline{B_{i+1}C}\|>\|\overrightarrow{p_0p_i}\|$, thus $A'_i$ is lying between $C$ and $A_0$ as $\overrightarrow{p_0p_i}$ is of type a). Moreover, $B'_i$ is on the boundary of $P_{L_1}$, it cannot be contained in the interior of the triangle $A_{i+1}B_{i+1}C$. If it is contained in $\overline{A_{i+1}B_{i+1}}$, then $\overline{A_{i+1}B_{i+1}}$ would be contained in an edge of $P_{L_1}$ as this chord contains three different points on the border of $P_{L_1}$, which is impossible.  Consequently, the two chords $\overline{A'_iB'_i}$ and $\overline{A_{i+1}B_{i+1}}$ must intersect and the intersection point is different from their endpoints.

Let $D=B_{i+1}B'_i\cap A_0A'_0$, then one sees easily it lies on the lower right of $A_{i+1}$. Let $E_i=A_0+\overrightarrow{p_ip_0}, E_{i+1}=A_0+\overrightarrow{p_{i+1}p_0}$,  $F=B_{i+1}E_i\cap A_0A'_0$ and  $G= E_{i+1}E_i\cap A_0A'_0$, then one sees easily the relative positions of $D, F$ and $G$ are just as pictured in the figure below. One computes $\|A_0G\|=\tfrac{\sqrt{2}}{a+b-c-d}\leq \sqrt{2}$, hence $G$ and $D$ are contained in the edge $e_{\alpha+\beta}(P_{L_1})$. In particular, this would imply $\overline{B_{i+1}D}$ is contained in an edge of $P_{L_1}$ with $D$ one of its endpoints. On the other hand, note that  $D$ lies strictly on the upper left of $F$ as $A'_i$ cannot coincide with $A_0$ since $\overrightarrow{p_0p_i}$ is of type a),  then we would deduce $G$ lies outside $P_{L_1}$, a contradiction.

\begin{figure}[H]
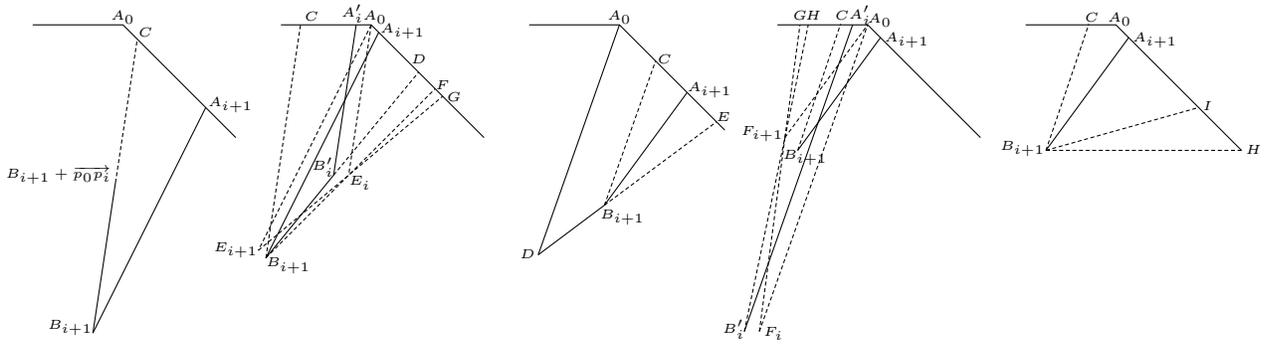

\begin{tikz1}

\draw (-4.5,0)--(-3.3,0)--(-1.8,-1.5);
\draw (-2.2,-1.1)--(-3.7,-4.1);
\draw (-3.7,-4.1)--(-3.4,-2.1);
\draw[dash pattern=on 1.5pt off 1pt](-3.4,-2.1)--(-3.11,-0.19);

\node at(-3.3,0.1) {\relsize{-10}$A_0$};
\node at(-3,-0.1) {\relsize{-10}$C$};
\node at (-1.87,-1.06){\relsize{-10}$A_{i+1}$};
\node at (-4,-4.03){\relsize{-10}$B_{i+1}$};
\node at (-4.16,-2){\relsize{-10}$B_{i+1}+\overrightarrow{p_0p_i}$};
------------------------------------------------------------------------------------
\draw (-1.2,0)--(0,0)--(1.5,-1.5);
\draw[dash pattern=on 1.5pt off 1pt](0,0)--(-1.5,-3);
\draw[dash pattern=on 1.5pt off 1pt](0,0)--(-0.3,-2);

\draw(-0.2,0)--(-0.5,-2);
\draw(0.1,-0.1)--(-1.4,-3.1);

\draw[dash pattern=on 1.5pt off 1pt](-1.5,-3)--(-0.3,-2)--(0.95,-0.95);
\draw[dash pattern=on 1.5pt off 1pt](-1.4,-3.1)--(-0.94,0);
\draw(-1.4,-3.1)--(-0.5,-2);
\draw[dash pattern=on 1.5pt off 1pt](-0.5,-2)--(0.63,-0.63);

\draw[dash pattern=on 1.5pt off 1pt](-1.4,-3.1)--(-0.3,-2)--(0.85,-0.85);

\node at(0.05,0.1) {\relsize{-10}$A_0$};
\node at (-0.15,-2.1){\relsize{-10}$E_i$};
\node at (-1.8,-3){\relsize{-10}$E_{i+1}$};
\node at (0.63,-0.45){\relsize{-10}$D$};
\node at (-0.79,0.12){\relsize{-10}$C$};
\node at(0.95,-0.8) {\relsize{-10}$F$};

\node at(-0.25,0.15) {\relsize{-10}$A'_i$};
\node at (-0.64,-1.89){\relsize{-10}$B_i'$};

\node at(0.42,-0.08) {\relsize{-10}$A_{i+1}$};
\node at (-1.1,-3.2){\relsize{-10}$B_{i+1}$};
\node at (1.1,-0.95){\relsize{-10}$G$};
---------------------------------------------------------------------
\draw(2.1,0)--(3.3,0)--(4.7,-1.4);
\draw(2.22,-3.06)--(3.3,0);
\draw (4.2,-0.9)--(3.1,-2.4);
\draw [dash pattern=on 1.5pt off 1pt] (3.1,-2.4)--(4.59,-1.29);
\draw [dash pattern=on 1.5pt off 1pt] (3.1,-2.4)--(3.78,-0.48);
\draw (2.22,-3.06)--(3.1,-2.4);
\node at(3.3,0.1) {\relsize{-10}$A_0$};
\node at (3.9,-0.45){\relsize{-10}$C$};
\node at (4.5,-0.87){\relsize{-10}$A_{i+1}$};
\node at (3.34,-2.56){\relsize{-10}$B_{i+1}$};
\node at (2.08,-3.04){\relsize{-10}$D$};
\node at (4.68,-1.22){\relsize{-10}$E$};
---------------------------------------------------------------------------
\draw(5.4,0)--(6.6,0)--(8.1,-1.5);
\draw (6.77,-0.17)--(5.67,-1.67);
\draw  [dash pattern=on 1.5pt off 1pt](5.16,-4.08)--(6.6,0);
\draw  [dash pattern=on 1.5pt off 1pt](5.5,-1.5)--(6.6,0);

\draw  [dash pattern=on 1.5pt off 1pt](5.5,-1.5)--(5.16,-4.08);

\draw  [dash pattern=on 1.5pt off 1pt](5.5,-1.5)--(5.7,0);

\draw  [dash pattern=on 1.5pt off 1pt](4.96,-4.08)--(5.5,-1.5)--(5.81,0);

\draw (4.96,-4.08)--(6.4,0);

\draw [dash pattern=on 1.5pt off 1pt] (5.67,-1.67)--(6.24,0.01);

\node at(6.76,0.06) {\relsize{-10}$A_0$};
\node at (6.26,0.13){\relsize{-10}$C$};
\node at (5.7,0.13){\relsize{-10}$G$};
\node at (5.88,0.13){\relsize{-10}$H$};
\node at (6.5,0.13){\relsize{-10}$A'_i$};
\node at (7.12,-0.2){\relsize{-10}$A_{i+1}$};
\node at (5.77,-1.77){\relsize{-10}$B_{i+1}$};

\node at (4.82,-4.05){\relsize{-10}$B'_i$};

\node at (5.35,-4.1){\relsize{-10}$F_i$};
\node at (5.2,-1.45){\relsize{-10}$F_{i+1}$};
--------------------------------------------------------------------------------------------
\draw(8.7,0)--(9.9,0)--(11.57,-1.67);
\draw (10.07,-0.17)--(8.97,-1.67);
\draw  [dash pattern=on 1.5pt off 1pt](8.97,-1.67)--(11.57,-1.67);
\draw [dash pattern=on 1.5pt off 1pt] (8.97,-1.67)--(9.54,0.01);
\draw [dash pattern=on 1.5pt off 1pt] (8.97,-1.67)--(11,-1.1);

\node at(9.95,0.08) {\relsize{-10}$A_0$};
\node at (9.58,0.1){\relsize{-10}$C$};

\node at (10.42,-0.2){\relsize{-10}$A_{i+1}$};
\node at (8.66,-1.65){\relsize{-10}$B_{i+1}$};
\node at (11.73,-1.65){\relsize{-10}$H$};
\node at (11.12,-1.08){\relsize{-10}$I$};

\end{tikz1}

\caption{The case a)$_i$-d)$_{i+1}$. } \label{fig:M1}
\end{figure}

Now we consider the case $\|\overrightarrow{p_0p_i}\|>\|\overrightarrow{p_0p_{i+1}}\|$. Let $C=l_{\bl{u}_i}(B_{i+1})\cap(A_0A'_0\cup l_{\bl{u}_h}(A_0))$, then we claim $C$ cannot be contained in $A_0A'_0$. Let $D$ be the other endpoint of the chord of $P_{L_1}$ passing through $A_0$ and paralleling with $\overrightarrow{p_0p_i}$, then  by assumption $\|\overline{A_0D}\|=\mathcal{L}_{\overrightarrow{p_0p_i}}(A_0)>\mathcal{L}_{\overrightarrow{p_0p_i}}(C)=\|\overline{B_{i+1}C}\|$. Let $E=DB_{i+1}\cap A_0A'_0$, then by the same reason as in the case $b)_i-e)_{i+1}$ and moreover $\|\overrightarrow{p_0p_i}\|<\|\overrightarrow{p_0p_{i+1}}\|$ we can deduce $E$ lies stirctly on the lower right of $A_{i+1}$.  For sufficiently small $\epsilon>0$, we have $\mathcal{L}_{\overrightarrow{p_0p_{i+1}}}(A_{i+1}+\epsilon(-\bl{u}_h+\bl{u}_v))<\|\overline{A_{i+1}B_{i+1}}\|=\|\overrightarrow{p_0p_{i+1}}\|$, which contradicts with the condition that $\overrightarrow{p_0p_{i+1}}$ is of type d). 

By our argument above $C$ is contained in $l_{\bl{u}_h}(A_0)$. Similar as the case $\|\overrightarrow{p_0p_i}\|<\|\overrightarrow{p_0p_{i+1}}\|$ above, one proves easily $C$ is contained in $e_{\beta}(P_{L_1})$. Besides, since $\|\overrightarrow{p_0p_i}\|>\|\overrightarrow{p_0p_{i+1}}\|$, by (\ref{abcd}) we have $c\geq a, d>b$, hence $\|\overline{B_{i+1}C}\|<\|\overrightarrow{p_0p_i}\|$. As $\overrightarrow{p_0p_i}$ is of type a), either both or neither of $A_i$ and $A'_i$ lie between $A_0$ and $C$. For the former case, let $F_i=A_0+\overrightarrow{p_ip_0}, F_{i+1}=A_0+\overrightarrow{p_{i+1}p_0}$, $G=F_iF_{i+1}\cap l_{\bl{u}_h}(A_0)$ and  $H=B'_iF_{i+1}\cap l_{\bl{u}_h}(A_0)$. One computes easily $\overline{A_0G}=\tfrac{1}{d}\leq 1$, thus $H$ is contained in $e_{\beta}(P_{L_1})$. On the other hand, note that $A_{i+1}\neq A_0$ as $\overrightarrow{p_0p_{i+1}}$ is of type d), $B_{i+1}$ is easily seen to lie inside the interior of the triangle $A_{i+1}B'_iH$, a contradiction.

Now assuming neither $A_i$ or $A'_i$ is contained in $\overline{CA_0}$. Let $H=l_{\bl{u}_h}(B_{i+1})\cap A_0A_0'$, then in order that there exists a point $x$ of $P_{L_1}$ on the left of $B_{i+1}C$ such that $\mathcal{L_{\overrightarrow{p_0p_i}}}(x)=\|\overrightarrow{p_0p_i}\|>\|\overline{B_{i+1}C}\|$, there must be some edge of $P_{L_1}$ containing $B_{i+1}$ whose extension intersects with $A_0A'_0$ at some point between $A_{i+1}$ and $H$. However, this implies $\mathcal{L}_{\overrightarrow{p_0p_{i+1}}}(A_{i+1}+\epsilon(\bl{u}_h-\bl{u}_v))<\|\overrightarrow{p_0p_{i+1}}\|$ (if $A_{i+1}\neq A'_0$) or $\mathcal{L}_{\overrightarrow{p_0p_{i+1}}}(A_{i+1}-\epsilon\bl{u}_v)<\|\overrightarrow{p_0p_{i+1}}\|$ (if $A_{i+1}=A'_0$) for all sufficiently small $\epsilon>0$, again contradicting with the fact that $\overrightarrow{p_0p_{i+1}}$ is of type d).

2.2.2.4 The case b)$_i$-d)$_{i+1}$.

The proof will be given only in the case when $\|\overrightarrow{p_0p_i}\|<\|\overrightarrow{p_0p_{i+1}}\|$. The case when the inequality is reversed is omitted as its proof is similar.

As is shown in the figure below, there are three possibilities for the relative positions of $\overline{A'_iB'_i}$, $\overline{A_{i+1}B_{i+1}}$.  When $A'_i$ lies on the upper left of $A_{i+1}$ and the two chords have no common points, let $C=l_{\overrightarrow{p_0p_i}}\cap A_0A'_0$, then $\mathcal{L}_{\overrightarrow{p_0p_i}}(C)>\|\overline{A'_iB'_i}\|=\|\overrightarrow{p_0p_i}\|$, which contradicts with the condition that $\overrightarrow{p_0p_i}$ is of type b).

Now assume $A'_i$ lies on the upper left of $A_{i+1}$ such that $\overline{A'_iB'_i}$, $\overline{A_{i+1}B_{i+1}}$ has a common point, as is shown in the middle figure below. Let $C=B_{i+1}B'_i\cap A_0A'_0$, then similar as in the suitation we met in the cases b)$_i$-e)$_{i+1}$ and a)$_i$-d)$_{i+1}$ one verifies easily $A'_0$ must lie on the upper left of $C$.  Let $D=A_{i+1}+\overrightarrow{p_ip_0}$, $E=l_{\bl{u}_v}(A_{i+1})\cap B_{i+1}D$ and $F=l_{\bl{u}_v}(A_{i+1})\cap B_{i+1}B'_i$. Then since  $A_i'\neq A_0'$ as $\overrightarrow{p_0p_i}$ is of type b), one sees easily $\|\overline{A_{i+1}E}\|>\|\overline{A_{i+1}F}\|\geq\|\overline{A'_0G}\|\geq\|e_{\alpha}(P_{L_1})\|\geq 1$.  On the other hand, one computes $\|\overline{A_{i+1}E}\|=\tfrac{1}{a-c}\leq 1$, a contradiction.

Now we consider the case when $A_{i+1}$ lies on the upper left of $A'_i$. Let $C=l_{\overrightarrow{p_0p_{i+1}}}(B'_i)\cap A_0A'_0$, then one sees easily $A'_0$ lies on the upper left of $C$.  Let $D$ be the other endpoint of the chord of $P_{L_1}$ through $A'_0$ with direction $\overrightarrow{p_0p_{i+1}}$, then we have $\|\overline{A'_0D}\|>\|\overrightarrow{p_0p_{i+1}}\|$. Let $E=DB'_i\cap A_0A'_0$, then one sees easily $E$ on the lower right of $C$. Let $F=A'_i+\overrightarrow{p_{i+1}p_0}$, $G=l_{\bl{u}_v}(A'_i)\cap FB'_i, H=l_{\bl{u}_v}(A'_0)\cap DB'_i$, then $\|\overline{A'_iG}\|>\|\overline{A'_0H}\|\geq\|e_{\alpha}(P_{L_1})\|\geq1$. However, one computes $\|\overline{A'_iG}\|=\tfrac{1}{a-c}\leq 1$, again a contradiction.

\begin{figure}[H]
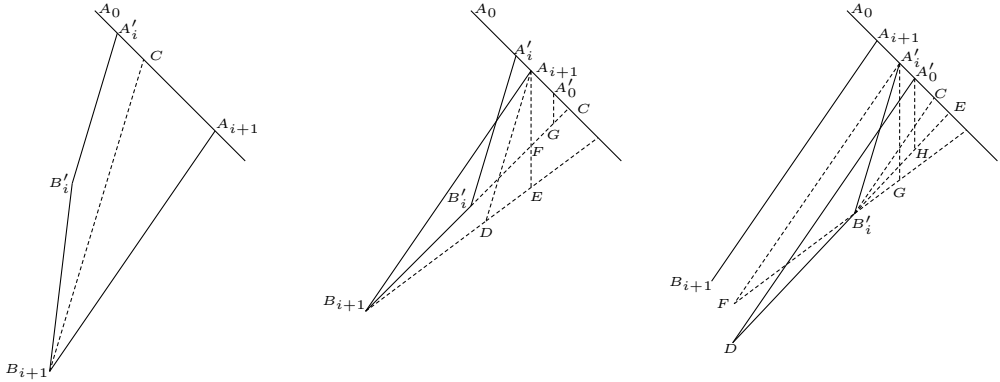


\begin{tikz1}

\draw (-5,0)--(-3,-2);
\draw(-4.7, -0.3)--(-5.3, -2.3);
\draw(-3.4, -1.6)--(-5.6,-4.8);
\draw (-5.6,-4.8)--(-5.3, -2.3);
\draw[dash pattern=on 1.5pt off 1pt] (-5.6,-4.8)--(-4.35,-0.65);

\node at(-4.8,0) {\relsize{-20}$A_0$};
\node at(-3.1,-1.55) {\relsize{-20}$A_{i+1}$};
\node at(-4.18,-0.6) {\relsize{-20}$C$};
\node at(-5.9,-4.8) {\relsize{-20}$B_{i+1}$};
\node at(-4.55,-0.25) {\relsize{-20}$A'_i$};
\node at(-5.47,-2.3) {\relsize{-20}$B'_i$};
-----------------------------------------------------------------------------

\draw (0,0)--(2,-2);
\draw(0.6, -0.6)--(0, -2.6);
\draw(0.8, -0.8)--(-1.4,-4);
\draw[dash pattern=on 1.5pt off 1pt](0.8,-0.8)--(0.2,-2.8);
\draw(-1.4,-4)--(0,-2.6);
\draw[dash pattern=on 1.5pt off 1pt](-1.4,-4)--(0.2,-2.8);
\draw[dash pattern=on 1.5pt off 1pt](0,-2.6)--(1.3,-1.3);
\draw[dash pattern=on 1.5pt off 1pt](0.2,-2.8)--(1.69,-1.69);
\draw [dash pattern=on 1.5pt off 1pt](0.8,-0.8)--(0.8,-2.35);
\draw[dash pattern=on 1.5pt off 1pt] (1.1, -1.1)--(1.1,-1.53);

\node at(0.2,0) {\relsize{-20}$A_0$};
\node at(0.2,-2.95) {\relsize{-20}$D$};
\node at(1.15,-0.8) {\relsize{-20}$A_{i+1}$};
\node at(1.26,-1.03) {\relsize{-20}$A'_0$};
\node at(0.87,-2.47) {\relsize{-20}$E$};
\node at(0.89,-1.87) {\relsize{-20}$F$};
\node at(1.1,-1.65) {\relsize{-20}$G$};
\node at(1.5,-1.3) {\relsize{-20}$C$};
\node at(0.7,-0.53) {\relsize{-20}$A'_i$};

\node at(-1.68,-3.9) {\relsize{-20}$B_{i+1}$};
\node at(-0.18,-2.5) {\relsize{-20}$B'_i$};
----------------------------------------------------
\draw (5,0)--(7,-2);
\draw(5.7, -0.7)--(5.1, -2.7);

\draw(5.4, -0.4)--(3.2,-3.6);
\draw(5.9, -0.9)--(3.48,-4.42);
\draw (3.48,-4.42)--(5.1, -2.7);
\draw [dash pattern=on 1.5pt off 1pt](5.1, -2.7)--(6.16,-1.16);
\draw [dash pattern=on 1.5pt off 1pt](3.5, -3.9)--(5.1, -2.7);
\draw[dash pattern=on 1.5pt off 1pt] (5.7, -0.7)--(3.5, -3.9);
\draw[dash pattern=on 1.5pt off 1pt] (5.1, -2.7)--(6.59,-1.59);
\draw[dash pattern=on 1.5pt off 1pt] (5.1, -2.7)--(6.36,-1.36);

\draw[dash pattern=on 1.5pt off 1pt] (5.7, -0.7)--(5.7,-2.25);
\draw[dash pattern=on 1.5pt off 1pt] (5.9, -0.9)--(5.9,-1.86);

\node at(5.2,0) {\relsize{-20}$A_0$};
\node at(6.5,-1.28) {\relsize{-20}$E$};
\node at(6.25,-1.1) {\relsize{-20}$C$};
\node at(5.85,-0.63) {\relsize{-20}$A'_i$};
\node at(6.06,-0.82) {\relsize{-20}$A'_0$};
\node at(5.69,-0.35) {\relsize{-20}$A_{i+1}$};
\node at(3.45,-4.5) {\relsize{-20}$D$};
\node at(6,-1.9) {\relsize{-20}$H$};
\node at(5.7,-2.42) {\relsize{-20}$G$};
\node at(3.35,-3.9) {\relsize{-20}$F$};
\node at(5.2,-2.85) {\relsize{-20}$B'_i$};

\node at(2.93,-3.65) {\relsize{-20}$B_{i+1}$};
\end{tikz1}

\caption{The case b)$_i$-d)$_{i+1}$, $\|\protect\overrightarrow{p_0p_i}\|<\|\protect\overrightarrow{p_0p_{i+1}}\|$.  } \label{fig:M1}
\end{figure}

\end{proof}

\bigskip

\emph{Step 3. Proof of (I2)}\hfill

The notations for some points and vectors in the following are different from those used in Step 2. 

\begin{proof}
Next we prove (I2), i.e. the surjectivity of $\phi_{L_1, L_2(D_{\sigma})}$ implies that of $\phi_{L_1, L_2}$ for $L_1$ with $L_1.D_{\sigma}=1$. We will regard $P_{L_2}$ as the polygon obtained from $P_{L_2(D_{\sigma})}$ by cutting off a stripe from the upper right. By assumption, we have
\begin{equation}\label{cr3}P_{L_2(D_{\sigma})}=\bigcup_{\chi\in P_{L^{-1}_1\otimes L_2(D_{\sigma})}\cap M}\chi+P_{L_1}.
\end{equation}

Let \[\Xi=\{\chi\in P_{L^{-1}_1\otimes L_2(D_{\sigma})}\cap M\:|\:\chi+P_{L_1}\:\text{is not contained in}\:P_{L_2}\},\]
then to prove the surjectivity of (\ref{f12})  we only need to show for each $\chi\in\Xi$
\begin{equation}\label{2inc}(\chi+P_{L_1})\cap P_{L_2}\subseteq\bigcup_{\chi'\in (P_{L_1^{-1}\otimes L_2}\cap M)}\chi'+P_{L_1},
\end{equation}
i.e. $\chi+P_{L_1}$ can be covered by lattice translations of $P_{L_1}$ contained in $P_{L_2}$.

Firstly we note if $\chi\in \Xi$, then either $\chi\in e_{\alpha+\beta}(P_{L_1^{-1}\otimes L_2(D_{\sigma})})$ or $\chi$ is the common vertex of $e_{\alpha}(P_{L_1^{-1}\otimes L_2(D_{\sigma}) })$ and $e_{\beta}(P_{L_1^{-1}\otimes L_2(D_{\sigma})})$ when $\|e_{\alpha+\beta}(P_{L^{-1}_1\otimes L_2(D_{\sigma})})\|=0$. Besides, since
\[L_1^{-1}\otimes L_2(D_{\sigma}).D_{\alpha}=L_2(D_{\sigma}).D_{\alpha}-L_1.D_{\alpha}=L_2.D_{\alpha}+1-L_1.D_{\alpha}\geq1,\]
 we have $e_{\alpha}(P_{L_1^{-1}\otimes L_2(D_{\sigma}) })\geq1$ and similarly $e_{\beta}(P_{L_1^{-1}\otimes L_2(D_{\sigma})})\geq1$. Therefore, $-\bl{u}_v+\chi$ and $-\bl{u}_h+\chi$ are both contained in $P_{L^{-1}_1\otimes L_2}$.

If $-\alpha, -\beta\in\Sigma_X(1)$ then one can prove directly 
\[(\chi+P_{L_1})\cap P_{L_2}\subseteq (-\bl{u}_h+\chi+P_{L_1})\:{\textstyle \bigcup}\:(-\bl{u}_v+\chi+P_{L_1}).
\]
The proof for the cases when $-\alpha\notin\Sigma_X(1)$ or $-\beta\notin\Sigma_X(1)$ are the same and next we will assume $-\beta\notin\Sigma_X(1)$.

Under the assumption above, by Lemma \ref{1q1c} $\chi+P_{L_1}$ has a unique lowest point and we will denote it by $C$. Let $\bl{u}_{-1}, \bl{u}_1$  (in counter-clockwise order) be the two primitive edge vectors of $\chi+P_{L_1}$ with their common initial point at $C$. Let $C_h=-\bl{u}_h+C$, $C_v=-\bl{u}_v+C$ and $D=l_{\bl{u}_1}(C_h)\cap l_{\bl{u}_{-1}}(C_v)$, then one checks easily $D$ is also a lattice point. Indeed, let 
\begin{equation}\label{abcd}\bl{u}_{-1}=a\bl{u}_h+b\bl{u}_v,\;\;\; \bl{u}_1=c\bl{u}_h+d\bl{u}_v, 
\end{equation}
one gets
\begin{equation}\label{dc}\overrightarrow{DC}=a(c+d)\bl{u}_h+d(a+b)\bl{u}_v=a\bl{u}_1+d\bl{u}_{-1}.
\end{equation}

We claim $D$ and the translation $\overrightarrow{CD}+\chi+P_{L_1}$ are contained in $P_{L_2}$. Let $P$ be the convex hull of $-\bl{u}_h+\chi+P_{L_1}, -\bl{u}_v+\chi+P_{L_1}$ together with $D$, then by our argument above $P$ is a lattice polygon and its normal fan is just $\Sigma_X$.  Therefore, to prove our claim it suffices to show $D$ is contained in $P_{L_2}$. 
Let $\tilde{C}$ be the lowest vertex of $P_{L_2}$, $\tilde{C}_{-1}=l_{\bl{u}_1}(C_h)\cap\bigcup_{0\leq i\leq k}e_{\alpha_i}(P_{L_2})$, $\tilde{C}_1=l_{\bl{u}_{-1}}(C_v)\cap\bigcup_{0\leq j\leq l}e_{\beta_j}(P_{L_2})$. As is shown in the figure on the left below,  $D$ lies inside the convex hull of $C_h, C_v, \tilde{C}_{-1}$ and $\tilde{C}_1$, hence contained in $P_{L_2}$.

 Now by  (\ref{dc}), as is shown by the figure on the right below, the lattice points corresponding to $-\bl{u}_{-1}+\chi$ and $-\bl{u}_1+\chi$ are contained in the convex hull of the lattice points $-\bl{u}_h+\chi, -\bl{u}_v+\chi, \overrightarrow{CD}+\chi$ as there are no other lattices points in the triangle with vertices $\chi, -\bl{u}_h+\chi$ and $-\bl{u}_v+\chi$. Therefore we have 
$-\bl{u}_{-1}+\chi, -\bl{u}_1+\chi\in P_{L_1^{-1}\otimes L_2}$ since the same is true for $-\bl{u}_h+\chi, -\bl{u}_v+\chi$ and $\overrightarrow{CD}+\chi$.

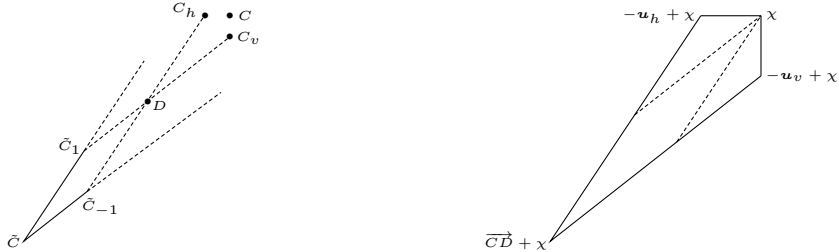
 \begin{tikz1}

\draw  (0.84,0.66)--(0,0)--(0.8,1.2); 

\draw  [dash pattern=on 1.5pt off 1pt](0.8,1.2)--(1.6,2.4); 
\draw  [dash pattern=on 1.5pt off 1pt](0.84,0.66)--(2.62,1.98); 

\draw  [dash pattern=on 1.5pt off 1pt](0.84,0.66)--(2.4,3); 
\draw  [dash pattern=on 1.5pt off 1pt](0.8,1.2)--(2.74,2.72);

\node at (-0.18,0)  (C) {\relsize{-20} $\tilde{C}$}; 
\node at (0.57,1.27)  (C) {\relsize{-20} $\tilde{C}_1$}; 
\node at (0.98,0.45)  (C) {\relsize{-20} $\tilde{C}_{-1}$}; 

\node at (2.1,3.08)  (C) {\relsize{-20} $C_h$}; \node at (2.41,3)  (C) {\relsize{2}$\bl{\cdot}$};

\node at (2.95,2.7)  (C) {\relsize{-20} $C_v$}; \node at (2.74,2.72)  (C) {\relsize{2}$\bl{\cdot}$};

\node at (2.9,3)  (C) {\relsize{-20} $C$}; \node at (2.74,3)  (C) {\relsize{2}$\bl{\cdot}$};

\node at (1.76,1.8)  (C) {\relsize{-20} $D$}; \node at (1.65,1.85)  (C) {\relsize{2}$\bl{\cdot}$};
-----------------------------------------------------------------------------------------------------------------------------------------------
\draw  (9,3)--(7,0)--(9.8,2.2); 
\draw  [dash pattern=on 1.5pt off 1pt](8.12,1.68)--(9.8,3); 
\draw  [dash pattern=on 1.5pt off 1pt](8.68,1.32)--(9.8,3); 
\draw (9,3)--(9.8,3);
\draw (9.8,2.2)--(9.8,3);

\node at (9.9,3)  (C) {\relsize{-20} $\chi$}; 
\node at (8.4,3)  (C) {\relsize{-20} $-\bl{u}_h+\chi$}; 
\node at (10.3,2.2)  (C) {\relsize{-20} $-\bl{u}_v+\chi$}; 
\node at (6.5,0)  (C) {\relsize{-20} $\overrightarrow{CD}+\chi$}; 
\end{tikz1}

Next we will show 
\begin{equation}\label{2c}(\chi+P_{L_1})\cap P_{L_2}\subseteq         (-\bl{u}_h+\chi+P_{L_1})\:{\textstyle \bigcup}\:(-\bl{u}_v+\chi+P_{L_1})\:{\textstyle \bigcup}\:(-\bl{u}_{-1}+\chi+P_{L_1})\:{\textstyle \bigcup}\:(-\bl{u}_1+\chi+P_{L_1}).
\end{equation}

For saving notations, we will drop $\chi$ from $\chi+P_{L_1}$ from now on. Next we label some contact points of $P_{L_1}$ that will be used in the sequel. As is shown in the figure on the right below, let $A, B$ be the endpoints of $e_{\alpha+\beta}(P_{L_1})$, then one of the contact point with respect to $-\bl{u}_h$ is $A_{-h}=-\bl{u}_{h}+A$.  The other contact point will be denoted by $E$, as is shown  in the figure on the left below.  
Similarly, we let $B_{-v}=-\bl{u}_v+B$ and $F$ be the contact points of $P_{L_1}$ in the direction of $-\bl{u}_v$, $A_{-v}=-\bl{u}_v+A$ and $G$ be the contact points of $-\bl{u}_h+P_{L_1}$ in the direction of $\bl{u}_h-\bl{u}_v$. In the figure on the left below, we let $C_{-1}=-\bl{u}_{-1}+C$ and $C_1=-\bl{u}_1+C$.

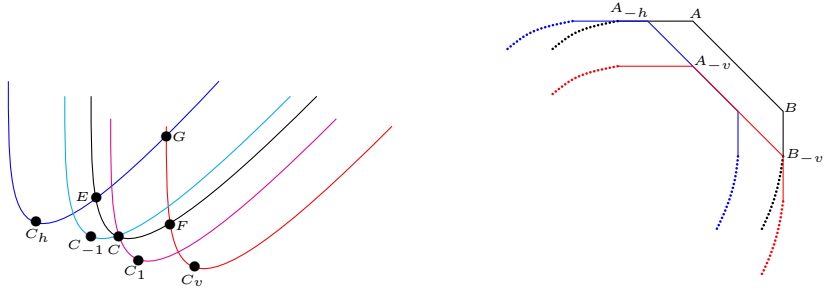
\begin{tikz1}

\draw (0,3) .. controls (0,1) and (0,0) .. (3,3);

\draw[blue] (-1.1,3.2) .. controls (-1.1,1.2) and (-1.1,0.2) .. (1.7,3.2);

\draw[red] (1,2.6) .. controls (1,0.6) and (1,-0.4) .. (4,2.6);

\draw[cyan] (-0.35,3) .. controls (-0.35,1) and (-0.35,0) .. (2.65,3);

\draw[magenta] (0.26,2.7) .. controls (0.26,0.7) and (0.26,-0.3) .. (3.26,2.7);

\node at (0.27,0.96)  (C) {\relsize{-20} $C$}; \node at (0.37,1.13)  (C) {$\bullet$};

\node at (-0.75,1.15)  (C) {\relsize{-20} $C_h$}; \node at (-0.73,1.32)  (C) {$\bullet$};

\node at (1.35,0.55)  (C) {\relsize{-20}$C_v$}; \node at (1.38,0.73)  (C) {$\bullet$}; 

\node at (0.55,0.65)  (C) {\relsize{-20} $C_1$}; \node at (0.63,0.81)  (C) {$\bullet$};

\node at (-0.1,0.98)  (C) {\relsize{-20} $C_{-1}$}; \node at (0,1.12)  (C) {$\bullet$}; 

\node at (-0.16,1.67)  (C) {\relsize{-20} $E$}; \node at (0.07,1.65)  (C) {$\bullet$};

\node at (1.17,1.28)  (F) {\relsize{-20} $F$}; \node at (1.05,1.28)  (C) {$\bullet$}; 

\node at (1.13,2.47)  (G) {\relsize{-20} $G$}; \node at (1,2.45)  (C) {$\bullet$}; 
=========================================================================================
\draw(7,4)--(8,4)--(9.2,2.8)--(9.2,2.2);
\draw [line width=1pt, line cap=round, dash pattern=on 0pt off 1.5\pgflinewidth](7,4) .. controls (6.47,3.9) and (6.4,3.85) .. (6.1,3.6); 
\draw [line width=1pt, line cap=round, dash pattern=on 0pt off 1.5\pgflinewidth](9.2,2.2) .. controls (9.18,2.1) and (9.16,1.7) .. (8.9,1.2); 
---------------------------------------------------------------------------------------------------
\draw[color=blue](6.4,4)--(7.4,4)--(8.6,2.8)--(8.6,2.2);
\draw[line width=1pt, line cap=round, dash pattern=on 0pt off 1.5\pgflinewidth, color=blue] (6.4,4) .. controls (5.87,3.9) and (5.8,3.85) .. (5.5,3.6); 
\draw [line width=1pt, line cap=round, dash pattern=on 0pt off 1.5\pgflinewidth, color=blue](8.6,2.2) .. controls (8.58,2) and (8.56,1.7) .. (8.3,1.2); 
-------------------------------------------------------------------------------------------------
\draw[color=red](7,3.4)--(8,3.4)--(9.2,2.2)--(9.2,1.6);
\draw [line width=1pt, line cap=round, dash pattern=on 0pt off 1.5\pgflinewidth, color=red](7,3.4) .. controls (6.47,3.3) and (6.4,3.25) .. (6.1,3); 
\draw [line width=1pt, line cap=round, dash pattern=on 0pt off 1.5\pgflinewidth, color=red](9.2,1.6) .. controls (9.18,1.5) and (9.16,1.1) .. (8.9,0.6);

\node at (8,4.1)  (C) {\relsize{-20} $A$}; 

\node at (9.26,2.85)  (C) {\relsize{-20} $B$}; 

\node at (7.12,4.13)  (C) {\relsize{-20} $A_{-h}$}; 


\node at (8.23,3.43)  (C) {\relsize{-20} $A_{-v}$}; 

\node at (9.46,2.2)  (C) {\relsize{-20} $B_{-v}$};

\end{tikz1}
Note that it is by no means self-evident that $E$ (or $F$) should reside in a position lower than $G$. However, we will prove this is indeed the case.
Let $\mathscr{P}_{1, -1}$ be the angular region bounded by $r_{\bl{u}_{-1}}(C)$ and $r_{\bl{u}_1}(C)$. By a short computation, one finds easily the chord of $\mathscr{P}_{1, -1}$ with direction $\bl{u}_h-\bl{u}_v$ and length $\|\bl{u}_h-\bl{u}_v\|$ is lying over the chord with direction $\bl{u}_h$ (resp. $\bl{u}_v$) and length $\|\bl{u}_h\|$ (resp. $\|\bl{u}_v\|$), then by Lemma \ref{ip} $G$ lies over both $E$ and $F$.

Now one sees easily the complement of $(-\bl{u}_h+P_{L_1})\cup(-\bl{u}_v+P_{L_1})$ in $P_{L_2}$ is  contained in the polygon represented by the curvilinear quadrilateral $GECF$, hence we only need to prove it is contained in $(-\bl{u}_{-1}+P_{L_1})\cup(-\bl{u}_1+P_{L_1})$. As can be seen from the figure on the left above, it suffices to show

1) $-\bl{u}_1+P_{L_1}$ and $-\bl{u}_h+P_{L_1}$ (resp. $-\bl{u}_v+P_{L_1}$) has a common point lying between $E$ and $G$ (resp. $C_v$ and $F$).

2) $-\bl{u}_{-1}+P_{L_1}$ and $-\bl{u}_h+P_{L_1}$ (resp. $-\bl{u}_v+P_{L_1}$) has a common point lying between $C_h$ and $E$ (resp. $F$ and $G$);


\indent The proof of the two claims are similar and we will only provide the first one. 

Let $\mathscr{P}_{1, -1}$ be the cone bounded by $r_{\bl{u}_{-1}}(C)$ and  $r_{\bl{u}_1}(C)$. One sees easily the chord of $\mathscr{P}_{1, -1}$ with direction
$\bl{u}_1-\bl{u}_h$ and length $\|\bl{u}_1-\bl{u}_h\|$ (resp. direction
$\bl{u}_1-\bl{u}_v$ and length $\|\bl{u}_1-\bl{u}_v\|$) shares an endpoint with but lies on the upper side (resp. lower side) of the chord with direction $\bl{u}_h$ and length $\|\bl{u}_h\|$ (resp. direction $\bl{u}_v$ and length $\|\bl{u}_v\|$). Furthermore, by using (\ref{abcd}) and the fact that $ad-bc=1$
one checks easily the chord of $\mathscr{P}_{1, -1}$ with direction $\bl{u}_v-\bl{u}_h$ and length $\|\bl{u}_v-\bl{u}_h\|$  lies above the chord with direction $\bl{u}_1-\bl{u}_h$ and length $\|\bl{u}_1-\bl{u}_h\|$. 
Then 1) follows from Lemma \ref{ip}, hence the proof of (I2) and Theorem \ref{sfhn} is complete.
 \end{proof}
 
\bigskip
  
Following notations in Step 2.2, for any $\alpha_i, 1\leq i\leq k$ and $\beta_j, 1\leq j\leq l$ there are always two contact points of $P_{L_1}$ with respect to the corresponding primitive edge vector. In general these vectors are not necessarily to be of type c). 
However, by our proof for the case when $p_0\in\partial P_{L_1^{-1}\otimes L_2}$ we can deduce this is indeed the case in certain circumstance. 

\begin{pro}

For $\bl{u}_{\alpha_i}$ (resp. $\bl{u}_{\beta_j}$) if there exists $\bl{u}_{\beta_j}$ (resp. $\bl{u}_{\alpha_i}$) such that the following condition are satisfied,  then it
is of type c).

i) the slope of $\bl{u}_{\alpha_i}$ is larger than that of $\bl{u}_{\beta_j}$;

ii)$\|\bl{u}_{\alpha_i}\|<\|\bl{u}_{\beta_j}\|$ (resp. $\|\bl{u}_{\beta_j}\|<\|\bl{u}_{\alpha_i}\|$);

iii) $\bl{u}_{\alpha_i}$ and $\bl{u}_{\beta_j}$ forms a basis of $M$.

\end{pro}

\begin{proof}
We only prove the conclusion for $\bl{u}_{\beta_j}$ and the other case can be proved similarly. Following the notations of in the proof of Claim \ref{cl2.0} for $p_0\in\partial P_{L_1^{-1}\otimes L_2}$, under the condition of the proposition, there exists only one lattice point in $P_{L_1^{-1}\otimes L_2}\cap R$ hence $\bl{u}_{\beta_j}=\overrightarrow{p_0p_1}$ must be of type c).
\end{proof}

\begin{rem} As a byproduct of our proof in Step 2.2, the conclusions (\ref{k1}) and (\ref{l1}) can be reformulated as follows.
Let $P, Q$ be lattice polygons such that $Q$ is smooth and the normal fan of $Q$ is a refinement of that of $P$, then for any lattice point $p$ in the interior $P$, there always exist lattice points $q_1, q_2\in P$ (which might be the same) such that $\overrightarrow{pq_1}$ and $\overrightarrow{pq_2}$ parallel with two edges of $Q$. It might be interesting to know whether we still have similar results in three dimension.
\end{rem}

\bibliographystyle{amsplain}
\bibliography{xhBibTex}
\end{document}